\documentclass[11pt]{article}
\usepackage[dvips]{graphicx}
\usepackage[dvips,letterpaper,margin=1in]{geometry}
\usepackage{color}
\usepackage{amsmath}
\usepackage{amsfonts}
\usepackage{amssymb}
\usepackage{amsthm}
\usepackage{newlfont}
\usepackage{epsfig}
\usepackage{multirow}
\usepackage{setspace}
\usepackage{hyperref}
\usepackage{cite}
\usepackage{bbm}
\usepackage{wrapfig}

\begin{document}

\newtheorem{assumption}{Assumption}[section]
\newtheorem{definition}{Definition}[section]
\newtheorem{lemma}{Lemma}[section]
\newtheorem{proposition}{Proposition}[section]
\newtheorem{theorem}{Theorem}[section]
\newtheorem{corollary}{Corollary}[section]
\newtheorem{remark}{Remark}[section]

\small

\title{A Computational Approach to Steady State Correspondence of Regular and Generalized Mass Action Systems}
\author{Matthew D. Johnston\\ \\Department of Mathematics\\University of Wisconsin-Madison\\480 Lincoln Dr., Madison, WI 53706\\email: mjohnston3@wisc.edu}
\date{}
\maketitle


\tableofcontents

\begin{abstract}
\small
It has been recently observed that the dynamical properties of mass action systems arising from many models of biochemical reaction networks can be derived by considering the corresponding properties of a related generalized mass action system. The correspondence process known as network translation in particular has been shown to be useful in characterizing a system's steady states. In this paper, we further develop the theory of network translation with particular focus on a subclass of translations known as improper translations. For these translations, we derive conditions on the network topology of the translated network which are sufficient to guarantee the original and translated systems share the same steady states. We then present a mixed-integer linear programming (MILP) algorithm capable of determining whether a mass action system can be corresponded to a generalized system through the process of network translation.
\end{abstract}



\bigskip

\section{Introduction}
\label{introduction}

Many biochemical and industrial processes can be represented graphically as networks of simultaneously occurring chemical reactions. Under simplifying assumptions such spatially homogeneity and mass action kinetics, the dynamical behavior of these \emph{chemical reaction networks} can be modeled mathematically by systems of autonomous polynomial ordinary differential equations known as \emph{mass action systems}. 


Motivated by the growth of systems biology, there has been significant recent interest in characterizing the long-term and steady state properties of such systems. A recent addition to this field has been the study of \emph{generalized chemical reaction networks}, which was introduced by M\"{u}ller and Regensburger in \cite{M-R}. A generalized chemical reaction network is given by a chemical reaction network together with an additional set of vertices which are in one-to-one correspondence with the vertices of the original network. The dynamics of these generalized networks are then given by a \emph{generalized mass action system}, where the first set of vertices controls the stoichiometry of the system (i.e. the reaction vectors), and the second set controls the kinetic rates (i.e. the reaction monomials). For example, consider the generalized network
\begin{equation}
\label{9919}
X_1 \; \; \; \cdots \; \; \; \; \; \; \; X_1 \; \; \; \mathop{\stackrel{k_1}{\begin{array}{c} \vspace{-0.25cm} \longleftarrow \vspace{-0.3cm} \\ \longrightarrow \end{array}}}_{k_2} \; \; \; X_2 \; \; \; \; \; \; \; \cdots \; \; \; 2X_2
\end{equation}
where the dotted lines denote the correspondence between the stoichiometric vertices and the kinetic vertices. The dynamic formulation of the corresponding generalized mass action system is the same as for a regular one except that we substitute the monomial $x_2^2$ corresponding to $2X_2$ in the place of the monomial $x_2$ corresponding to $X_2$. Although the theoretical study of generalized systems is in its early stages, several substantial results are known, including results sufficient to guarantee the existence of ``complex-balanced-like'' steady states, and results guaranteeing the uniqueness of such states within stoichiometric compatibility classes \cite{M-R,M-F-R-C-S-D}.

It was noted by Johnston in \cite{J1} that dynamical and steady state properties of classical mass action systems can often be determined by first making a suitable correspondence with a generalized mass action system. For example, consider the regular network
\begin{equation}
\label{9929}
X_1 \; \stackrel{k_1}{\longrightarrow} \; X_2, \hspace{1.5cm} 2X_2 \; \stackrel{k_2}{\longrightarrow} \; X_1 + X_2.
\end{equation}
Despite the difference in appearance and network structure between (\ref{9919}) and (\ref{9929}), it can be easily verified that they share the same governing set of differential equations. Johnston introduced a correspondence process called \emph{network translation} and was able to identify two subcategories: proper translations and improper translations. A translation is said to be \emph{proper} if there is a one-to-one correspondence between the source vertices of the original network and those of the translated network; otherwise, it is said to be \emph{improper}. For proper translations, the original and generalized systems are known to be \emph{dynamically equivalent} (Lemma 2, \cite{J1}) while for improper translations supplemental conditions are known which allow rate constants to be selected so that the original and generalized systems share the same steady states (Lemma 4, \cite{J1}). Johnston also gave conditions which are sufficient to guarantee the existence of \emph{toric steady states} as introduced by P\'{e}rez Mill\'{a}n \emph{et al.} in \cite{M-D-S-C} (Theorem 5, \cite{J1}). The method of network translation has since been applied to characterize the steady states of processive multisite phosphorylation networks by Conradi and Shiu in \cite{C-S}.

Two important questions were left open in \cite{J1} which we address in the current work:
\begin{enumerate}
\item[(Q1)] Given an improper translation, are there sufficient conditions on the structure of the \emph{translated reaction graph alone} which guarantee steady state equivalence of the original and translated systems?
\end{enumerate}
While sufficient conditions were given in \cite{J1} for guaranteeing steady state equivalence of the two systems, the conditions depended upon an algebraic combination of rate constants which may be difficult to compute in practice. In Section \ref{translationsection}, we improve upon this result by presenting conditions on the \emph{translated reaction graph alone} which are sufficient to guarantee such a correspondence can be made (Theorem \ref{maintheorem}). This is in keeping with the general flavor of so-called \emph{chemical reaction network theory} (CRNT) which has placed considerable emphasis on dynamical results which follow from properties of the underlying network structure.

After answering (Q1), we consider the following more fundamental question:
\begin{enumerate}
\item[(Q2)] Given a mass action system, can we algorithmically determine the structure underlying a generalized system which is either dynamically or steady state equivalent to the original system? 
\end{enumerate}
It was noted in \cite{J1} that, in practice, we do not have the structure of the translated network give to us; rather we must \emph{find} it. Even for networks of only moderate size, computing this structure by hand alone can be extremely difficult. While an algorithm for constructing translations was presented in \cite{J1}, it was not directly amenable to computational implementation as it required a full enumeration of all possible cyclic combinations of reactions on the network's stoichiometric generators. There was also no guarantee that the translation would satisfiy desirable network properties such as being weakly reversible or having a low deficiency.

In Section \ref{milpsection}, we recast this fundamental question as a mixed-integer linear programming (MILP) problem. This framework has been previously used within CRNT to determine dynamically equivalent and linearly conjugate network structures in the papers of Szederk\'{e}nyi and various collaborators in \cite{Sz2,Sz-H,Sz-H-P,Sz-H-T,J-S4,J-S5,J-S6,R-S-H}. The algorithm we present here is capable of determining the structure of the translated chemical reaction network, ensuring steady state equivalence may be made in accordance with Theorem \ref{maintheorem}, and also guaranteeing weak reversibility and a minimal deficiency is attained according to the results of \cite{J-S4,J-S6}. In Section \ref{applications}, we apply the computational algorithm to a pair of models drawn from the mathematical biochemistry literature to determine a generalized mass action system with the same steady states \cite{Sh-F,Dasgupta,Karp}. 

\section{Background}
\label{backgroundsection}

In this section, we present the required background information on CRNT in both the classical and generalized setting.

\subsection{Chemical Reaction Networks}

The central object of study in this paper is the following.

\begin{definition}
\label{crn}
A \textbf{chemical reaction network} is a triple $\mathcal{N} = (\mathcal{S},\mathcal{C},\mathcal{R})$ where:
\begin{enumerate}
\item
The \textbf{species set} $\mathcal{S} = \left\{ X_1, \ldots, X_n \right\}$ consists of the individual (chemical) species $X_i$ capable of undergoing chemical change.
\item
The \textbf{complex set} $\mathcal{C} = \left\{ C_1, \ldots, C_m \right\}$ consists of linear combinations of the species, i.e. terms of the form $C_j = \sum_{i=1}^n y_{ji} X_i$, $j=1,\ldots, m$. The values $y_{ji} \in \mathbb{Z}_{\geq 0}$ are called stoichiometric coefficients and each complex $C_j$ is associated with a stoichiometric vector $y_j = (y_{j1}, y_{j2}, \ldots, y_{jn})$. It is assumed that the complexes are stoichiometrically distinct, i.e. $y_i \not= y_j$ for $i \not= j$.
\item
The \textbf{reaction set} $\mathcal{R} \subseteq \mathcal{C} \times \mathcal{C}$ consists of ordered pairs $(C_i,C_j)$ where $C_i, C_j \in \mathcal{C}$. It is also common to represent reactions in the form $C_i \to C_j$.
\end{enumerate}
\end{definition}

\begin{remark}
It is typical in CRNT to assume that $(i)$ every species appears in at least one complex, $(ii)$ every complex appears in at least one reaction (as either a reactant or product), and (iii) there are no self-reactions (i.e. reactions of the form $C_i \to C_i$). To accommodate the computational processes used in Section \ref{milpsection}, it will be occasionally necessary to violate condition (ii). These exceptions will be noted in the text.
\end{remark}

\begin{remark}
It will be convenient to allow the complex set $\mathcal{C}$ to correspond to the underlying index set, i.e. we will let $\mathcal{C} = \left\{ 1, \ldots, m \right\}$ and allow $i \in \mathcal{C}$ to stand in for $C_i \in \mathcal{C}$. We will also allow the ordered index pair $(i,j) \in \mathcal{R}$ to represent the reaction $C_i \to C_j$. 
\end{remark}




It is natural to interpret chemical reaction networks as directed graphs $G(V,E)$ where the vertex set is given by the complexes (i.e. $V = \mathcal{C}$) and the edge set is given by the reactions (i.e. $E = \mathcal{R}$). Two complexes $C_i$ and $C_j$ are said to be \emph{connected} if there is a sequence of complexes such that $C_i = C_{\nu(1)} \leftrightarrow C_{\nu(2)} \leftrightarrow \cdots \leftrightarrow C_{\nu(l)} = C_j$ where $C \leftrightarrow C'$ if $C \to C'$ or $C' \to C$. If there is such a chain where all the reactions are of the form $C \to C'$, we say there is a \emph{path} from $C_i$ to $C_j$. The maximal sets of connected complexes are called \emph{linkage classes} and are denoted $\mathcal{L} = ( L_1, \ldots, L_\ell)$ where $\ell = | \mathcal{L} |$. Two complexes $C_i$ and $C_j$ are said to be \emph{strongly connected} if, given a path from $C_i$ to $C_j$, there is a path from $C_j$ to $C_i$. The maximal sets of strongly connected complexes are called a \emph{strong linkage classes}. A network is said to be \emph{weakly reversible} if the linkage classes and strong linkage classes coincide.

To each reaction $C_i \to C_j$ we associate the \emph{reaction vector} $y_j - y_i \in \mathbb{Z}^n$ which keeps track of the change in the number of each species as a result of the reaction. The span of the reaction vectors is called the \emph{stoichiometric subspace} and is denoted $S = \mbox{span} \left\{ y_j - y_i \; | \; (i,j) \in \mathcal{R} \right\}$. The dimension of the stoichiometric subspace is denoted $s = \mbox{dim}(S)$.

A network parameter which has been particularly well studied in the literature is the \emph{deficiency} \cite{F1,H,F2,Fe2,Fe4}.

\begin{definition}
\label{deficiency}
The \textbf{deficiency} of a chemical reaction network $\mathcal{N} = (\mathcal{S},\mathcal{C},\mathcal{R})$ is given by
\[\delta = m - \ell - s\]
where $m$ is the number of stoichiometrically distinct complexes (i.e. $m = |\mathcal{C}|)$, $\ell$ is the number of linkage classes (i.e. $\ell = | \mathcal{L} |$), and $s$ is the dimension of the stoichiometric subspace (i.e. $s = \mbox{dim}(S)$).
\end{definition}

\subsection{Reaction-Weighted Networks and Mass Action Systems}
\label{reactionsection}

A common kinetic assumption for chemical reaction networks is \emph{mass action kinetics}, which states that the rate of a reaction is proportional to the product of the concentrations of the reacting species. For instance, if a reaction $C_i \to C_j$ has the form $X_1 + X_2 \to \cdots$, then the associated rate function would be $[\mbox{rate}] = k(i,j) [X_1][X_2]$ where $k(i,j) > 0$ is the \emph{rate constant} (i.e. proportionality constant) of the reaction. Other kinetic assumptions are also frequently used, especially in the mathematical biochemistry literature, including Michaelis-Menten kinetics \cite{M-M} and Hill kinetics \cite{Hi}.

It is therefore natural to associate to every reaction $(i,j) \in \mathcal{R}$ a \emph{reaction-weight} $k(i,j) > 0$. We formally define the following.
\begin{definition}
Suppose $\mathcal{N} = (\mathcal{S},\mathcal{C},\mathcal{R})$ is a chemical reaction network. We will say that $\mathcal{K} = \{ k(i,j) \; | \; i,j \in \mathcal{C} \}$ is a \textbf{reaction-weight set} if $k(i,j) > 0$ if $(i,j) \in \mathcal{R}$ and $k(i,j) = 0$ if $(i,j) \not\in \mathcal{R}$. We further define the \textbf{reaction-weighted chemical reaction network} associated with $\mathcal{N}$ and $\mathcal{K}$ to be $\mathcal{N}(\mathcal{K}) = (\mathcal{S},\mathcal{C},\mathcal{R},\mathcal{K})$.
\end{definition}
\noindent It is common to incorporate the reaction-weights $k(i,j)$ into the reaction graph as edge weights. This gives rise to an edge-weighted reaction graph $G(V,E(\mathcal{K}))$. For instance, we write
\[\hspace{0.3cm} \underline{G(V,E)}: \hspace{5.4cm} \underline{G(V,E(\mathcal{K}))}:\]\\[-0.9cm]
\[C_1 \begin{array}{c} \vspace{-0.25cm} \longleftarrow \vspace{-0.3cm} \\ \longrightarrow \end{array} C_2 \longrightarrow C_3 \longleftarrow C_4 \hspace{3cm} C_1 \mathop{\stackrel{k(1,2)}{\begin{array}{c} \vspace{-0.25cm} \longleftarrow \vspace{-0.3cm} \\ \longrightarrow \end{array}}}_{k(2,1)} C_2 \stackrel{k(2,3)}{\longrightarrow} C_3 \stackrel{k(4,3)}{\longleftarrow} C_4\]
for the unweighted and weighted reaction graphs of $\mathcal{N}$, respectively.


Defining $\mathbf{x} = (x_1, x_2, \ldots, x_n) \in \mathbb{R}_{\geq 0}^n$ to be the vector of species concentrations, the \emph{mass action system} corresponding to a reaction-weighted chemical reaction network $\mathcal{N}(\mathcal{K}) = (\mathcal{S},\mathcal{C},\mathcal{R},\mathcal{K})$ is given by the system of ordinary differential equations
\begin{equation}
\label{de2}
\frac{d\mathbf{x}}{dt} = Y \cdot A(\mathcal{K}) \cdot \Psi(\mathbf{x})
\end{equation}
where
\begin{enumerate}
\item
The \emph{complex matrix} $Y \in \mathbb{Z}_{\geq 0}^{n \times m}$ is the matrix with columns $Y_{\cdot,i} = y_i$.
\item
The \emph{kinetic} or \emph{Kirchhoff matrix} $A(\mathcal{K}) \in \mathbb{R}^{m \times m}$ is the matrix with entries 
\begin{equation}
\label{kirchhoff}
[A(\mathcal{K})]_{i,j} = \left\{ \begin{array}{ll} \displaystyle{-\sum_{l=1}^m k(i,l)}, \; \; \; & \mbox{for } i = j, \\ k(j,i), & \mbox{for } i \not= j, \end{array} \right.
\end{equation}
for $i,j=1, \ldots, m.$
\item
The \emph{mass action vector} $\Psi(\mathbf{x}) \in \mathbb{R}_{\geq 0}^{m}$ is the vector with entries $[\Psi(\mathbf{x})]_i = \mathbf{x}^{y_i} = \prod_{j=1}^n x_j^{y_{ij}}$, $i=1, \ldots, m$.
\end{enumerate}

\noindent It is known that trajectories of any mass action system are restricted to \emph{stoichiometric compatibility classes} $\mathsf{C}_{\mathbf{x}_0} = (\mathbf{x}_0 + S) \cap \mathbb{R}_{>0}^n$ for all $\mathbf{x}_0 \in \mathbb{R}_{> 0}^n$ \cite{V-H}. 

\begin{remark}
Note that $A(\mathcal{K})$ explicitly relates the topology of the weighted reaction graph to the dynamics. In particular, an off-diagonal element $[A(\mathcal{K})]_{i,j}$ is non-zero if and only if there is a reaction in the network from $C_j$ to $C_i$. 
\end{remark}

\begin{remark}
It is tempting to automatically correspond reaction-weighted networks $(\mathcal{S},\mathcal{C},\mathcal{R},\mathcal{K})$ with mass action systems (\ref{de2}). The theory developed in Section \ref{translationsection}, however, will necessitate the construction of reaction-weighted chemical reaction networks which do not have meaningful interpretations as mass action systems. We will use the notation $\mathcal{B}$ to denote reaction-weight sets which do not necessarily correspond to the kinetic rate constants in a corresponding mass action system.
\end{remark}



\subsection{Generalized Chemical Reaction Networks}

An alternative to mass action kinetics is \emph{power-law formalism}, where the powers of the kinetic terms in the governing equations (\ref{de2}) are allowed to take (potentially non-integer) powers which are not necessarily implied by the stoichiometry of the network \cite{Sa}. A recent graph-based extension of this is the concept of a \emph{generalized chemical reaction network} \cite{M-R}.

\begin{definition}
A \textbf{generalized chemical reaction network} $\mathcal{N} = (\mathcal{S},\mathcal{C},\mathcal{C}_K,\mathcal{R})$ is a chemical reaction network $(\mathcal{S},\mathcal{C},\mathcal{R})$ together with a set of kinetic complexes $\mathcal{C}_K$ which are in one-to-one correspondence with the elements of $\mathcal{C}$.
\end{definition}

When permitted by space, we denote the correspondence between the stoichiometric and kinetic complexes with dotted lines. For example, we write
\begin{equation}
\label{network32}
7 X_1 + X_2 \; \; \; \cdots \; \; \; \; \; \; \; X_1 \; \; \; \begin{array}{c} \vspace{-0.25cm} \longleftarrow \vspace{-0.3cm} \\ \longrightarrow \end{array} \; \; \; X_2 + X_3 \; \; \; \; \; \; \; \cdots \; \; \; X_3
\end{equation}
to imply that the stoichiometric complex $C_1 = X_1$ is associated with the kinetic complex $(C_K)_1 = 7X_1 + X_2$ and that the stoichiometric complex $C_2 = X_2 + X_3$ is associated with the kinetic complex $(C_K)_2 = X_3$. We define properties of the reaction graph $(\mathcal{S},\mathcal{C},\mathcal{R})$ as we do for a standard reaction network. For example, this network has the stoichiometric subspace $S = \mbox{span} \{ (-1,1,1) \}$ and $\delta = 0$. A reaction graph for $(\mathcal{S},\mathcal{C}_K,\mathcal{R})$ can also be defined. We do this by substituting the kinetic complexes for the stoichiometric complexes. For the example network (\ref{network32}), we have
\[7X_1 + X_2 \; \; \; \begin{array}{c} \vspace{-0.25cm} \longleftarrow \vspace{-0.3cm} \\ \longrightarrow \end{array} \; \; \; X_3.\]
We define the \emph{kinetic-order subspace} $S_K$ and the \emph{kinetic-order deficiency} $\delta_K$ as the corresponding quantities for the reaction graph of $(\mathcal{S},\mathcal{C}_K,\mathcal{R})$. For this example, we have $S_K = \mbox{span} \{ (-7,-1,1) \}$ and $\delta_K = 0$.

Given a reaction-weight set $\mathcal{K}$, we define the \emph{generalized reaction-weighted chemical reaction network} associated with $\mathcal{N}$ and $\mathcal{K}$ to be $\mathcal{N}(\mathcal{K}) = (\mathcal{S},\mathcal{C},\mathcal{C}_K,\mathcal{R},\mathcal{K})$. The \emph{generalized mass action system} corresponding to $\mathcal{N}(\mathcal{K})$ is given by
\begin{equation}
\label{gde}
\frac{d\mathbf{x}}{dt} = Y \cdot A(\mathcal{K}) \cdot \Psi_K(\mathbf{x})
\end{equation}
where $\Psi_K(\mathbf{x})$ has entries $[\Psi_K(\mathbf{x})]_i = \mathbf{x}^{(y_K)_i}$, $i =1, \ldots, m$. In other words, a generalized mass action is the mass action system (\ref{de2}) with the monomials $\mathbf{x}^{y_i}$ replaced by the monomials $\mathbf{x}^{(y_K)_i}$. For example, given the reaction-weight set $\mathcal{K} = \{ k(1,2), k(2,1) \}$, the generalized mass action system corresponding to the network (\ref{network32}) is
\[\frac{dx_1} {dt} = -\frac{dx_2}{dt} = -\frac{dx_3}{dt} = -k(1,2) x_1^7x_2 + k(2,1) x_3.\]
Notice that the stoichiometry of the network comes from the stoichiometric complexes $\mathcal{C}$ but the monomials come from the kinetic complexes $\mathcal{C}_K$. Results regarding the existence and location of steady states of generalized mass action systems are contained in \cite{M-F-R-C-S-D,M-R} but will not be summarized here.

\subsection{Kinetically-Relevant Complexes}
\label{complexsection}

It is possible for a source complex to appear in the network $\mathcal{N}(\mathcal{K})$ but not appear in the corresponding mass action system. For example, consider the network
\begin{equation}
\label{34543}
2X_1 \mathop{\stackrel{k(1,2)}{\begin{array}{c} \vspace{-0.25cm} \longleftarrow \vspace{-0.3cm} \\ \longrightarrow \end{array}}}_{k(2,1)} X_1+X_2 \mathop{\stackrel{k(2,3)}{\begin{array}{c} \vspace{-0.25cm} \longleftarrow \vspace{-0.3cm} \\ \longrightarrow \end{array}}}_{k(3,2)} 2X_2.
\end{equation}
For $k(1,2)=k(2,1)=k(2,3)=k(3,2)=1$ we have $(k(2,1)-k(2,3))x_1x_2 = 0$ so that $x_1x_2$ does not appear in (\ref{de2}). For the theory developed in Section \ref{translationsection} we will be interested only in those complexes for which the coefficients of the corresponding monomials $\mathbf{x}^{y_i}$ or $\mathbf{x}^{(y_K)_i}$ do not vanish in (\ref{de2}) or (\ref{gde}). We therefore introduce the following.


\noindent  

\begin{definition}
Consider a regular or generalized reaction-weighted chemical reaction network ($\mathcal{N}(\mathcal{K}) = (\mathcal{S},\mathcal{C},\mathcal{R},\mathcal{K})$ or $\mathcal{N}(\mathcal{K}) = (\mathcal{S},\mathcal{C},\mathcal{C}_K,\mathcal{R},\mathcal{K})$, respectively). We define the \textbf{kinetically-relevant complexes} of $\mathcal{N}(\mathcal{K})$, $\mathcal{C}(\mathcal{K}) \subseteq \mathcal{C},$ to be the set of $i = 1, \ldots, m,$ such that
\begin{equation}
\label{12321}
\mathop{\sum_{j=1}^m}_{j \not= i} k(i,j) \; (y_j - y_i) \not= \mathbf{0}.
\end{equation}
\end{definition}


\noindent Note that $\mathcal{C}(\mathcal{K})$ may depend upon both the structure of $\mathcal{N}(\mathcal{K})$ and the reaction-weighting set $\mathcal{K}$. For example, in (\ref{34543}) we have $\mathcal{C}(\mathcal{K}) = \{ 1, 3\}$ if we choose $k(1,2)=k(2,1)=k(2,3)=k(3,2)=1$; however, we have $\mathcal{C}(\mathcal{K}) = \{ 1, 2, 3\}$ if we choose $k(1,2)=k(2,1)=k(3,2) = 1$ and $k(2,3) = 2$.


\section{Reaction-Weighted Translated Chemical Reaction Networks}
\label{translationsection}



It was observed in \cite{J1} that mass action systems (\ref{de2}) may have related representations as generalized mass action systems (\ref{gde}). In cases where the network underlying the generalized mass action system is better structured (e.g. weakly reversible, lower deficiency, etc.) it may be beneficial to analyze the generalized system rather than the classical one. Consider the following example.\\

\noindent \emph{Example 1:} Consider the reaction-weighted chemical reaction network $\mathcal{N}(\mathcal{K}) = (\mathcal{S},\mathcal{C},\mathcal{R},\mathcal{K})$ and the reaction-weighted generalized chemical reaction network $\tilde{\mathcal{N}}(\tilde{\mathcal{K}}) = (\tilde{\mathcal{S}},\tilde{\mathcal{C}},\tilde{\mathcal{C}}_K,\tilde{\mathcal{R}},\tilde{\mathcal{K}})$ given respectively by:
\begin{equation}
\label{example1}
\begin{array}{c}
\underline{\mathcal{N}(\mathcal{K}) = (\mathcal{S},\mathcal{C},\mathcal{R},\mathcal{K})} \hspace{3cm} \underline{\tilde{\mathcal{N}}(\tilde{\mathcal{K}}) = (\tilde{\mathcal{S}},\tilde{\mathcal{C}},\tilde{\mathcal{C}}_K,\tilde{\mathcal{R}},\tilde{\mathcal{K}})} \hspace{1cm} \vspace{0.2cm} \\
\begin{array}{|c|} \hline \begin{split} X_1 \; & \stackrel{k_1}{\longrightarrow} \; 2X_1 \\ X_1 + X_2 \; & \stackrel{k_2}{\longrightarrow} \; 2X_2 \\ X_2 \; & \stackrel{k_3}{\longrightarrow} \; \emptyset \end{split} \\ \hline \end{array} \hspace{0.4cm}
\begin{array}{c} \vspace{-0.4cm} \\ \Longleftrightarrow \end{array} \hspace{0.4cm} 
\begin{array}{|c|} \hline X_1 \; \; \; \cdots \; \; \; \emptyset \; \stackrel{\tilde{k}_1}{\longrightarrow} \; X_1 \; \; \; \cdots \; \; \; X_1 + X_2 \\ {}_{\tilde{k}_3} \nwarrow \hspace{0.2cm} \swarrow_{\tilde{k}_2} \hspace{1.2cm}\\ \hspace{0.45cm} X_2 \; \; \; \cdots \; \; \; X_2 \\ \hline \end{array}
\end{array}
\end{equation}
\noindent It can be easily verified by expanding (\ref{de2}) or (\ref{gde}), respectively, that the mass action systems correponding to $\mathcal{N}(\mathcal{K})$ and the generalized mass action systems corresponding to $\tilde{\mathcal{N}}(\tilde{\mathcal{K}})$ are identical if we take $k_1 = \tilde{k}_1$, $k_2 = \tilde{k}_2$, and $k_3 = \tilde{k}_3$.

It was noted in \cite{J1} that the process of corresponding $\mathcal{N}(\mathcal{K})$ to $\tilde{\mathcal{N}}(\tilde{\mathcal{K}})$ can be visualized by ``translating'' the complexes of each reaction. For this example, we have
\begin{equation}
\label{translationexample}
\begin{split} X_1 \; &\longrightarrow \; 2X_1 \; \; \; \; \; (-X_1) \; \; \; \; \; \; \; \; \; \; \; \; \;\;\;\; \; \; \; \; \; \; \emptyset \; \longrightarrow \; X_1 \\ X_1 + X_2 \; & \longrightarrow \; 2 X_2 \; \; \; \; \; (- X_2) \; \; \; \; \; \Longrightarrow \; \; \; \; \; \; X_1 \; \longrightarrow \; X_2\\ X_2 \; & \longrightarrow \; \emptyset \hspace{0.37in} (+ \emptyset) \hspace{0.83in} X_2 \; \longrightarrow \; \emptyset.\end{split}
\end{equation}
Notice that this process does not change the reaction vectors, and that we may preserve the monomials in (\ref{de2}) by associating the reactant complexes of the original reactions as the kinetic complexes of the new ones (e.g. associate $X_1$ (left) as the kinetic complex of $\emptyset$ (right), etc.). If we transfer the reaction-weights with the reactions, we arrive at the generalized reaction-weighted network in (\ref{example1}). Notice that $\tilde{\mathcal{N}}(\tilde{\mathcal{K}})$ is weakly reversible while $\mathcal{N}(\mathcal{K})$ is not. This will be one of our primary network properties when understanding ``better'' versus ``poorer'' structure. \hfill $\square$ \\


A further class of systems for which (\ref{de2}) and (\ref{gde}) do not coincide but for which the steady states are identical was also identified in \cite{J1} (see Example 2 in Section \ref{resolvabilitysection}). We introduce the following.

\begin{definition}
\label{dynamicallyequivalent}
Let $\mathcal{N}(\mathcal{K}) = (\mathcal{S},\mathcal{C},\mathcal{R},\mathcal{K})$ denote a reaction-weighted chemical reaction network with corresponding mass action system (\ref{de2}) and $\tilde{\mathcal{N}}(\tilde{\mathcal{K}}) = (\tilde{\mathcal{S}},\tilde{\mathcal{C}},\tilde{\mathcal{C}}_K,\tilde{\mathcal{R}},\tilde{\mathcal{K}})$ denote a generalized reaction-weighted chemical reaction network with corresponding generalized mass action system (\ref{gde}). We will say that $\mathcal{N}(\mathcal{K})$ and $\tilde{\mathcal{N}}(\tilde{\mathcal{K}})$ are:
\begin{enumerate}
\item
\textbf{dynamically equivalent} if (\ref{de2}) and (\ref{gde}) coincide; and
\item 
\textbf{steady state equivalent} if (\ref{de2}) and (\ref{gde}) have the same steady states.
\end{enumerate}
\end{definition}
\noindent We can see that the reaction-weighted networks in (\ref{example1}) of Example 1 are dynamically equivalent.

The author of \cite{J1} called the process outlined in (\ref{translationexample}) \emph{network translation}. In this paper, we adopt a modified definition of network translation which explicitly takes reaction-weights into account.


\begin{definition}
\label{translation}
Consider a reaction-weighted chemical reaction network $\mathcal{N}(\mathcal{K}) = (\mathcal{S}, \mathcal{C}, \mathcal{R},\mathcal{K})$ with reaction-weight set $\mathcal{K} = \{ k(i,j) \; | \; i,j = 1, \ldots, q \}$ and kinetically-relevant complex set $\mathcal{C}(\mathcal{K})$, and a reaction-weighted generalized chemical reaction network $\tilde{\mathcal{N}}(\tilde{\mathcal{B}}) = (\tilde{\mathcal{S}},\tilde{\mathcal{C}},\tilde{\mathcal{C}}_K,\tilde{\mathcal{R}},\tilde{\mathcal{B}})$ with reaction-weight set $\tilde{\mathcal{B}} = \{ \tilde{b}(i',j') \; | \; i',j' = 1, \ldots, \tilde{m} \}$ and kinetically-relevant complex set $\tilde{\mathcal{C}}(\tilde{\mathcal{B}})$. We say $\tilde{\mathcal{N}}(\tilde{\mathcal{B}})$ is a \textbf{reaction-weighted translation} of $\mathcal{N}(\mathcal{K})$ if:
\begin{enumerate}
\item
There is a surjection $h: \mathcal{C}(\mathcal{K}) \mapsto \tilde{\mathcal{C}}(\tilde{\mathcal{B}})$ such that, for every $i \in \mathcal{C}(\mathcal{K})$, there are values $\lambda(i,j') \geq 0$, satisfying:
\begin{enumerate}
\item
$\lambda(i,j') > 0$ implies $(h(i),j') \in \tilde{\mathcal{R}}$;
\item
$\displaystyle{\sum_{ \{i | h(i) = i' \} } \lambda(i,j') = \tilde{b}(i',j');}$ and
\item
$\displaystyle{\mathop{\sum_{j=1}^m}_{j \not= i} k(i,j) (y_j - y_i)  = \mathop{\sum_{j'=1}^{\tilde{m}}}_{j' \not= h(i)} \lambda(i,j') \left(\tilde{y}_{j'} - \tilde{y}_{i'}\right).}$
\end{enumerate}

\item
There is an injection $h_K: \tilde{\mathcal{C}}(\tilde{\mathcal{B}}) \mapsto \mathcal{C}(\mathcal{K})$ so that $h(h_K(i'))=i'$ and $(\tilde{C}_K)_{i'} = C_{h_K(i')}$ for all $i' \in \tilde{\mathcal{C}}(\tilde{\mathcal{B}})$.
\end{enumerate}
The process of finding a generalized network $\tilde{\mathcal{N}}$ which is a reaction-weighted translation of $\mathcal{N}$ is called \textbf{reaction-weighted network translation}.
\end{definition}



To interpret Definition \ref{translation}, we notice that if we sum property $1(c)$ over $i \in \mathcal{C}(\mathcal{K})$ such that $h(i) = i'$ then we have
\begin{equation}
\label{equation2}
\sum_{\{i| h(i) = i'\}} \mathop{\sum_{j=1}^m}_{j \not= i} k(i,j) (y_j-y_i) = \mathop{\sum_{j'=1}^{\tilde{m}}}_{j' \not= i'} \tilde{b}(i',j') (\tilde{y}_{j'} - \tilde{y}_{i'}).
\end{equation}
That is, we may interpret property $1.$ as allowing us to shift reactant complexes in complex space so long as we maintain the \emph{net flux} out of each kinetically-relevant complex in the translation. The technical conditions of property $1.$ guarantee that each translated complex has its flux represented in the network $\tilde{\mathcal{N}}(\tilde{\mathcal{B}})$ which may not be guaranteed by (\ref{equation2}) alone due to cancellation. Property $2.$ requires that we preserve the original source complex as the kinetic complex of the corresponding complex in the translation. The resulting reaction-weighted translated chemical reaction network draws its kinetic complexes from the source complexes of the original network, but may have a significantly different reaction graph. \\

\noindent \emph{Example 1:} We make the assignments $C_1 = X_1$, $C_2 = X_1 + X_2$, $C_3 = X_2$, $C_4 = 2X_1$, $C_5 = 2X_2$, $C_6 = \emptyset$, $\tilde{C}_1 = \emptyset$, $\tilde{C}_2 = X_1$, and $\tilde{C}_3 = X_2$. We can then satisfy the requirements on $h$ and $h_K$ given in Definition \ref{translation} by taking $h(1) = 1$, $h(2) = 2$, $h(3) = 3$, $h_K(1) = 1$, $h_K(2)=2$, and $h_K(3)=3$ so that $(\tilde{C}_K)_1 = C_1 = X_1$, $(\tilde{C}_K)_2 = C_2 = X_1 +X_2$, and $(\tilde{C}_K)_3 = C_3 = X_2$ (property $2.$). The conditions of property $1.$ may be satisfied by taking $\tilde{b}(1,2)=\lambda(1,2) = k_1$, $\tilde{b}(2,3) = \lambda(2,3) = k_2$, and $\tilde{b}(3,1) = \lambda(3,1) = k_3$, and we are done. \hfill $\square$ 


\begin{remark}
Following the conventions of \cite{J1}, we will distinguish objects and sets related to translations with the tilde notation $(\tilde{\cdot})$, e.g. $\tilde{L} \in \tilde{\mathcal{L}}$ for linkage classes, $\tilde{m} = | \tilde{\mathcal{C}} |$ for the number of complexes, etc. In particular, we will denote the structural and kinetic deficiencies of translations by $\tilde{\delta}$ and $\tilde{\delta}_K$, respectively, and denote the kinetic-order subspace by $\tilde{S}_K$. Wherever possible, we will distinguish the indices of the translated complexes by primes, e.g. $i', j' = 1, \ldots, \tilde{m}$, $(i',j') \in \tilde{\mathcal{R}}$, etc. We note that this notation differs from that used in \cite{M-R} for generalized chemical reaction networks.
\end{remark}


\begin{remark}
In general, the reaction-weighting set $\tilde{\mathcal{B}}$ in Definition \ref{translation} consists of computational constructs which do not necessarily correspond to the reaction-weights for any meaningful generalized mass action system. We will reserve the symbol $\tilde{\mathcal{K}}$ for reaction-weighting sets for which the reaction-weighted generalized network $\tilde{\mathcal{N}}(\tilde{\mathcal{K}})$ is either dynamically or steady state equivalent to the original reaction-weighted network $\mathcal{N}(\mathcal{K})$.
\end{remark}

The stoichiometric and kinetic-order subspaces $\tilde{S}$ and $\tilde{S}_K$ for translated chemical reaction networks are characterized by the following result.

\begin{lemma}[Lemma 1, \cite{J1}]
\label{lemma232}
Suppose $\tilde{\mathcal{N}}(\tilde{\mathcal{B}}) = (\tilde{\mathcal{S}},\tilde{\mathcal{C}},\tilde{\mathcal{C}}_K,\tilde{\mathcal{R}},\tilde{\mathcal{B}})$ is a reaction-weighted translation of a reaction-weighted chemical reaction network $\mathcal{N}(\mathcal{K}) = (\mathcal{S},\mathcal{C},\mathcal{R},\mathcal{K})$. Then, if $\tilde{\mathcal{N}}$ is weakly reversible, the stoichiometric subspaces $S$ of $\mathcal{N}$ and $\tilde{S}$ of $\tilde{\mathcal{N}}$ coincide and the kinetic-order subspace $\tilde{S}_K$ of $\tilde{\mathcal{N}}$ is given by
\begin{equation}
\label{tildeS}
\tilde{S}_K = \mbox{span} \left\{ (\tilde{y}_K)_{i'} - (\tilde{y}_K)_{j'} \; | \; \; i',j' \in \tilde{L}_\theta,  \theta =1, \ldots, \tilde{\ell} \right\}
\end{equation}
where $\tilde{L}_\theta, \theta=1, \ldots, \tilde{\ell}$, are the linkage classes of $\tilde{\mathcal{N}}$.
\end{lemma}

\begin{proof}
The result follows from the proof of Lemma 1 in \cite{J1} and the fact that, since the network is weakly reversible, the kinetic and stoichiometric subspaces of $\tilde{\mathcal{N}}$ coincide by Corollary 1 of \cite{H-F}. (Note here that we define the kinetic subspace as in \cite{H-F} and that this is not the same object as the kinetic-order subspace $\tilde{S}_K$.)
\end{proof}

\subsection{Proper Reaction-Weighted Translations}
\label{propersection}

An important subset of reaction-weighted translations is the following, which is modified from Definition 7 in \cite{J1} to accommodate reaction-weights.

\begin{definition}
\label{proper}
Consider a reaction-weighted chemical reaction network $\mathcal{N}(\mathcal{K}) = (\mathcal{S},\mathcal{C},\mathcal{R},\mathcal{K})$ and a reaction-weighted translation $\tilde{\mathcal{N}}(\tilde{\mathcal{B}}) = (\tilde{\mathcal{S}},\tilde{\mathcal{C}},\tilde{\mathcal{C}}_K,\tilde{\mathcal{R}},\tilde{\mathcal{B}})$. We will say $\tilde{\mathcal{N}}(\tilde{\mathcal{B}})$ is a \textbf{proper reaction-weighted translation} of $\mathcal{N}(\mathcal{K})$ if $h: \mathcal{C}(\mathcal{K}) \mapsto \tilde{\mathcal{C}}(\tilde{\mathcal{B}})$ is injective as well as surjective. A reaction-weighted translation $\tilde{\mathcal{N}}(\tilde{\mathcal{B}})$ will be called \textbf{improper} if it is not proper.
\end{definition}

\noindent That is, a reaction-weighted translation is proper if every kinetically-relevant complex in $\mathcal{N}(\mathcal{K})$ corresponds to exactly one kinetically-relevant complex in $\tilde{\mathcal{N}}(\tilde{\mathcal{B}})$. Notice that, if $\tilde{\mathcal{N}}$ is proper, properties $1(a-c)$ in Definition \ref{translation} and (\ref{equation2}) are equivalent. For proper translations, we also have $h_K = h^{-1}$.

The following result is modified from a result proved in \cite{J1}.

\begin{lemma}[Lemma 2, \cite{J1}]
\label{lemma231}
Suppose $\tilde{\mathcal{N}}(\tilde{\mathcal{B}}) = (\tilde{\mathcal{S}},\tilde{\mathcal{C}},\tilde{\mathcal{C}}_K,\tilde{\mathcal{R}},\tilde{\mathcal{B}})$ is a proper reaction-weighted translation of the reaction-weighted chemical reaction network $\mathcal{N}(\mathcal{K}) = (\mathcal{S},\mathcal{C},\mathcal{R},\mathcal{K})$. Then the reaction-weighted network $\tilde{\mathcal{N}}(\tilde{\mathcal{K}}) = (\tilde{\mathcal{S}},\tilde{\mathcal{C}},\tilde{\mathcal{C}}_K,\tilde{\mathcal{R}},\tilde{\mathcal{K}})$ with $\tilde{\mathcal{K}} = \tilde{\mathcal{B}}$ is dynamically equivalent to $\mathcal{N}(\mathcal{K})$.
\end{lemma}

\begin{proof}
The result follows immediately from the proof of Lemma 2 in \cite{J1}, the observation that properties $1(a-c)$ in Definition \ref{translation} and (\ref{equation2}) coincide for proper translations, and Definition \ref{dynamicallyequivalent}.
\end{proof}

\noindent \emph{Example 1:} It can be easily seen that the translation scheme (\ref{translationexample}) in Example 1 results in a proper translation (\ref{example1}) for any reaction-weightings $k_1, k_2,$ and $k_3$. It was previously noted that the two reaction-weighted networks have the same dynamics. This is consistent with the application of Lemma \ref{lemma231}. \hfill $\square$ \\


\subsection{Improper Reaction-Weighted Translations}
\label{resolvabilitysection}


It was noted in \cite{J1} that any generalized mass action system (\ref{gde}) corresponding to an \emph{improper} reaction-weighted translation $\tilde{\mathcal{N}}(\tilde{\mathcal{B}}) = (\tilde{\mathcal{S}},\tilde{\mathcal{C}},\tilde{\mathcal{C}}_K,\tilde{\mathcal{R}},\tilde{\mathcal{B}})$ must necessarily differ from the mass action system (\ref{de2}) corresponding to the original network $\mathcal{N}(\mathcal{K}) = (\mathcal{S},\mathcal{C},\mathcal{R},\mathcal{K})$. A result analogous to Lemma \ref{lemma231} is therefore not possible. Nevertheless, conditions were given in \cite{J1} under which a \emph{rescaled} reaction-weighting set $\tilde{\mathcal{K}}$ could be constructed so that $\mathcal{N}(\mathcal{K})$ and $\tilde{\mathcal{N}}(\tilde{\mathcal{K}})$ shared the same steady state set. Consider the following example.\\

\noindent \emph{Example 2:} Consider the reaction-weighted chemical reaction network $\mathcal{N}(\mathcal{K})$ with reaction-weight set $\mathcal{K} = \{ k_i > 0 \; | \; i=1, \ldots, 14 \}$ corresponding to the reactions as labeled:
\begin{equation}
\label{system3}
\begin{split}
& \displaystyle{X_1 \mathop{\stackrel{1}{\rightleftarrows}}_{2} X_2 \mathop{\stackrel{3}{\rightleftarrows}}_{4} X_3 \stackrel{5}{\rightarrow} X_4} \\
& \displaystyle{X_4 + X_5 \mathop{\stackrel{6}{\rightleftarrows}}_{7} X_6 \stackrel{8}{\rightarrow} X_2+X_7} \\
& \displaystyle{X_3 + X_7 \mathop{\stackrel{9}{\rightleftarrows}}_{10} X_8 \stackrel{11}{\rightarrow} X_3 + X_5} \\
& \displaystyle{X_1 + X_7 \mathop{\stackrel{12}{\rightleftarrows}}_{13} X_9 \stackrel{14}{\rightarrow} X_1 + X_5.}
\end{split}
\end{equation}
This network has been studied by Shinar and Feinberg in \cite{Sh-F} and by P\'{e}rez Mill\'{a}n \emph{et al.} in \cite{M-D-S-C}. (Further details are contained in the Supplemental Material.) It was noted by Johnston in \cite{J1} that the translation scheme 
\begin{equation}
\label{2translation}
\begin{split}
& \displaystyle{X_1 \mathop{\stackrel{1}{\rightleftarrows}}_{2} X_2 \mathop{\stackrel{3}{\rightleftarrows}}_{4} X_3 \stackrel{5}{\rightarrow} X_4} \; \; \; \; \; \; \; \; \; \; \; \; \; \; \; \; \; \; \; \; (+ X_1 + X_3 + X_5) \\
& \displaystyle{X_4 + X_5 \mathop{\stackrel{6}{\rightleftarrows}}_{7} X_6 \stackrel{8}{\rightarrow} X_2 + X_7} \; \; \; \; \; \; \; \; \; \; \; \; (+ X_1 + X_3) \\
& \displaystyle{X_3 + X_7 \mathop{\stackrel{9}{\rightleftarrows}}_{10} X_8 \stackrel{11}{\rightarrow} X_3 + X_5} \; \; \; \; \; \; \; \; \; \; \; \; (+ X_1 + X_2)\\
& \displaystyle{X_1 + X_7 \mathop{\stackrel{12}{\rightleftarrows}}_{13} X_9 \stackrel{14}{\rightarrow} X_1 + X_5 \; \; \; \; \; \; \; \; \; \; \; \; (+ X_2 + X_3)}
\end{split}
\end{equation}
yields the following reaction-weighted translation $\tilde{\mathcal{N}}(\tilde{\mathcal{B}})$, where $\tilde{b}_i = k_i$, $i=1,\ldots, 14$:
\begin{equation}
\label{system4}
\begin{split} & 2X_1 + X_3 + X_5 \mathop{\stackrel{1}{\rightleftarrows}}_{2} X_1 + X_2 + X_3 + X_5 \mathop{\stackrel{3}{\rightleftarrows}}_{4} X_1 + 2X_3 + X_5 \\ & \; \; \; \; \; \; \; \; \; \; \; \; \; \; \; \; \; \; \; \; \; \; \; \; \nearrow \hspace{-0.1cm} {}_{14} \; \; \; \; \; \; \; \; \; \; \; \; \; \uparrow_{11} \; \; \; \; \; \; \; \; \; \; \; \; \; \; \; \; \; \; \; \; \; \; \; \; \; \; \; \downarrow_{5} \\ & \; \; X_2 + X_3 + X_9 \; \; \; \; \; \; \; \; \; \; \; X_1 + X_2 + X_8 \; \; \; \; \; \; \; X_1 + X_3 + X_4 + X_5 \\ & \; \; \; \; \; \; \; \; \; \; \; \; \; \; \; \; \; \; \; {}_{12} \hspace{-0.1cm} \nwarrow \hspace{-0.15cm} \searrow \hspace{-0.1cm} {}^{13} \; \; \; \; \; \; \; \; \; \; \; {}_9 \uparrow \downarrow {}_{10} \; \; \; \; \; \; \; \; \; \; \; \; \; \; \; \; \; \; \; \; \; \; {}_7 \uparrow \downarrow {}_6 \\ & \; \; \; \; \; \; \; \; \; \; \; \; \; \; \; \; \; \; \; \; \; \; \; \; \; \; \; \; X_1 + X_2 + X_3 + X_7 \stackrel{8}{\leftarrow} X_1 + X_3 + X_6. \end{split}
\end{equation}
\noindent The translated network $\tilde{\mathcal{N}}(\tilde{\mathcal{B}})$ is improper since the source complexes $X_3 + X_7$ and $X_1 + X_7$ are both translated to $X_1 + X_2 +X_3 + X_7$ but we may only keep one as the corresponding kinetic complex. Notice that, regardless of the choice of kinetic complex corresponding to $X_1 + X_2 + X_3 + X_7$, the generalized system (\ref{gde}) corresponding to (\ref{system4}) is not dynamically equivalent to the system (\ref{de2}) corresponding to (\ref{system3}).

It was shown in \cite{J1} that, if we choose $X_3 + X_7$ as the kinetic complex of the stoichiometric complex $X_1 + X_2 + X_3 + X_7$, the reaction-weighted networks $\mathcal{N}(\mathcal{K})$ and $\tilde{\mathcal{N}}(\tilde{\mathcal{K}})$ given in (\ref{system3}) and (\ref{system4}), respectively, are \emph{steady state equivalent} for $\tilde{k}_i = k_i$, $i=1, \ldots, 14$, $i \not= 12$, and 
\begin{equation}
 \label{3100}
\tilde{k}_{12} = \left( \frac{k_2(k_4+k_5)}{k_1k_3} \right) k_{12}.
\end{equation}
In other words, the systems (\ref{de2}) and (\ref{gde}) coincide at steady state after a rescaling of the rate parameter $k_{12}$. Notice importantly that the set $\tilde{\mathcal{K}}$ does not satisfy (\ref{equation2}), and that substituting the set $\tilde{\mathcal{B}}$ in (\ref{gde}) does not produce a system which is steady state equivalent with (\ref{de2}). That is, while corresponding to the same network structure, the reaction-weight sets $\tilde{\mathcal{K}}$ and $\tilde{\mathcal{B}}$ serve distinct and non-interchangeable functions.
\hfill $\square$ \\


Algebraic conditions on the reaction-weight set $\tilde{\mathcal{B}}$ which are sufficient to guarantee such a rescaling can be made were derived in \cite{J1}. The conditions were called \emph{resolvability} conditions, which we do not reproduce here (some details are contained in Appendix \hyperlink{appendixa}{A}). Instead, we consider the following broader definition.

\begin{definition}
\label{resolvable}
Let $\tilde{\mathcal{N}}(\tilde{\mathcal{B}}) = (\tilde{\mathcal{S}},\tilde{\mathcal{C}},\tilde{\mathcal{C}}_K,\tilde{\mathcal{R}},\tilde{\mathcal{B}})$ denote an improper reaction-weighted translation of the reaction-weighted chemical reaction network $\mathcal{N}(\mathcal{K}) = (\mathcal{S},\mathcal{C},\mathcal{R},\mathcal{K})$. We will say that $\mathcal{N}(\mathcal{K})$ and  $\tilde{\mathcal{N}}(\tilde{\mathcal{B}})$ are \textbf{steady state resolvable} if there is a reaction-weight set $\tilde{\mathcal{K}}$ such that $\mathcal{N}(\mathcal{K})$ and $\tilde{\mathcal{N}}(\tilde{\mathcal{K}})$ are steady state equivalent.
\end{definition}



\noindent \emph{Example 2:} We can see that the reaction-weighted networks $\mathcal{N}(\mathcal{K})$ and $\tilde{\mathcal{N}}(\tilde{\mathcal{B}})$ are \emph{steady state resolvable} since a reaction-weight set $\tilde{\mathcal{K}}$ with the same structure as $\tilde{\mathcal{B}}$ may be selected so that $\mathcal{N}(\mathcal{K})$ and $\tilde{\mathcal{N}}(\tilde{\mathcal{K}})$ are \emph{steady state equivalent}. \hfill $\square$ \\

\subsection{Sufficient Conditions for Steady State Resolvability}

In this section, we consider the following problem: given a reaction-weighted chemical reaction network $\mathcal{N}(\mathcal{K})$ and a translation $\tilde{\mathcal{N}}(\tilde{\mathcal{B}})$, are there sufficient conditions on the \emph{reaction graph of the translation alone} which guarantee that $\tilde{\mathcal{N}}(\tilde{\mathcal{B}})$ is steady state resolvable with $\mathcal{N}(\mathcal{K})$? This approach differs from that taken in \cite{J1}, where the resolvability conditions were algebraic in nature. We will answer the question affirmatively with Theorem \ref{maintheorem}. We will use Example 2 introduced in Section \ref{resolvabilitysection} as a running example.

We begin by introducing the following definitions.


\begin{definition}
\label{cs}
Suppose $\tilde{\mathcal{N}}(\tilde{\mathcal{B}}) = (\tilde{\mathcal{S}},\tilde{\mathcal{C}},\tilde{\mathcal{C}}_K,\tilde{\mathcal{R}},\tilde{\mathcal{B}})$ is a reaction-weighted improper translation of a reaction-weighted chemical reaction network $\mathcal{N}(\mathcal{K}) = (\mathcal{S},\mathcal{C},\mathcal{R},\mathcal{K})$. Then:
\begin{enumerate}
\item
The \textbf{improper complex set} $\tilde{\mathcal{C}}_I \subseteq \tilde{\mathcal{C}}(\tilde{\mathcal{B}})$ is given by
\begin{equation}
\label{conflictedset}
\tilde{\mathcal{C}}_I = \{ k' \in \tilde{\mathcal{C}}(\tilde{\mathcal{B}}) \; | \; h(i)=h(j)=k' \mbox{ for some } i,j \in \mathcal{C}(\mathcal{K}), i \not= j\}.
\end{equation}
\item
The $k'$-\textbf{unresolved complex set} $h^{-1}(k') \subseteq \mathcal{C}(\mathcal{K})$ is given by
\begin{equation}
\label{unresolved}
h^{-1}(k') = \{ i \in \mathcal{C}(\mathcal{K}) \; | \; h(i) = k' \mbox{ where } k' \in \tilde{\mathcal{C}}_I \}.
\end{equation}
\item
The \textbf{improper subspace} $\tilde{S}_I$ of $\tilde{\mathcal{N}}(\tilde{\mathcal{B}})$ is given by
\begin{equation}
\label{impropersubspace}
\tilde{S}_I = \mbox{span} \left\{ y_j - y_i \; | \; i,j \in h^{-1}(k') \mbox{ where } k' \in \tilde{\mathcal{C}}_I \right\}.
\end{equation}
\end{enumerate}
\end{definition}
\noindent Note that the definition of the improper subspace $\tilde{S}_I$ differs notationally from the corresponding definition in \cite{J1} (Definition 9). It can easily be checked that the two definitions are equivalent.\\

\noindent \emph{Example 2:} Consider the reaction-weighted network $\mathcal{N}(\mathcal{K})$ given by (\ref{system3}) and the generalized reaction-weighted network given by $\tilde{\mathcal{N}}(\tilde{\mathcal{B}})$. We index the complexes of $\mathcal{N}(\mathcal{K})$ according to:
\[\begin{split} &C_1 = X_1, \; C_2 = X_2, \; C_3 = X_3, \; C_4 =X_4 + X_5, \; C_5 = X_6 ,\\
&C_6 = X_3 + X_7, \; C_7 = X_8, \; C_8 = X_1 + X_7, \; C_9 = X_9,\\
&C_{10} = X_4, \; C_{11} = X_2 + X_7, \; C_{12} = X_3 + X_5, \; C_{13} = X_1 + X_5.\end{split}\]
and the complexes of $\tilde{\mathcal{N}}(\tilde{\mathcal{B}})$ according to:
\[\begin{split}& \tilde{C}_1 = 2X_1 + X_3 + X_5, \; \tilde{C}_2 = X_1 + X_2 + X_3 + X_5, \; \tilde{C}_3 = X_1 + 2X_3 + X_5,\\
& \tilde{C}_4 = X_1 + X_3 + X_4 + X_5, \; \tilde{C}_5 = X_1 + X_3 + X_6, \; \tilde{C}_6 = X_1 + X_2 + X_3 + X_7,\\
&\tilde{C}_7 = X_1 + X_2 + X_8, \; \tilde{C}_8 = X_2 + X_3 + X_9.\end{split}\]
We furthermore index the kinetic complex set $\tilde{\mathcal{C}}_K$ according to:
\begin{equation}
\label{examplekinetic}
\begin{split} & (\tilde{C}_K)_1 = X_1, \; (\tilde{C}_K)_2 = X_2, \; (\tilde{C}_K)_3 = X_3, \; (\tilde{C}_K)_4 = X_4 + X_5, \\ & (\tilde{C}_K)_5 = X_6, \; (\tilde{C}_K)_6 = X_3 + X_7, \; (\tilde{C}_K)_7 = X_8, \; (\tilde{C}_K)_8 = X_9.\end{split}
\end{equation}
Notice that we have chosen $(\tilde{C}_K)_6 = C_6 = X_3 + X_7$ but could have chosen $(\tilde{C}_K)_6 = C_8 = X_1 + X_7$ by property $2.$ of Definition \ref{translation}. Since we have $h(6) = 6$ and $h(8) = 6$, it follows by (\ref{conflictedset}), (\ref{unresolved}), and (\ref{impropersubspace}), that $\tilde{\mathcal{C}}_I = \{ 6 \}$, $h^{-1}(6) = \{ 6, 8 \}$, and $\tilde{S}_I = \mbox{span} \{ y_8 - y_6 \} = \mbox{span} \{ (1,0,-1,0,0,0,0,0,0) \}.$ \hfill $\square$ \\


The relationship between the kinetic-order subspace $\tilde{S}_K$ and the improper subspace $\tilde{S}_I$ was shown in \cite{J1} to be crucial to obtaining steady state resolvability of $\tilde{\mathcal{N}}(\tilde{\mathcal{B}})$. We omit the algebraic details here. We instead introduce the following. (The connection between these definitions and conditions to resolvability as defined in \cite{J1} is contained in Appendix \hyperlink{appendixa}{A}.)

\begin{definition}
\label{resolving2}
Let $\tilde{\mathcal{N}}(\tilde{\mathcal{B}}) = (\tilde{\mathcal{S}},\tilde{\mathcal{C}},\tilde{\mathcal{C}}_K,\tilde{\mathcal{R}},\tilde{\mathcal{B}})$ be an improper reaction-weighted translation of a reaction-weighted chemical reaction network $\mathcal{N}(\mathcal{K}) = (\mathcal{S},\mathcal{C},\mathcal{R},\mathcal{K})$. Suppose furthermore that $\tilde{\mathcal{N}}$ is weakly reversible and that $\tilde{S}_I \subseteq \tilde{S}_K$. Then we say $\tilde{\mathcal{C}}_R \subseteq \tilde{\mathcal{C}}(\tilde{\mathcal{B}})$ is a \textbf{resolving complex set} of $\tilde{\mathcal{N}}(\tilde{\mathcal{B}})$ if, for every $i,j \in h^{-1}(k')$ where $k' \in \tilde{\mathcal{C}}_I$, there is a set of constants $c(i',j'),$ $i',j' = 1, \ldots, \tilde{m}$, $i' < j'$, such that:
\begin{enumerate}
\item
$\displaystyle{y_j - y_i = \mathop{\sum_{i',j'=1}^{\tilde{m}}}_{i' < j'} c(i',j') (y_{h_K(j')}-y_{h_K(i')})}$;
\item
$c(i',j') \not= 0$ implies $i',j' \in \tilde{L}_{\theta}$ for some linkage class $\tilde{L}_{\theta}$ of $\tilde{\mathcal{N}}(\tilde{\mathcal{B}})$; and
\item
$c(i',j') \not= 0$ implies $i', j' \in \tilde{\mathcal{C}}_R$.
\end{enumerate}
\end{definition}


\noindent \emph{Example 2:} Notice that we have
\begin{equation}
\label{11111}
y_8 - y_6 = (1,0,-1,0,0,0,0,0) = (\tilde{y}_K)_1 - (\tilde{y}_K)_3.
\end{equation}
It follows that we may satisfy condition $1.$ of Definition \ref{resolving2} by choosing $c(1,3) = 1$ and $c(i',j')$ for all other $i',j' = 1, \ldots, 8$. We may therefore take $\tilde{\mathcal{C}}_R = \{ 1,3 \}$ as our resolving constant set.

Intuitively, at steady state we may ``resolve'' the competition between the two complexes translated to $\tilde{C}_6$ by appealing to the resolving kinetic complexes $(\tilde{C}_K)_1 = X_1$ and $(\tilde{C}_K)_3 = X_3$. Rearranging condition (\ref{11111}) gives
\begin{equation}
\label{11112}
\mathbf{x}^{y_8} =  \left( \frac{\mathbf{x}^{(\tilde{y}_K)_1}}{\mathbf{x}^{(\tilde{y}_K)_3}} \right) \mathbf{x}^{y_6} \; \; \; \Longrightarrow \; \; \; x_1 x_7 = \left( \frac{x_1}{x_3} \right) x_3 x_7.
\end{equation}
The key insight is the monomials $x_1$, $x_3$, and $x_3x_7$ correspond to kinetic complexes in (\ref{examplekinetic}) while the monomial $x_1 x_7$ does not. This is the monomial which needs to be ``resolved'' since it appears in the original equations (\ref{de2}) but not in the generalized equations (\ref{gde}). \hfill $\square$

\begin{remark}
The resolving complex set $\tilde{\mathcal{C}}_R$ corresponds to the kinetic complexes which are required to related the vectors in $\tilde{S}_I$ to those in $\tilde{S}_K$. Note that $\tilde{S}_I \subseteq \tilde{S}_K$ gives a sufficient condition for $\tilde{\mathcal{C}}_R \not= \emptyset$ by Lemma 3 of \cite{J1}. Condition 2 follows from Lemma 3 of \cite{J1} and Lemma \ref{lemma231} here.
\end{remark}

We now state the main technical result of the paper. The proof can be found in Appendix \hyperlink{appendixb}{B}. We also present there an alternative statement of the Theorem which may be more intuitive to some readers (Lemma \ref{mainlemma}). The statement presented here is more amenable to the computation procedure of Section \ref{milpsection}.

\begin{theorem}
\label{maintheorem}
Let $\tilde{\mathcal{N}}(\tilde{\mathcal{B}}) = (\tilde{\mathcal{S}},\tilde{\mathcal{C}},\tilde{\mathcal{C}}_K,\tilde{\mathcal{R}},\tilde{\mathcal{B}})$ denote an improper reaction-weighted translation of a reaction-weighted chemical reaction network $\mathcal{N}(\mathcal{K}) = (\mathcal{S},\mathcal{C},\mathcal{R},\mathcal{K})$. Suppose that $\tilde{\mathcal{N}}$ is weakly reversible, $\tilde{\delta} = 0$, and $\tilde{S}_I \subseteq \tilde{S}_K$. Suppose furthermore that there are complex sets $\tilde{\mathcal{C}}^*, \tilde{\mathcal{C}}^{**} \subseteq \tilde{\mathcal{C}}$, and reaction sets $\tilde{\mathcal{R}}^* \subseteq \tilde{\mathcal{R}}$ and $\tilde{\mathcal{R}}^{**} \subseteq \tilde{\mathcal{C}}^{**} \times \tilde{\mathcal{C}}^*$ such that:
\begin{enumerate}
\item
$\tilde{\mathcal{C}}_I \subseteq \tilde{\mathcal{C}}^*$ and $\tilde{\mathcal{C}}_R \cap \tilde{\mathcal{C}}^* = \emptyset$;
\item
$(i',j') \in \tilde{\mathcal{R}}^*$ if and only if $i' \in \tilde{\mathcal{C}}^*$ and $(i',j') \in \tilde{\mathcal{R}}$;
\item
$|\tilde{\mathcal{C}}^{**}| = |\tilde{\mathcal{L}}^*|$ where $\tilde{\mathcal{L}}^*$ is the set of linkage classes of the network $(\tilde{\mathcal{S}},\tilde{\mathcal{C}}^* \cup \tilde{\mathcal{C}}^{**},\tilde{\mathcal{R}}^* \cup \tilde{\mathcal{R}}^{**})$; and
\item
The network $(\tilde{\mathcal{S}},\tilde{\mathcal{C}}^* \cup \tilde{\mathcal{C}}^{**},\tilde{\mathcal{R}}^* \cup \tilde{\mathcal{R}}^{**})$ is weakly reversible.
\end{enumerate}
Then $\mathcal{N}(\mathcal{K})$ and $\tilde{\mathcal{N}}(\tilde{\mathcal{B}})$ are steady state resolvable.
\end{theorem}

The conditions of Theorem \ref{maintheorem} may be understood in the following way. We construct a network $(\tilde{\mathcal{S}},\tilde{\mathcal{C}}^* \cup \tilde{\mathcal{C}}^{**},\tilde{\mathcal{R}}^* \cup \tilde{\mathcal{R}}^{**})$ which tracks paths from complexes in $\tilde{\mathcal{C}}_I$ to complexes in $\tilde{\mathcal{C}}_R$. The four technical conditions of Theorem \ref{maintheorem} guarantee that:
\begin{enumerate}
\item[(1-2)]
We consider all possible paths which originate at complexes in $\tilde{\mathcal{C}}_I$ and force them to stop if they reach a complex in $\tilde{\mathcal{C}}_R$ (although they may stop earlier).
\item[(3-4)]
By continuing these paths, we attempt to construct a component (i.e. linkage class) which has a unique sink. If such a component can be constructed, this sink may then be connected to the rest of the component (by a reaction in $\tilde{\mathcal{R}}^{**}$) to create a weakly reversible network.
\end{enumerate}

\noindent The property of reaching a unique sink before passing through any complex in $\tilde{\mathcal{C}}_R$ is key to the proof of Lemma \ref{mainlemma} for guaranteeing reaction-weights exist for which $\tilde{\mathcal{N}}(\tilde{\mathcal{K}})$ is steady state equivalent with $\mathcal{N}(\mathcal{K})$. The full statement of Lemma \ref{mainlemma}, and a proof is its equivalence to Theorem \ref{maintheorem}, are given in Appendix \hyperlink{appendixb}{B}.\\

\begin{figure}[h]
\centering
\includegraphics[width=8.5cm]{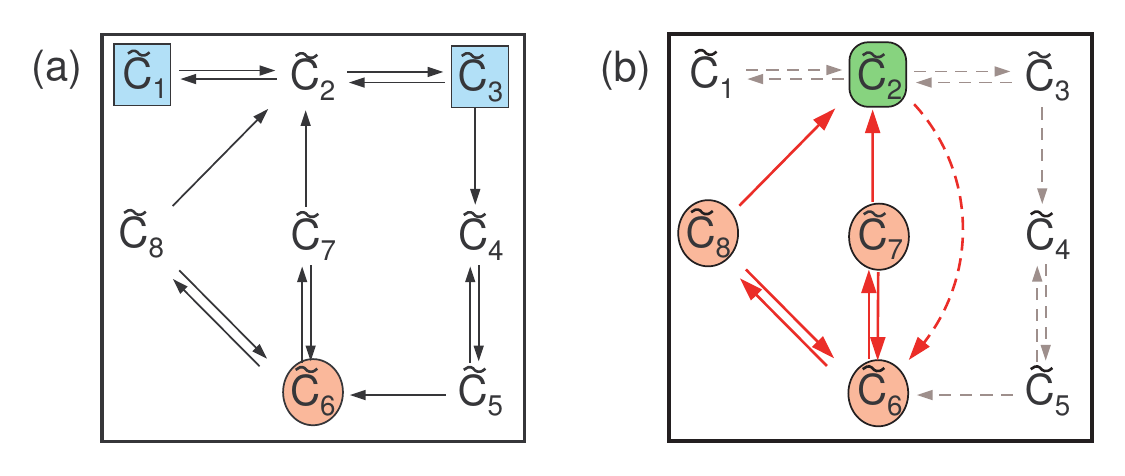}
\caption{\small In (a), we have the improper reaction-weighted network $\tilde{\mathcal{N}}(\tilde{\mathcal{B}}) = (\tilde{\mathcal{S}},\tilde{\mathcal{C}},\tilde{\mathcal{C}}_K,\tilde{\mathcal{R}},\tilde{\mathcal{B}})$ corresponding to Example 2. Highlighted are the improper complex set $\tilde{\mathcal{C}}_I = \{ 6 \}$ (pink) and resolving complex set $\tilde{\mathcal{C}}_R = \{ 1, 3\}$ (blue). In (b), we have the network $(\tilde{\mathcal{S}},\tilde{\mathcal{C}}^{*} \cup \tilde{\mathcal{C}}^{**},\tilde{\mathcal{R}}^{*} \cup \tilde{\mathcal{R}}^{**})$ where $\tilde{\mathcal{C}}^* = \{ 6, 7, 8 \}$ (pink), $\tilde{\mathcal{C}}^{**} = \{ 2 \}$, $\tilde{\mathcal{R}}^* = \{ (6,7),(6,8),(7,2),(7,6),(8,2),(8,6)\}$ (solid red arrows), and $\tilde{\mathcal{R}}^{**} = \{ (2,6) \}$ (dotted red arrow). It is clear that $| \tilde{\mathcal{C}}^{**} | = | \tilde{\mathcal{L}}^* |$ and $(\tilde{\mathcal{S}},\tilde{\mathcal{C}}^{*} \cup \tilde{\mathcal{C}}^{**},\tilde{\mathcal{R}}^{*} \cup \tilde{\mathcal{R}}^{**})$ is weakly reversible. Notice that $\tilde{\mathcal{C}}^{*}$ does contain any complex in $\tilde{\mathcal{C}}_R$ (although it is permissible to have $\mathcal{C}^{**}$ contain such a complex). Also notice that reactions in $\tilde{\mathcal{R}}^{**}$ need not be in the original network, nor be singletons, but that, by the construction of $\tilde{\mathcal{C}}^{**}$, they do need to originate at a sink of linkage class in $\tilde{\mathcal{L}}^{*}$. Since $\tilde{\delta} = 0$, Theorem \ref{maintheorem} applies so that we are guaranteed $\tilde{\mathcal{N}}(\tilde{\mathcal{B}})$ and $\mathcal{N}(\mathcal{K})$ are steady state resolvable.}
\label{figure1}
\end{figure}

\noindent \emph{Example 2:} The required sets $\tilde{\mathcal{C}}^*$, $\tilde{\mathcal{C}}^{**}$, $\tilde{\mathcal{R}}^*$, and $\tilde{\mathcal{R}}^{**}$ for application of Theorem \ref{maintheorem} are given in Figure \ref{figure1}(b). 

\section{Mixed-Integer Linear Programming Framework}
\label{milpsection}

As noted in the Section \ref{introduction}, when attempting to apply Definition \ref{translation} we do not have the network structure of the translation $\tilde{\mathcal{N}}(\tilde{\mathcal{B}})$ given to us; rather, we much \emph{find} it. In \cite{J1}, Johnston presented a heuristic algorithm for determining network translations based on the network's decomposition in elementary flux modes. This method, however, required determining the network's stoichiometriometric generators and then enumerating all possible reaction cycles on the support of these generators. In the case of improper translations, it further required checking algebraic conditions on the network's reaction weights in order to guarantee resolvability. As such, it did not readily lend itself to computational algorithmization or implementation.

We instead adopt here the methodology introduced by Szederk\'{e}nyi in \cite{Sz2}. In that paper, the author introduced a method for determining dynamically equivalent realizations of mass action (or general polynomial) systems when the network structure of the desired realization is unknown. It was shown that the problem of determining a realization with the greatest or least number of reactions (a \emph{dense} or \emph{sparse} realization, respectively) could be formulated as a \emph{mixed-integer linear programming} (MILP) problem.  In subsequent papers, Szederk\'{e}nyi and various collaborators gave additional constraint sets capable of restricting to detailed and complex balanced mass action systems \cite{Sz-H,Sz-H-P}, weakly reversible mass action systems \cite{Sz-H-T,J-S4}, and linearly conjugate mass action systems \cite{J-S4,J-S5,J-S6}.

In this section, we build upon this framework to correspond mass action systems to generalized mass action systems through reaction-weighted network translation. In particular, we detail the logical equivalences corresponding to Definition \ref{translation} and Theorem \ref{maintheorem}. The corresponding MILP code is contained in Appendix \hyperlink{appendixc}{C}. We will also need to re-iterate the results of Johnston \emph{et al.} in \cite{J-S4} and \cite{J-S6}, respectively, pertaining to weak reversibility and minimizing the deficiency of realized networks.

\subsection{Initialization of MILP procedure}

Suppose we have a reaction-weighted chemical reaction network $\mathcal{N}(\mathcal{K}) = (\mathcal{S},\mathcal{C},\mathcal{R},\mathcal{K})$ and wish to determine a reaction-weighted translation $\tilde{\mathcal{N}}(\tilde{\mathcal{B}}) = (\tilde{\mathcal{S}},\tilde{\mathcal{C}},\tilde{\mathcal{C}}_K,\tilde{\mathcal{R}},\tilde{\mathcal{B}})$. We first reorder the complexes of $\mathcal{C}$ so that the first $q = | \mathcal{C}(\mathcal{K}) |$ complexes correspond to be kinetically-relevant complexes. We let $\tilde{m}$ denote the number of potential kinetically-relevant complexes $\tilde{\mathcal{C}}(\tilde{\mathcal{B}})$. We initialize the following matrices:
\begin{enumerate}
\item[(1)]
The matrix $Y \in \mathbb{Z}_{\geq 0}^{q \times n}$ with entries $[Y]_{\cdot,i} = y_i$ where $y_{i}, i=1, \ldots, q$, are the stoichiometric vectors of the kinetically-relevant complexes $\mathcal{C}(\mathcal{K})$.
\item[(2)]
The matrix $\tilde{Y} \in \mathbb{Z}_{\geq 0}^{\tilde{m} \times n}$ with entries $[\tilde{Y}]_{\cdot,i} = \tilde{y}_{i'}$ where $\tilde{y}_{i'}, i' =1, \ldots, \tilde{m}$, are the stoichiometric vectors of the potential set of kinetically-relevant complexes $\tilde{\mathcal{C}}(\tilde{\mathcal{B}})$.
\item[(3)]
The matrix $M \in \mathbb{R}^{q \times n}$ with entries $[M]_{\cdot,i} = [Y \cdot A(\mathcal{K})]_{\cdot,i}$, $i=1, \ldots, q$, where $A(\mathcal{K})$ is the Kirchoff matrix of $\mathcal{N}(\mathcal{K})$. That is, it is the restriction of $Y \cdot A(\mathcal{K})$ to the kinetically-relevant complex set $\mathcal{C}(\mathcal{K})$.
\end{enumerate}

\noindent We note the following:
\begin{itemize}
\item
The kinetically-relevant complexes $\tilde{\mathcal{C}}(\tilde{\mathcal{B}})$ which compose $\tilde{Y}$ need not overlap with the kinetically-relevant complexes $\mathcal{C}(\mathcal{K})$ which compose $Y$ as they did in \cite{Sz2,Sz-H,Sz-H-P,Sz-H-T,J-S4,J-S5,J-S6,R-S-H}. We leave the selection of the candidate stoichiometric complexes as an avenue for future work.
\item
The complexes in $\tilde{Y}$ may not appear in any reaction selected by the computational algorithm and therefore may not appear in $\tilde{\mathcal{N}}(\tilde{\mathcal{B}})$. This is a slight abuse of convention within CRNT literature but will be allowed in the present context. It was shown in \cite{J-S6} that such this abuse of convention does not alter the deficiency of the network or the property of weak reversibility.
\item
In contrast to the results of \cite{J1}, the method presented here determines a translation \emph{for a specific set of chosen rate constants only}. In particular, the reaction-weights of $\mathcal{N}(\mathcal{K})$ must be numeric rather than symbolic. The numerical procedure presented here, however, may nevertheless inform subsequent symbolic analysis.
\end{itemize}

\subsection{Implementing Proper and Improper Translations}
\label{propertranslationssection}

In this section, we derive the necessary logical relations to guarantee that the conditions of Definition \ref{translation} are satisfied. We introduce decision variables $H[i,j'] \in \left\{ 0, 1 \right\}$, $i = 1, \ldots, q$, $j' = 1, \ldots, \tilde{m}$, and $\tilde{b}[i',j'] \geq 0,$ $i',j' = 1 ,\ldots, \tilde{m}$, $i' \not= j'$, so that
\begin{equation}
\label{equation1}
\left\{
\begin{array}{ll}
H[i,j'] = 1, \; \; \; \; \; \; \; \; & \mbox{ if and only if } h(i)=j'.\\
\tilde{b}[i',j'] > 0, & \mbox{ if and only if } (i',j') \in \tilde{\mathcal{R}} \\
\end{array}
\right.
\end{equation}

\noindent We can accommodate (\ref{equation1}) and (\ref{equation2}) with the constraint set (\ref{constraint1}) where $\tilde{A}(\tilde{\mathcal{B}})$ is the matrix with off-diagonal entries $[\tilde{A}(\tilde{\mathcal{B}})]_{i',j'} = \tilde{b}[j',i']$, $i' \not= j'$, and $H \in \mathbb{Z}^{q \times \tilde{m}}$ is the matrix with entries $H_{i,j'} = H[i,j']$. We can further restrict to proper translations by imposing the constraint set (\ref{constraint2}).

In order to satisfy properties $1(a-c)$ of Definition \ref{translation}, we introduce variables $\lambda[i,j'] \geq 0$, $i=1,\ldots, q$, $j' = 1, \ldots, \tilde{m}$, such that
\begin{equation}
\label{2344}
\begin{array}{ll}
\lambda[i,j'] > 0, \; \; \; \; \; \; \; \; & \mbox{if and only if } \lambda(i,j') > 0.\\
\end{array}
\end{equation}
We can accommodate (\ref{2344}) with the constraint set (\ref{constraint}) where $\Lambda \in \mathbb{R}^{\tilde{m} \times q}$ is the matrix with entries $\Lambda_{j',i} = \lambda[i,j']$. Notice that this constraint set is only distinct from (\ref{constraint1}) if we are allowing improper translations. Consquently, if we are interested only in proper translations, we use (\ref{constraint1}) and (\ref{constraint2}), and if we are interested in improper translations (or do not care which is attained) we use (\ref{constraint1}) and (\ref{constraint}).

\subsection{Implementing Weak Reversibility}
\label{deficiencysection}

In this section, we reiterate the results of \cite{J-S4} and \cite{J-S6}, respectively, for guaranteeing that the translation is weak reversibility and that it has the minimal structural deficiency.

In order to guarantee $\tilde{\mathcal{N}}(\tilde{\mathcal{B}})$ is weakly reversible, we introduce decision variables $\tilde{w}[i',j'] \geq 0$, $i',j' = 1, \ldots, \tilde{m},$ $i' \not= j'$, so that
\begin{equation}
\label{equation12121}
\left\{
\begin{array}{ll}
\tilde{w}[i',j'] > 0, \; \; \; \; \; \; \; \; & \mbox{ if and only if } (i',j') \in \tilde{\mathcal{R}}\\
\mathbf{1} \cdot \tilde{A}(\tilde{\mathcal{W}}) = \mathbf{0}, & \mbox{ for } \mathbf{1} = (1,\ldots,1), \; \mathbf{0} = (0,\ldots,0)\\
\tilde{A}(\tilde{\mathcal{W}}) \cdot \mathbf{1} = \mathbf{0}, &
\end{array}
\right.
\end{equation}
where $\tilde{A}(\tilde{\mathcal{W}})$ is the matrix with off-diagonal entries $[\tilde{A}(\tilde{\mathcal{W}})]_{i',j'} = \tilde{w}[j',i']$, $i' \not= j'$. The matrix $\tilde{A}(\tilde{\mathcal{W}})$ has the same structure as $\tilde{A}(\tilde{\mathcal{B}})$ but has been scaled along its columns (for details, see \cite{J-S4}). The logical requirements (\ref{equation12121}) can be accommodated by the constraint set (\ref{weaklyreversible}).

We now introduce decision variables capable of calculating the deficiency of a chemical reaction network. It was observed in \cite{J-S6} that $m$ and $s$ are fixed prior to the optimization begin, so that to determine the deficiency it suffices to calculate the number of linkage classes. It also follows by the well-known property $\delta = m - s - \ell \geq 0$ that $\ell \leq m - s$. Since we have $\tilde{s}=s$ for weakly reversible network translations by Lemma \ref{lemma232}, it is sufficient to allow at most $\tilde{\ell} = \tilde{m}-s$ linkage classes. Following \cite{J-S6}, we introduce decision variables $\gamma[i',\theta] \in \left\{ 0, 1 \right\}$, $i' = 1, \ldots, \tilde{m}$, $\theta=1, \ldots, \tilde{m}-s$, and $\tilde{L}[\theta] \in [0, 1]$, $\theta=1, \ldots, 
\tilde{m}-s$, so that
\begin{equation}
\label{5103}
\left\{
\begin{array}{ll}
\gamma[i',\theta] = 1, \; \; \; \; \; \; \; \; & \mbox{ if and only if } i' \in \tilde{L}_{\theta} \\
\tilde{L}[\theta] = 1, & \mbox{ if and only if } \tilde{L}_{\theta} \not= \emptyset \\
\tilde{w}[i',j'] > 0, & \mbox{ for } i' \not= j' \mbox{ implies } i',j' \in \tilde{L}_{\theta} \mbox{ for some } \theta \in \{ 1, \ldots, \tilde{m}-s \}
\end{array}
\right.
\end{equation}
where $\tilde{L}_{\theta}$, $\theta=1, \ldots, \tilde{m}-s$, are the (potential) linkage classes of $\tilde{\mathcal{N}}(\tilde{\mathcal{B}})$.

The variables $\gamma[i',\theta]$ keep track of which complexes are assigned to which linkage class while the variables $\tilde{L}[\theta]$ keep track of whether a particular linkage classes has complexes in it. It is worth noting that unused complexes in the potential complex set are assigned to their own isolated linkage classes. This is a slight abuse of chemical reaction network convention but will be allowed in the present context. It was noted in \cite{J-S6} that allowing isolated linkage classes does not alter the network property of weak reversibility or the value of the deficiency. The final requirement of (\ref{5103}) guarantees that no reaction may proceed between complexes on different linkage classes.

It was also noted in \cite{J-S6} that the assignment of complexes to linkage classes is not unique since any permutation of the assignment of linkage classes corresponds to the same network. This can be a significant problem for the efficiency of mixed integer programming methods. We therefore require that partition structure, if it can be found, is unique. This uniqueness requirement and (\ref{5103}) can be accommodated with the constraint set (\ref{partition}). (See \cite{J-S6} for a rigorous justification of these constraints.)


We may now find the weakly reversible reaction-weighted translated chemical reaction network with the underlying reaction network with the minimal deficiency by optimizing (\ref{deficiency}) over the constraint (\ref{constraint1}), (\ref{constraint}), (\ref{weaklyreversible}), and (\ref{partition}). If we are only interested in proper reaction-weighted translations, we may substitute the constraint set (\ref{constraint2}) in place of (\ref{constraint}).

\subsection{Implementing steady state Resolvability}
\label{implementingsection}


In this section, we develop constraint sets which guarantee that the conditions of Theorem \ref{maintheorem} are satisfied for improper translations. We will divide this into the three steps.\\

\noindent \emph{Step 1: Determine constants $c(i',j')$ consistent with Definitions \ref{cs} and \ref{resolving2}:} It will not be necessary to assign decision variables to track $\tilde{\mathcal{C}}_I$ and $\tilde{\mathcal{C}}_R$ specifically. We will instead build conditions which will accurately determine the constants $c(i',j')$ in Definition \ref{resolving2}. Note first of all, however, that the complex vectors relevant to condition $1.$ of Definition \ref{resolving2} are found in the matrix $Y$ rather than $\tilde{Y}$. We therefore define the variables $c[i,j] \geq 0$, $i,j = 1, \ldots, q$, $i \not= j$, and require that
\[c[i,j] > 0 \mbox{ if and only if } h(i) = i', h(j) = j', \mbox{ and either } c(i',j') \not= 0 \mbox{ or } c(j',i') \not= 0.\]
To track the improper and resolving complex sets, $\tilde{\mathcal{C}}_I$ and $\tilde{\mathcal{C}}_R$, we introduce the variables
\begin{equation}
\label{9293}
\left\{
\begin{array}{ll}
\delta_I[i,j] \in \left\{ 0, 1 \right\}, \; \; \; & i,j =1 ,\ldots, q, i < j\\
\delta_K[i,j] \in \left\{ 0, 1 \right\}, \; \; \; \; \; \; \; \; \; \; \; \; \;& i,j = 1, \ldots, q, i < j\\
\gamma_K[i,\theta] \in \left\{ 0, 1 \right\}, \; & i =1, \ldots, q, \theta=1, \ldots, \tilde{m}-s.
\end{array} 
\right.
\end{equation}
The variables $\delta_I[i,j]$ track the supports of the complexes in $Y$ which are mapped through $h$ to $\tilde{\mathcal{C}}_I$ (left-hand sides of condition 1. of Definition \ref{resolving2}) while variables $\delta_R[i,j]$ track the supports of the complexes in $Y$ which are mapped through $h$ to $\tilde{\mathcal{C}}_R$ (right-hand sides of condition 1. of Definition \ref{resolving2}). The variables $\gamma_K[i,\theta]$ correspond the linkage classes in $\tilde{\mathcal{N}}(\tilde{\mathcal{B}})$ to the supports of the complexes in $Y$ so that condition $2.$ of Definition \ref{resolving2} may be imposed.

In order to limit the number of variables in the system, we attempt to satisfy condition $1.$ of Definition \ref{resolving2} simultaneously over all pairs $i,j \in h^{-1}(k')$ where $k' \in \tilde{\mathcal{C}}_I$. We introduce a stochastic parameter $v[i,j] \in [\sqrt{\epsilon},1/\sqrt{\epsilon}]$, $i,j = 1, \ldots, q$, $i < j$, and consider the $v[i,j]$-weighted sum of the conditions in condition $1.$ of Definition \ref{resolving2}. The introduction of the parameter stochastic parameter $v[i,j]$ makes it almost certain that linear dependence does not become an issue when summing over the left-hand sides of condition $1$ of Definition \ref{resolving2}. The parameters are chosen over the range $[ \sqrt{\epsilon},1/\sqrt{\epsilon}]$ rather than the more natural $[\epsilon,1/\epsilon]$ for numerical stability.



In order to satisfy the requirements of Definition \ref{cs} and \ref{resolving2}, we require the following logical relations:
\begin{equation}
\label{crazy}
\left\{
\begin{array}{ll}
\delta_I[i,j] = 1, & \mbox{ if and only if } h(i)=k' \mbox{ and } h(j)=k' \mbox{ for some } k' \in \tilde{\mathcal{C}}_I \\
\delta_K[i,j] = 1, & \mbox{ if and only if } c(i,j) > 0 \mbox{ or } c(j,i) > 0\\
\gamma_K[i,\theta] = 1, & \mbox{ if and only if } h(i)=k' \mbox{ and } k' \in \tilde{L}_\theta
\end{array}
\right.
\end{equation}

\noindent We can accommodate the requirements of (\ref{crazy}) with the constraint set (\ref{constraint3}).\\ 

\noindent \emph{Step 2: (conditions $(1-2)$ of Theorem \ref{maintheorem}):} We introduce the decision variables $\tilde{C}^*[i'] \in \{ 0, 1\}$, $i' = 1, \ldots, \tilde{m}$, and $\tilde{b}^*[i',j'] \geq 0$, $i', j' =1, \ldots, \tilde{m}$, $i' \not= j'$, and impose
\begin{equation}
\label{rewq}
\left\{
\begin{array}{ll}
\tilde{C}^*[i'] = 1, & \mbox{if and only if } i' \in \tilde{\mathcal{C}}^*\\
\tilde{b}^*[i',j'] > 0,  \; \; \; \; \; \; \; \; \; \; \; & \mbox{if and only if } (i',j') \in \tilde{\mathcal{R}}^*.\\
\end{array}
\right.
\end{equation}
We want $\tilde{\mathcal{C}}^*$ to restrict the supports of $\tilde{\mathcal{C}}_I$ and $\tilde{\mathcal{C}}_R$ according to condition 1. of Theorem \ref{maintheorem}. We also want $\tilde{\mathcal{C}}^*$ and the reaction network $\tilde{\mathcal{R}}$ to restrict $\tilde{\mathcal{R}}^*$ according to condition 2. of Theorem \ref{maintheorem}. We can accomplish this with the constraint set (\ref{constraint4}).\\

\noindent \emph{Step 3: (conditions $(3-4)$ of Theorem \ref{maintheorem}):} We introduce the decision variables
\begin{equation}
\label{93284972}
\left\{
\begin{array}{ll}
\tilde{C}^{**}[i'] \in \{ 0, 1\}, & i' = 1, \ldots, \tilde{m}\\
\tilde{b}^{**}[i',j'] \geq 0, & i',j' =1, \ldots, \tilde{m}, i'\not=j'\\
\gamma^*[i',\theta] \geq 0, & i'=1,\ldots, \tilde{m}, \theta = 1, \ldots, \tilde{\ell}^*\\
\tilde{L}^*[\theta] \in \{ 0, 1 \}, & \theta = 1, \ldots, \tilde{\ell}^*
\end{array}
\right.
\end{equation}
where $\tilde{\ell}^*$ is a predetermined upper limit on the number of linkage classes of $(\tilde{\mathcal{S}},\tilde{\mathcal{C}}^* \cup \tilde{\mathcal{C}}^{**},\tilde{\mathcal{R}}^* \cup \tilde{\mathcal{R}}^{**})$. Note that this may be strictly larger than $\tilde{m} - s$. We impose that
\begin{equation}
\label{asdf}
\left\{
\begin{array}{ll}
\tilde{C}^{**}[i'] =1, & \mbox{ if and only if }  i' \in \tilde{\mathcal{C}}^{**}\\
\tilde{b}^{**}[i',j'] > 0, & \mbox{ if and only if } (i',j') \in \tilde{\mathcal{R}}^{**}\\
\gamma^*[i',\theta] = 1, \; \; \; \; \; \; \; \; & \mbox{ if and only if } i' \in \tilde{L}_{\theta}^* \\
\tilde{L}^*[\theta] = 1, & \mbox{ if and only if } \tilde{L}_{\theta}^* \not= \emptyset
\end{array}
\right.
\end{equation}
where $\tilde{L}_{\theta}^*$ is a linkage class of $(\tilde{\mathcal{S}},\tilde{\mathcal{C}}^* \cup \tilde{\mathcal{C}}^{**}, \tilde{\mathcal{R}}^* \cup \tilde{\mathcal{R}}^{**})$. That is, the variables $\tilde{C}^{**}[i']$ keep track of the complexes in $\tilde{\mathcal{C}}^{**}$ while the variables $\tilde{b}^{**}[i',j']$ keep track of the structure of $\tilde{\mathcal{R}}^{**}$. The variables $\gamma^*[i',\theta]$ and $\tilde{L}^*[\theta]$ keep track of the linkage classes of $(\tilde{\mathcal{S}},\tilde{\mathcal{C}}^* \cup \tilde{\mathcal{C}}^{**},\tilde{\mathcal{R}}^* \cup \tilde{\mathcal{R}}^{**})$ (see Section \ref{deficiencysection}). We can accommodate the requirements of conditions (3-4) of Theorem \ref{maintheorem} with the constraint set (\ref{constraint5}). In order to limit the size of the components in $(\tilde{\mathcal{S}},\tilde{\mathcal{C}}^* \cup \tilde{\mathcal{C}}^{**},\tilde{\mathcal{R}}^* \cup \tilde{\mathcal{R}}^{**})$, we can additionally optimize over (\ref{optimal}).

\section{Applications}
\label{applications}

In this section, we apply the computational methodology of Section \ref{milpsection} to two examples drawn from the mathematical biology literature.

The first network was considered earlier as Example 2 in Section \ref{resolvabilitysection}. The model was original introduced as a candidate EnvZ/OmpR signaling pathway mechanism in \emph{escherichia coli} by Shinar and Feinberg in the Supporting Online Material of \cite{Sh-F}. The model was shown to be steady state equivalent to a generalized reaction network in \cite{J1}. The second network is modified from a model of the PFK-2/FBPase-2 mechanism in mammals which was originally presented by Dasgupta \emph{et al.} in \cite{Dasgupta,Karp}. The application of network translation to this model is novel. All computations were performed with the GNU Linear Programming Kit (GLPK) on the author's personal use Toshiba Satellite laptop (AMD Quad-Core A6-Series APU, 6GB RAM). Full details of the computations are contained in the Supplemental Material.

\subsection{Application I: EnvZ-OmpR Mechanism}


Reconsider the mechanism given in Example 2 in Section \ref{resolvabilitysection}. We now apply the computational process presented in Section \ref{milpsection}. The details of the initialization are contained in the the Supplementary Material. We note here, however, that we have initialized the rate constants stochastically within the range $k_i \in [\sqrt{\epsilon},1/\sqrt{\epsilon}]$, $i=1, \ldots, 14$, rather than chosing them to be fixed constants. The code was run 25 times, with an average time to completion of $2.788$ seconds and a standard deviation of $1.4898$ seconds. In each case, the algorithm successfully found the weakly reversible network translation given in Figure \ref{figure1}(b).

This is consistent with the translation obtained in \cite{J1}. To further check the consistency of the code, we observe that it returned the sets $\tilde{\mathcal{C}}_I = \{ 6 \}$, $\tilde{\mathcal{C}}_R = \{ 1, 3\}$, $\tilde{\mathcal{C}}(6) = \{ 6, 7, 8 \}$, $\tilde{\mathcal{R}}^* = \{ (6,7),(6,8),(7,2),(7,6),(8,2),\\(8,6) \}$, $\tilde{\mathcal{C}}^{**} = \{ 2 \}$, and $\tilde{\mathcal{R}}^{**} = \{ (2,6) \}$. This is consistent with the application of Theorem \ref{maintheorem} to the reaction-weighted translation $\tilde{\mathcal{N}}(\tilde{\mathcal{B}})$ (see Figure \ref{figure1}(c)). Since the network has $\tilde{\delta}=0$, it follows by Theorem \ref{maintheorem} that the network is steady state resolvable. (Further methodology for characterizing the steady state set is contained in the Supplemental Material and in \cite{J1}.)

\subsection{Application II: PFK-2/FBPase-2 Mechanism}

\begin{figure}[h]
\centering
\includegraphics[width=12cm]{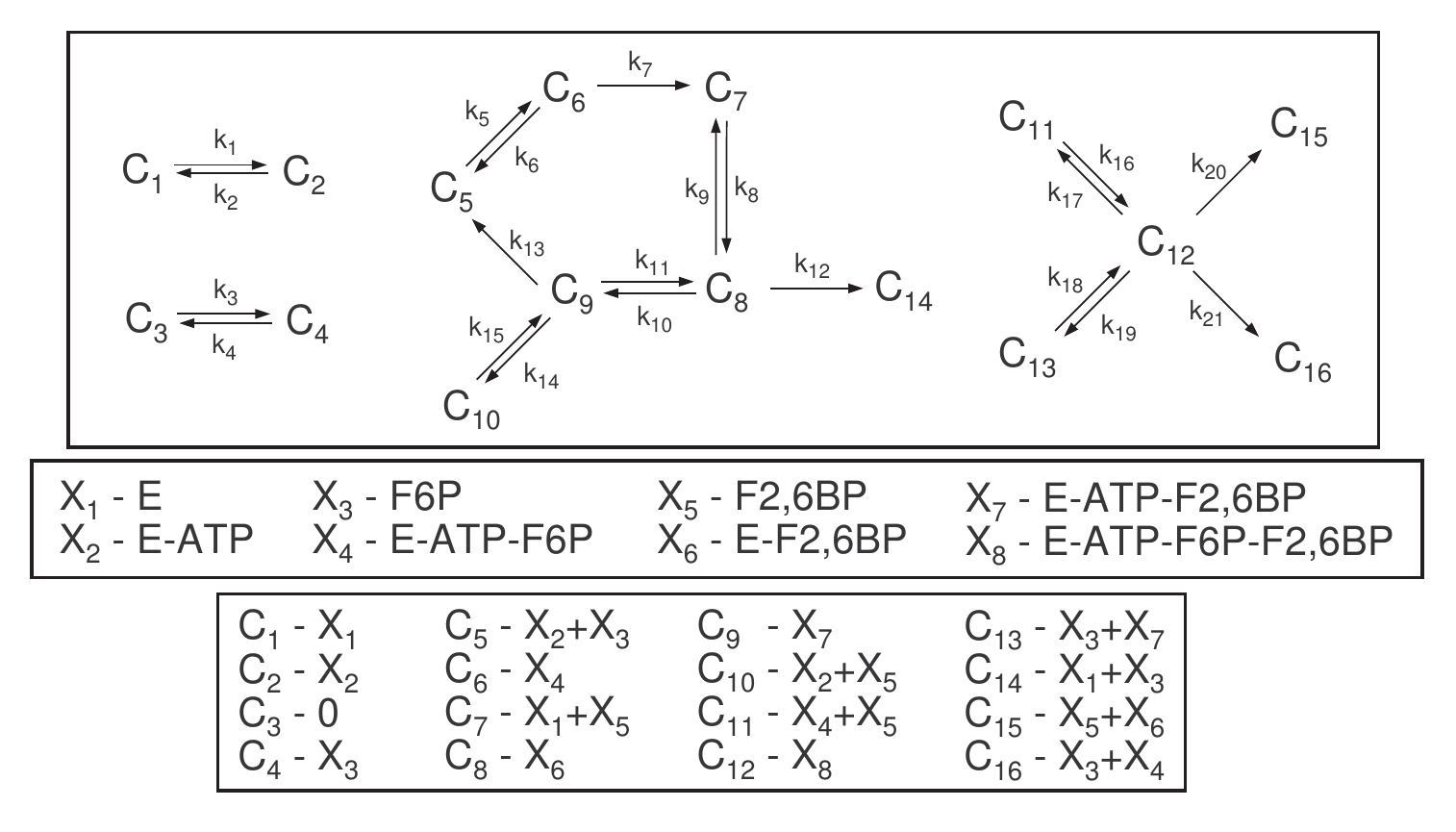}
\caption{Candidate PFK-2/FBPase-2 mechanism in mammalian cells.}
\label{figure3}
\end{figure}

Consider the following hypothetical PFK-2/FBPase-2 mechanism contained in Figure \ref{figure3}. This model is based on one proposed in \cite{Dasgupta,Karp} but differs in the reversible reaction pair $C_3 \rightleftarrows C_4$ which corresponds to $\emptyset \rightleftarrows X_3$. Our mechanism therefore allows for inflow and outflow of Fructose 6-phosphate ($F6P$). We defer biochemical justification and analysis of this mechanism to \cite{Dasgupta,Karp}.

We now apply the computational algorithm of Section \ref{milpsection}. We first simplify the model by assuming that $k_{19}=k_{20}$ and initializing the remaining reaction-weights stochastically from the range $k_i \in [\sqrt{\epsilon}, 1/\sqrt{\epsilon}]$, $i=1,\ldots, 20$. The code was run successfully 25 times with an average completion time of $6.604$ seconds seconds and a standard deviation of $2.6871$ seconds.

\begin{figure}[h!]
\centering
\includegraphics[width=12cm]{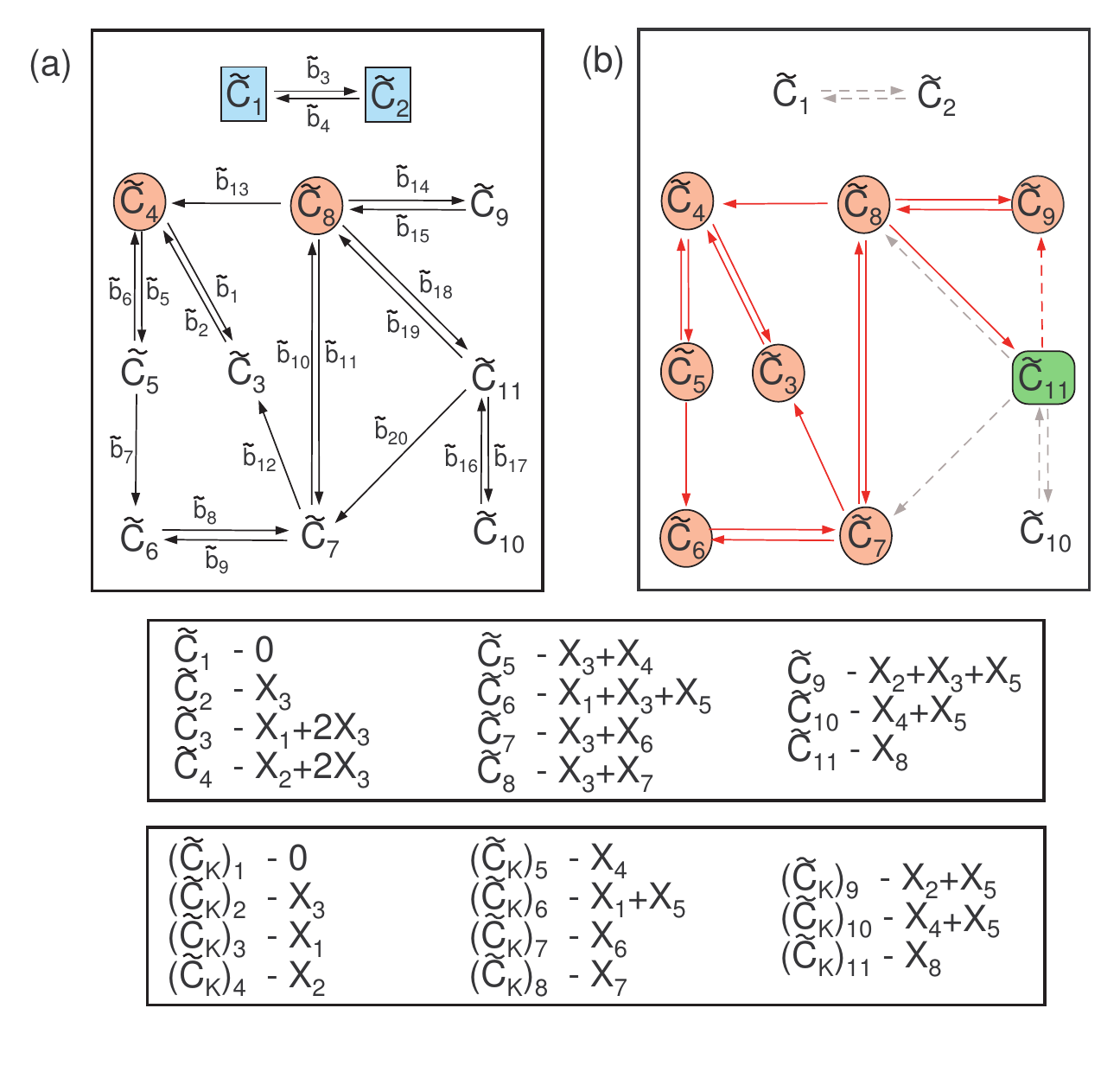}
\caption{In (a), we have the computationally-determined reaction-weighted translation $\tilde{\mathcal{N}}(\tilde{\mathcal{K}})$ for the PFK-2/FBPase-2 mechanism contained in Figure \ref{figure3}. In (b), we identify  $\tilde{\mathcal{C}}^* = \{3,4,5,6,7,8,9\}$ (pink), $\tilde{\mathcal{C}}^{**} = \{ 11 \}$ (green), and the reaction sets $\tilde{\mathcal{R}}^*$ (solid red arrows) and $\tilde{\mathcal{R}}^{**} = \{ (11,9) \}$ (dashed red arrow). The network $(\tilde{\mathcal{S}},\tilde{\mathcal{C}}^* \cup \tilde{\mathcal{C}}^{**}, \tilde{\mathcal{R}}^* \cup \tilde{\mathcal{R}}^{**})$ is clearly weakly reversible.}
\label{figure4}
\end{figure}

A recurring network structure for the translation was the one contained in Figure \ref{figure4}(a). Note that both $C_2 = X_2$ and $C_5 = X_2 + X_3$ are translated to $\tilde{C}_4 = X_2 + 2X_3$, and both $C_9 = X_7$ and $C_{13} = X_3 + X_7$ are translated to $\tilde{C}_8 = X_3 + X_7$. The reaction-weighted translation is therefore improper. The algorithm returned the sets $\tilde{\mathcal{C}}_I = \{ 4,8 \}$, $\tilde{\mathcal{C}}_R = \{ 1, 2\}$, $\tilde{\mathcal{C}}^* = \{ 3, 4, 5, 6, 7, 8, 9 \}$, $\tilde{\mathcal{R}}^* = \{ (3,4),(4,3), (4,5), (5,4),(5,6),(6,7), (8,4), (8,7),\\ (8,9), (8,11), (9,8) \}$, $\tilde{\mathcal{C}}^{**} = \{ 11 \}$, and $\tilde{\mathcal{R}}^{**} = \{ (11,9) \}$. Notice that, even though the technical conditions of Theorem \ref{maintheorem} are satisfied trivially  (see Figure \ref{figure4}(b)), the algorithm still constructs a weakly reversible component containing $\tilde{\mathcal{C}}_I$. Details of the computation are contained in the Supplementary Material.

Notice that we may not apply Theorem \ref{maintheorem} directly since the network has $\tilde{\delta} = 2$. Nevertheless, it can shown that ker$(\tilde{Y} \cdot A(\tilde{\mathcal{K}}))$ decomposes in such a way steady state equivalence may be guaranteed (see Supplementary Material). The generalized mass action system with the rate constants given in Table \ref{table1} has the same steady states as the original system. Note that, although the computational process requires numerical values for the reaction-weights, the insight gained from the process was able to guide a correspondence which can be shown to work for all reaction-weights.

\begin{table}[h!]
\caption{Reaction-weights for steady state equivalence of $\mathcal{N}(\mathcal{K})$ and $\tilde{\mathcal{N}}(\tilde{\mathcal{K}})$.}
\label{table1}
\begin{center}
\begin{tabular}[h]{|l|l|l|l|}
\hline
$\tilde{k}_1 = k_1$ & $\tilde{k}_{6} = k_6$ & $\tilde{k}_{11} = k_{11}$ & $\tilde{k}_{16} = k_{16}$ \\
$\tilde{k}_2 = k_2$ & $\tilde{k}_7 = k_7$ & $\tilde{k}_{12} = k_{12}$ & $\tilde{k}_{17} = k_{17} + k_{20}$ \\
$\tilde{k}_3 =k_3$ & $\tilde{k}_8 = k_8$ & $\tilde{k}_{13} = k_{13}$ & $\tilde{k}_{18} = \frac{k_3}{k_4} k_{18}$ \\
$\tilde{k}_4 = k_4$ & $\tilde{k}_9 = k_9$ & $\tilde{k}_{14} = k_{14}$ & $\tilde{k}_{19} = k_{19}$ \\
$\tilde{k}_5 = \frac{k_3}{k_4} k_5 $ & $\tilde{k}_{10} = k_{10}$ & $\tilde{k}_{15} = k_{15}$ & $\tilde{k}_{20} = k_{20}$ \\
\hline
\end{tabular}
\end{center}
\end{table}

\section{Conclusions and Future Work}
\label{conclusionsection}

In this paper, we have extended the theory of network translation \cite{J1} in two important ways:
\begin{itemize}
\item[(Q1)]
We have presented conditions which suffice to guarantee steady state resolvability of a reaction-weighted network $\mathcal{N}(\mathcal{K})$ and a reaction-weighted translation $\tilde{\mathcal{N}}(\tilde{\mathcal{B}})$ (Theorem \ref{maintheorem}). Importantly, these conditions are graph theoretic in nature and do not require an enumeration over all cycles on the support of the elementary flux modes as was previously required by \cite{J1}.
\item[(Q2)]
We have presented an algorithm for determining whether a reaction-weighted translation of a given chemical reaction network exists. This algorithm is implementable within the well-known MILP framework and is capable of imposing the technical conditions of Theorem \ref{maintheorem}. The code is contained in Appendix \hyperlink{appendixc}{C}.
\end{itemize}

There are numerous avenues open for future work in the study of network translations, and generalized mass action systems in general. The avenues specifically related to the work contained in this paper include:

\begin{enumerate}
\item
\emph{Algorithmic determination of optimal $\tilde{Y}$}: The MILP algorithm presented in Section \ref{milpsection} requires initialization of the matrix $\tilde{Y}$ consisting of potential stoichiometric complexes in the network $\tilde{\mathcal{N}}$. Without prior intuition, a suitable choice of these complexes may not be obvious. Nevertheless, this choice set should be kept small to maintain computation efficiency. Algorithmically determining a suitable set of potential stoichiometric complexes is therefore a primary concern moving forward.
\item
\emph{Simplifying constraint sets}: When not carefully posed, the algorithm presented in Section \ref{milpsection} may take significant time to complete. Numerical stability is also an issue for small values of $\epsilon$. While this is not unexpected as MILP optimization problems are known to be NP-hard, it is nevertheless an important task to simplify the code, and the conditions underlying resolvability, in order to make the algorithm computationally tractable for larger problems.
\item
\emph{Expansion of underlying theory}: The main result of this paper (Theorem \ref{maintheorem}) depends implicitly on the results regarding translations contained in \cite{J1} and those regarding generalized mass action systems contained in \cite{M-R}. It is anticipated that, as these nascent theories are further developed that the applications of computational approaches such as those contained in this paper will become necessary. We present in the Supplemental Material an example which illustrates one further avenue of research regarding the theory of network translation.
\end{enumerate}

\section*{Acknowledgments}

The author acknowledges helpful comments from Gabor Szed\'{e}rkenyi and David F. Anderson who suggested several significant improvements to the early drafts of this paper.

\section*{\hypertarget{appendixa}{Appendix A} (Resolvability)}


In order to make the connection between the results of \cite{J1}, Definition \ref{resolvable}, and Theorem \ref{maintheorem}, we briefly introduce here some background on resolvability. We begin by defining the following concept, which was introduced informally in Section \ref{backgroundsection}. 
\begin{definition}
Suppose $\mathcal{N} = (\mathcal{S},\mathcal{C},\mathcal{R})$ is a chemical reaction network. We say a subgraph $P = \{ \mathcal{C}_P,\mathcal{R}_P \}$ where $\mathcal{C}_P \subseteq \mathcal{C}$ and $\mathcal{R}_{P} \subseteq \mathcal{R}$ is a \textbf{path} from $C_i$ to $C_j$ if:
\begin{enumerate}
\item
there is an ordering $\{ \nu(1), \nu(2), \ldots, \nu(l) \}$ with all $\nu(i)$, $i=1,\ldots, l,$ distinct such that $C_i = C_{\nu(1)} \rightarrow C_{\nu(2)} \rightarrow \cdots \rightarrow C_{\nu(l)} = C_j$;
\item
$\mathcal{C}_P = \left\{ \nu(1), \nu(2), \ldots, \nu(l)\right\} \subseteq \mathcal{C}$; and
\item
$\mathcal{R}_P = \left\{  \big(\nu(1),\nu(2)\big), \ldots, \big(\nu(l-1),\nu(l) \big) \right\} \subseteq \mathcal{R}$.
\end{enumerate}
We will let $\mathcal{P}(i,j)$ denote the set of all paths from $C_i$ to $C_j$.
\end{definition}

Now consider the following.

\begin{definition}
\label{tree}
Suppose $\mathcal{N} =(\mathcal{S},\mathcal{C},\mathcal{R})$ is a chemical reaction network. We say a subgraph $T = \{\mathcal{C}_{T},\mathcal{R}_{T}\}$ where $\mathcal{C}_{T} \subseteq \mathcal{C}$ and $\mathcal{R}_{T} \subseteq \mathcal{R}$ is a \textbf{spanning i-tree} on $\mathcal{C}_T$ if:
\begin{enumerate}
\item
$\mathcal{R}_T$ spans $\mathcal{C}_{T}$;
\item
$T$ contains no directed or undirected cycles; and
\item
$T$ has a unique sink at $C_i$.
\end{enumerate}
We will let $\mathcal{T}(i)$ denote the set of all spanning $i$-trees on $\mathcal{C}_T$.
\end{definition}

In general, an arbitrary subset $\mathcal{C}_T \subseteq \mathcal{C}$ may not permit any spanning $i$-trees; however, if the network is weakly reversible, there is at least one spanning $i$-tree on the set $\mathcal{C}_T = L_{\theta}$ where $L_{\theta}$ is the linkage class which contains $C_i$. These are the components to which we will be interested in restricting. We may define the following for weakly reversible networks.
\begin{definition}
Consider a reaction-weighted chemical reaction network $\mathcal{N}(\mathcal{K}) = (\mathcal{S},\mathcal{C},\mathcal{R},\mathcal{K})$ which is weakly reversible. Then the \textbf{tree constant} for $i \in \{ 1, \ldots, m \}$ is given by
\begin{equation}
\label{treeconstant}
K_{i} = \sum_{T \in \mathcal{T}(i)} \prod_{(i^*,j^*) \in T} k(i^*,j^*)
\end{equation}
where $\mathcal{T}(i)$ is the set of spanning $i$-trees on the component $\mathcal{C}_T = L_{\theta}$ where $C_i \in L_{\theta}$.
\end{definition}

For example, for the network
\[\begin{array}{c}  \vspace{-0.1cm}C_1  \displaystyle{\mathop{\stackrel{k_1}{\begin{array}{c} \vspace{-0.25cm} \leftarrow \vspace{-0.3cm} \\ \rightarrow \end{array}}}_{k_2}} C_2 \\ {}_{k_5} \hspace{-0.15cm} \uparrow \; \; \; \; \; \; \; \; \downarrow_{k_3} \\ C_4 \; \stackrel{k_4}{\leftarrow} \; C_3 \end{array}\]
we have
\[K_1 = k_2k_4k_5+k_3k_4k_5\]
corresponding to the two spanning trees with unique sinks at $C_1$:
\[\begin{array}{c}  \vspace{-0.1cm}C_1  \displaystyle{\mathop{\stackrel{k_1}{\begin{array}{c} \vspace{-0.25cm} {\color{red}\leftarrow} \vspace{-0.3cm} \\ \rightarrow \end{array}}}_{{\color{red}k_2}}} C_2 \\ {}_{{\color{red}k_5}} \hspace{-0.05cm} {\color{red}\uparrow} \; \; \; \; \; \; \; \; \downarrow_{k_3} \\ C_4 \; \; {\color{red}\stackrel{k_4}{\leftarrow}} \; \; C_3 \end{array} \hspace{0.5cm} \mbox{and} \hspace{0.5cm} \begin{array}{c}  \vspace{-0.1cm}C_1  \displaystyle{\mathop{\stackrel{k_1}{\begin{array}{c} \vspace{-0.25cm} \leftarrow \vspace{-0.3cm} \\ \rightarrow \end{array}}}_{k_2}} C_2 \\ {\color{red}{}_{k_5} \hspace{-0.15cm} \uparrow} \; \; \; \; \; \; \; \; \; {\color{red}\downarrow_{k_3}} \\ C_4 \; \; {\color{red}\stackrel{k_4}{\leftarrow}} \; \; C_3. \end{array}\]

\begin{remark}
An immediate consequence of Definition \ref{tree} is that, given an $i$-tree which spans $\mathcal{C}_T$, there is a unique path from every $C_j \in \mathcal{C}_T$ to $C_i$. We will make use of this fact in the proofs contained in Appendix \hyperlink{appendixb}{B}.
\end{remark}

\begin{remark}
We will denote the tree constants of the translated reaction-weighted networks $\tilde{\mathcal{N}}(\tilde{\mathcal{B}}) = (\tilde{\mathcal{S}},\tilde{\mathcal{C}},\tilde{\mathcal{C}}_K,\tilde{\mathcal{R}},\tilde{\mathcal{B}})$ as $\tilde{B}_{i'}$, $i' = 1, \ldots, \tilde{m}$. Note also that the convention of referring to these algebraic constructs as ``tree constants'' is original to \cite{J1}.
\end{remark}

\section*{\hypertarget{appendixb}{Appendix B} (Proof of Theorem \ref{maintheorem})}

Before proving Theorem \ref{maintheorem}, we present the following equivalent result. The result may be more intuitive to many readers.

\begin{lemma}
\label{mainlemma}
Let $\tilde{\mathcal{N}}(\tilde{\mathcal{B}}) = (\tilde{\mathcal{S}},\tilde{\mathcal{C}},\tilde{\mathcal{C}}_K,\tilde{\mathcal{R}},\tilde{\mathcal{B}})$ denote an improper reaction-weighted translation of a reaction-weighted chemical reaction network $\mathcal{N}(\mathcal{K}) = (\mathcal{S},\mathcal{C},\mathcal{R},\mathcal{K})$. Suppose that $\tilde{\mathcal{N}}$ is weakly reversible, $\tilde{\delta} = 0$, $\tilde{S}_I \subseteq \tilde{S}_K$, and that there is a resolving complex set $\tilde{\mathcal{C}}_R$ satisfying $\tilde{\mathcal{C}}_I \cap \tilde{\mathcal{C}}_R = \emptyset$, where $\tilde{\mathcal{C}}_I$ is the improper complex set of $\tilde{\mathcal{N}}(\tilde{\mathcal{B}})$. Suppose furthermore that $\tilde{\mathcal{C}}_I$ and $\tilde{\mathcal{C}}_R$ satisfy the following property:
\begin{itemize}
\item[($*$)] If $p' \in \tilde{\mathcal{C}}_I$, then there is a $k' \in \tilde{\mathcal{C}}(\tilde{\mathcal{B}})$, $k' \not= p'$, such that, if $i' \in \tilde{\mathcal{C}}_R$ and $\tilde{P} \in \tilde{\mathcal{P}}(p',i')$, then $k' \in \tilde{\mathcal{C}}_{\tilde{P}}$.
\end{itemize}
Then $\mathcal{N}(\mathcal{K})$ and $\tilde{\mathcal{N}}(\tilde{\mathcal{B}})$ are steady state resolvable.
\end{lemma}

\begin{remark}
This results says that, given the technical requirement ($*$), the translations is resolvable if, for every improper complex there is a common complex such that every path from the improper complex to a resolving complex goes through the common complex. It is worth noting similarities in condition ($*$) and those of conditions (14-16) of \cite{Boros2013}, although no deeper connection is known to the author at present.
\end{remark}


\begin{proof}[Proof of Lemma \ref{mainlemma}]
Suppose that $\tilde{\mathcal{N}}(\tilde{\mathcal{B}}) = (\tilde{\mathcal{S}},\tilde{\mathcal{C}},\tilde{\mathcal{C}}_K,\tilde{\mathcal{R}},\tilde{\mathcal{B}})$ is an improper reaction-weighted translation of a reaction-weighted chemical reaction network $\mathcal{N}(\mathcal{K}) = (\mathcal{S},\mathcal{C},\mathcal{R},\mathcal{K})$ according to Definition \ref{translation}. Suppose furthermore that $\tilde{\mathcal{N}}(\tilde{\mathcal{B}})$ is weakly reversible, $\tilde{\delta} = 0$, and $\tilde{S}_I \subseteq \tilde{S}_K$. 


Since $\tilde{S}_I \subseteq \tilde{S}_K$, there is a non-empty resolving complex set $\tilde{\mathcal{C}}_R$ according to Definition \ref{resolving2}. Let $\tilde{B}_{i'}$, $i' = 1, \ldots, \tilde{m}$, denote the tree constants (\ref{treeconstant}) corresponding to the reaction-weighted reaction graph of $\tilde{\mathcal{N}}(\tilde{\mathcal{B}})$. 
Consider any pair $i,j \in h^{-1}(k')$ where $k' \in \tilde{\mathcal{C}}_I$, and define the ratios
\begin{equation}
\label{weird3}
\tilde{B}_{i,j} = \prod_{i',j'=1}^{\tilde{m}} \left( \frac{\tilde{B}_{j'}}{\tilde{B}_{i'}} \right)^{c(i',j')} = \prod_{\theta=1}^{\tilde{\ell}} \prod_{i',j' \in \tilde{L}_{\theta}} \left( \frac{\tilde{B}_{j'}}{\tilde{B}_{i'}} \right)^{c(i',j')}
\end{equation}
where the final decomposition into linkage classes can be made by condition 1(b) of Definition \ref{resolving2}. We will show that the (\ref{weird3}) does not depend on any rate constant from any complex in the set $\tilde{\mathcal{C}}_I$.

Fix a $\theta \in \{ 1, \ldots, \tilde{\ell} \}$ such that $\tilde{\mathcal{C}}_R \cap \tilde{L}_{\theta} \not= \emptyset$. Notice that condition 1(b) of Definition \ref{resolving2} implies that there are at least two $i', j' \in \tilde{\mathcal{C}}_R \cap \tilde{L}_{\theta}$ such that $j' \not= i'$. We consider two cases.\\

\emph{Case 1:} Suppose $\tilde{\mathcal{C}}_I \cap \tilde{L}_{\theta} = \emptyset$. Since the spanning $i$-trees only span $\tilde{L}_\theta$, it follows that for any $p' \in \tilde{\mathcal{C}}_I$ we have that
\begin{equation}
\label{weird99}
\prod_{i',j' \in \tilde{L}_{\theta}} \left( \frac{\tilde{B}_{j'}}{\tilde{B}_{i'}} \right)^{c(i',j')}
\end{equation}
does not depend on any reaction from any $p' \in \tilde{\mathcal{C}}_I$.\\

\emph{Case 2:} Suppose there is an $p' \in \tilde{\mathcal{C}}_I \cap \tilde{L}_{\theta}$. By assumption $2.$ of Lemma \ref{mainlemma}, there is a $k' \in \tilde{L}_{\theta}$ such that, for every path $\tilde{P} = \{ \tilde{\mathcal{C}}_{\tilde{P}},\tilde{\mathcal{R}}_{\tilde{P}} \}$ from $p'$ to $i'$ we have $k' \in \tilde{\mathcal{C}}_{\tilde{P}}$. Let $\tilde{\mathcal{P}}(p',k')$ and $\tilde{\mathcal{P}}(k',i')$ denote the set of all paths from $p'$ to $k'$ and from $k'$ to $i'$, respectively. Now define
\[\begin{split} \tilde{\mathcal{C}}(p',k') & = \bigcup_{\tilde{P} \in \tilde{\mathcal{P}}(p',k')} \tilde{\mathcal{C}}_{\tilde{P}}\\ \tilde{\mathcal{R}}(p',k') & = \bigcup_{\tilde{P} \in \tilde{\mathcal{P}}(p',k')} \tilde{\mathcal{R}}_{\tilde{P}}.\end{split} \]
That is to say, $\tilde{\mathcal{C}}(p',k')$ and $\tilde{\mathcal{R}}(p',k')$ are the set of all complexes and reactions, respectively, which are on a path from $p'$ to $k'$.

Let $\tilde{\mathcal{T}}(p',k')$ denote the set of all $k'$-trees which span $\tilde{\mathcal{C}}(p',k')$. Note that every path from a complex in $\tilde{\mathcal{C}}(p',k')$ to $i'$ goes through $k'$, and that no path from $k'$ to $i'$ passes through $\tilde{\mathcal{C}}(p',k')$ (since it would return to $k'$). It follows that we may write any $\tilde{T} \in \tilde{\mathcal{T}}(i')$ as
\begin{equation}
\label{weird15}
\tilde{T} = \tilde{T}^* \cup \tilde{P}^* \cup \tilde{X}^*
\end{equation}
where $\tilde{T}^* \in \tilde{\mathcal{T}}(p',k')$, $\tilde{P}^* \in \tilde{\mathcal{P}}(k',i')$, and $\tilde{X}^* \in \tilde{\mathcal{X}}(\tilde{P}^*)$, where $\tilde{\mathcal{X}}(\tilde{P}^*)$ the set of configuration of reactions which, for a given path $\tilde{P}^* \in \tilde{\mathcal{P}}(k',i')$, connect the remaining complexes in $\tilde{L}_\theta$ to either $\tilde{P}^*$ or  $\tilde{T}^*$. That is to say, we construct $\tilde{T} \in \tilde{\mathcal{T}}(i')$ by first selecting a $k'$-tree on the reduced complex set $\tilde{\mathcal{C}}(p',k')$ (i.e. $\tilde{T}^*$), then connecting $k'$ to $i'$ with a direct path (i.e. $\tilde{P}^*$), and then connecting the remaining complexes to this structure (i.e. $\tilde{X}^*$). Notice that $\tilde{T}^*$ and $\tilde{P}^*$ may be chosen independently, and that $\tilde{X}^*$ depends on the chosen path $\tilde{P}^*$ but \emph{not} on the tree $\tilde{T}^*$.

We now construct $\tilde{B}_{i'}$ by considering all possible trees $\tilde{T} \in \tilde{\mathcal{T}}(i')$ constructed by (\ref{weird15}). We have that
\begin{equation}
\label{Kiprime}
\tilde{B}_{i'} = \left( \sum_{\tilde{T}^* \in \tilde{\mathcal{T}}(p',k')} \prod_{(i^*,j^*) \in \tilde{T}^*} \tilde{b}(i^*,j^*) \right) \left[ \sum_{\tilde{P}^* \in \tilde{\mathcal{P}}(k',i')}  \prod_{(i^*,j^*) \in \tilde{P}^*} \tilde{b}(i^*,j^*) \left( \sum_{\tilde{X}^* \in \tilde{\mathcal{X}}(\tilde{P}^*)} \prod_{(i^*,j^*) \in \tilde{X}^*} \tilde{b}(i^*,j^*)\right) \right]
\end{equation}
Now consider any $j' \in \tilde{\mathcal{C}}_R \cap \tilde{L}_{\theta}$, $j' \not= i'$. Noting that every path from $p'$ to $j'$ also goes through $k'$, we have
\begin{equation}
\label{Kjprime}
\tilde{B}_{j'} = \left( \sum_{\tilde{T}^* \in \tilde{\mathcal{T}}(p',k')} \prod_{(i^*,j^*) \in \tilde{T}^*} \tilde{b}(i^*,j^*) \right) \left[ \sum_{\tilde{P}^* \in \tilde{\mathcal{P}}(k',j')}  \prod_{(i^*,j^*) \in \tilde{P}^*} \tilde{b}(i^*,j^*) \left( \sum_{\tilde{X}^* \in \tilde{\mathcal{X}}(\tilde{P}^*)} \prod_{(i^*,j^*) \in \tilde{X}^*} \tilde{b}(i^*,j^*)\right) \right].
\end{equation}
Note that in (\ref{Kiprime}), the arrangements $\tilde{\mathcal{X}}(\tilde{P}^*)$ depend on the paths $\tilde{P}^* \in \tilde{\mathcal{P}}(k',i')$ while in (\ref{Kjprime}), they depend on the paths $\tilde{P}^* \in \tilde{\mathcal{P}}(k',j')$. It is important, however, that neither depends on any reaction from a complex in $\tilde{\mathcal{C}}(p',k')$ (the support of $\tilde{T}^*$ in both cases).

After simplifying, it follows from (\ref{Kiprime}) and (\ref{Kjprime}) that we have
\begin{equation}
\label{weird2}
\frac{\tilde{B}_{j'}}{\tilde{B}_{i'}} = \displaystyle{\frac{\sum_{\tilde{P}^* \in \tilde{\mathcal{P}}(k',j')}  \prod_{(i^*,j^*) \in \tilde{P}^*} \tilde{b}(i^*,j^*) \left( \sum_{\tilde{X}^* \in \tilde{\mathcal{X}}(\tilde{P}^*)} \prod_{(i^*,j^*) \in \tilde{X}^*} \tilde{b}(i^*,j^*)\right) }{\sum_{\tilde{P}^* \in \tilde{\mathcal{P}}(k',i')}  \prod_{(i^*,j^*) \in \tilde{P}^*} \tilde{b}(i^*,j^*) \left( \sum_{\tilde{X}^* \in \tilde{\mathcal{X}}(\tilde{P}^*)} \prod_{(i^*,j^*) \in \tilde{X}^*} \tilde{b}(i^*,j^*)\right)}}
\end{equation}
which does not depend on any complex in $\tilde{\mathcal{C}}(p',k')$, and therefore does not depend on $p'$. Since $p' \in \tilde{\mathcal{C}}_I \cap \tilde{L}_{\theta}$ was chosen arbitrarily, it follows that (\ref{weird99}) 
does not depend on any reaction from $p' \in \tilde{\mathcal{C}}_I$. Now consider an arbitrary pair $i,j \in h^{-1}(k')$ where $k' \in \tilde{\mathcal{C}}_I$. Applying the result of either case $1$ of case $2$ to (\ref{weird99}), we have that $\tilde{B}_{i,j}$ does not depend on any reaction from any $p' \in \tilde{\mathcal{C}}_I$.\\

It remains to connect the form (\ref{weird3}) to steady state resolvability as defined by Definition \ref{resolvable}. We make the following notes regarding the relationship between the definitions given in this paper, and Definition 6, Definition 10, Definition 11, and Lemma 4 in \cite{J1}:
\begin{enumerate}
\item
Definition 6 (translation) and Definition 10 (resolvability) in \cite{J1} emphasize the translation of individual reactions, whereas the definitions here emphasize the net flux out of a given source complex given a particular reaction-weight set. Nevertheless, we can clearly see that (\ref{weird3}) not depending on any reaction from any complex in $\tilde{\mathcal{C}}_I$ is sufficient to imply it does not depend on any reaction from the set required of Definition 10 in \cite{J1}. It follows that a translation satisfying the conditions of Lemma \ref{mainlemma} is resolvable as defined by Definition 10 of \cite{J1}.
\item
Definition 11 (construction of reweighted network) assigns reaction weights by arbitrarily selecting a single complex $i^* \in h^{-1}(k')$ for each $k' \in \tilde{\mathcal{C}}_I$ so that $(\tilde{C}_K)_{k'} = C_{i^*}$. For all reactions from this complex, the rate constants remain the same. For every other source complex $i \in h^{-1}(k')$, the reaction is scaled by a factor of the form (\ref{weird3}). Since the network is resolvable by Definition 10 of (\cite{J1}), property $1.$ of Definition \ref{translation} guarantees we may rescale rate constants in the same way to construct $\tilde{\mathcal{N}}(\tilde{\mathcal{K}})$ without altering the network structure of $\tilde{\mathcal{N}}(\tilde{\mathcal{B}})$. (Notice that condition (\ref{equation2}) is not sufficient to accomplish this by itself, as reactions may sum to zero in (\ref{equation2}) when they are reweighted. That is, in general we need the full conditions of property $1.$ of Definition \ref{translation}.)
\item
Since $\tilde{\mathcal{N}}(\tilde{\mathcal{B}})$ and the $\tilde{\mathcal{N}}(\tilde{\mathcal{K}})$ constructed by Definition 10 of \cite{J1} have the same network structure and $\tilde{\delta} = 0$, it follows from Lemma 4 of \cite{J1} that the mass action system (\ref{de2}) corresponding to $\mathcal{N}(\mathcal{K})$ and the generalized mass action system (\ref{gde}) corresponding to $\tilde{\mathcal{N}}(\tilde{\mathcal{K}})$ have the same steady states. That is to say, $\mathcal{N}(\mathcal{K})$ and $\tilde{\mathcal{N}}(\tilde{\mathcal{B}})$ are steady state resolvable, and we are done.
\end{enumerate}

\end{proof}

We now prove the main result of the paper, Theorem \ref{maintheorem}.

\begin{proof}[Proof of Theorem \ref{maintheorem}]

It is sufficient to prove the equivalence of the technical condition ($*$) of Lemma \ref{mainlemma} and the four technical conditions of Theorem \ref{maintheorem}.\\

\noindent \emph{Lemma \ref{mainlemma} $\Longrightarrow$ Theorem \ref{maintheorem}:} Suppose condition ($*$) of Lemma \ref{mainlemma} holds. That is, for every $p' \in \tilde{\mathcal{C}}_I$, there is a $k' \in \tilde{\mathcal{C}}$ such that every path from $p'$ to a $i' \in \tilde{\mathcal{C}}_R$ goes through $k'$. For a given $p' \in \tilde{\mathcal{C}}_I$, define $k'(p')$ to be the corresponding $k'$ and define $\tilde{\mathcal{C}}^*(p')$ to be the set of all complexes in $\tilde{\mathcal{C}}$ which can be reached from $p'$ without passing through $k'(p')$. Note that, by assumption, $\tilde{\mathcal{C}}^*(p') \cap \tilde{\mathcal{C}}_R = \emptyset$ and $\tilde{\mathcal{C}}^*(p') \cap k'(p') = \emptyset$. We define
\[\tilde{\mathcal{C}}^* = \bigcup_{p' \in \tilde{\mathcal{C}}_I} \tilde{\mathcal{C}}^*(p') \; \; \; \mbox{ and } \; \; \;  \tilde{\mathcal{R}}^* = \mathop{\bigcup_{(i',j') \in \tilde{\mathcal{R}}}}_{i' \in \tilde{\mathcal{C}}^*}(i',j').\]
By construction, these sets satisfy conditions $1.$ and $2.$ of Theorem \ref{maintheorem}.

We now construct the supplemental sets $\tilde{\mathcal{C}}^{**}$ and $\tilde{\mathcal{R}}^{**}$. Notice that $\tilde{\mathcal{C}}^{**}$ may contain complexes $k'$ selected earlier but that there must be a path from such a complex to another $k'$. We therefore define
\[\tilde{\mathcal{C}}^{**} = \left[ \bigcup_{p' \in \tilde{\mathcal{C}}_I} k'(p') \right] \setminus \tilde{\mathcal{C}}^*.\]
We also define $\tilde{\mathcal{R}}^{**}$ to be the set of all pairs $(k',p')$ where (1) $k' \in \tilde{\mathcal{C}}^{**}$, and (2) for a given $k'$, $p' \in \tilde{\mathcal{C}}_I$ is such that there is a path from $p'$ to $k'$ in the network $(\tilde{\mathcal{S}},\tilde{\mathcal{C}}^* \cup \tilde{\mathcal{C}}^{**},\tilde{\mathcal{R}}^*)$. Notice that these pairs need not be in the network $\tilde{\mathcal{N}}(\tilde{\mathcal{B}})$.

By construction, each linkage class of $(\tilde{\mathcal{S}},\tilde{\mathcal{C}}^* \cup \tilde{\mathcal{C}}^{**},\tilde{\mathcal{R}}^*)$ has a unique sink. (Otherwise, we would contradict condition $2.$ of Lemma \ref{mainlemma}.) The addition of the reaction set $\tilde{\mathcal{R}}^{**}$ clearly makes this network weakly reversible so that we have satisfied condition $3.$ of Theorem \ref{maintheorem}. Condition $4.$ follows from the uniqueness of the sinks in each linkage class prior to adding $\tilde{\mathcal{R}}^{**}$, since these sinks correspond to complexes in $\tilde{\mathcal{C}}^{**}$, and we are done.\\

\noindent \emph{Theorem \ref{maintheorem} $\Longrightarrow$ Lemma \ref{mainlemma}:} Suppose that there are sets $\tilde{\mathcal{C}}^*$, $\tilde{\mathcal{R}}^*$, $\tilde{\mathcal{C}}^{**}$, and $\tilde{\mathcal{R}}^{**}$ which satisfy conditions $1-4.$ of Theorem \ref{maintheorem}. Take an arbitrary $p' \in \tilde{\mathcal{C}}_I$. By condition $1.$ of Theorem \ref{maintheorem}, we have that $p' \in \tilde{\mathcal{C}}^*$. By condition $3.$ and $4.$, we have that the each linkage class of the network $(\tilde{\mathcal{S}},\tilde{\mathcal{C}}^* \cup \tilde{\mathcal{C}}^{**}, \tilde{\mathcal{R}}^*)$ has a unique sink at some complex in $\tilde{\mathcal{C}}^{**}$. From condition $2.$, however, we have that every path from $p'$ to this complex in $\tilde{\mathcal{N}}(\tilde{\mathcal{B}})$ is contained in $(\tilde{\mathcal{S}},\tilde{\mathcal{C}}^* \cup \tilde{\mathcal{C}}^{**}, \tilde{\mathcal{R}}^*)$. Since $\tilde{\mathcal{C}}^* \cup \tilde{\mathcal{C}}_R = \emptyset$ by condition $1.$, we have that for every $p' \in \tilde{\mathcal{C}}_I$ there is a $k' \in \tilde{\mathcal{C}}$ (the identified element in $\tilde{\mathcal{C}}^{**}$) such that every path from $p'$ to any complex in $\tilde{\mathcal{C}}_R$ goes through $k'$. It follows that condition ($*$) of Lemma \ref{mainlemma} is satisfied, and we are done.
\end{proof}

\section*{\hypertarget{appendixc}{Appendix C} (Code for Section \ref{milpsection})} \label{appendixD}

The following code corresponds to that derived in Section \ref{milpsection}. We derive the code into four sections: \emph{parameters}, \emph{decision variables}, \emph{objective functions}, and \emph{constraint sets}.\\

\textbf{Parameters:}

{\footnotesize \begin{flalign}\tag{\textbf{Par}}\label{parameters}
&
\left\{
\begin{array}{p{4cm}p{9.5cm}}
$n$ & Number of chemical species\\
$q$ & Number of kinetically-relevant complexes in $\mathcal{N}(\mathcal{K})$\\
$\tilde{m}$ & Number of hypothetical stoichiometric complexes in $\tilde{\mathcal{N}}(\tilde{\mathcal{B}})$\\
$s$ & Dimension of stoichiometric subspace of $\mathcal{N}(\mathcal{K})$\\
$\tilde{m} - s$ & Upper bound on number of linkage classes in $\tilde{\mathcal{N}}(\tilde{\mathcal{B}})$\\
$\tilde{\ell}^*$ & Upper bound on number of linkage classes in $(\tilde{\mathcal{S}},\tilde{\mathcal{C}}^* \cup \tilde{\mathcal{C}}^{**},\tilde{\mathcal{R}}^* \cup \tilde{\mathcal{R}}^{**})$\\
$V \in (0,1)^{\tilde{m} \times \tilde{m}}$ & Matrix of uniform random variables, $V_{i,j} = v[i,j] \in [\sqrt{\epsilon},1/\sqrt{\epsilon}]$\\
$Y \in \mathbb{Z}_{\geq 0}^{q \times n}$ & Complex matrix of $\mathcal{N}(\mathcal{K})$ \\
$\tilde{Y} \in \mathbb{Z}_{\geq 0}^{\tilde{m} \times n}$ & Complex matrix of $\tilde{\mathcal{N}}(\tilde{\mathcal{B}})$\\
$M = Y \cdot A(\mathcal{K}) \in \mathbb{R}^{q \times n}$ & Weighted stoichiometric matrix of $\mathcal{N}(\mathcal{K})$
\end{array}
\right\}
&
\end{flalign}}

\textbf{Decision variables:}

{\footnotesize \begin{flalign}\tag{\textbf{Dec}}\label{decisionvariables}
&
\left\{
\begin{array}{p{8cm}p{5.5cm}}
$\displaystyle{H[i,j'] \in \left\{ 0, 1 \right\},}$ & $\displaystyle{i = 1, \ldots, q, j' = 1, \ldots, \tilde{m}}$\\
$\displaystyle{\lambda[i,j'] \geq 0,}$ & $\displaystyle{i =1, \ldots, q, j' = 1, \ldots, \tilde{m}}$\\
$\displaystyle{\tilde{w}[i',j'] \geq 0,}$ & $\displaystyle{i',j' = 1, \ldots, \tilde{m}, i' \not= j'}$\\
$\displaystyle{\tilde{b}[i',j'] \geq 0,}$ & $\displaystyle{i',j' = 1 ,\ldots, \tilde{m}, i' \not= j'}$\\
$\displaystyle{\tilde{b}^*[i',j'] \geq 0,}$ & $\displaystyle{i', j' =1, \ldots, \tilde{m}, i' \not= j'}$\\
$\displaystyle{\tilde{b}^{**}[i',j'] \geq 0,}$ & $\displaystyle{i',j' =1, \ldots, \tilde{m}, i'\not=j'}$\\
$\displaystyle{\tilde{C}^*[i'] \in \left\{ 0, 1 \right\},}$ & $\displaystyle{i' = 1, \ldots, \tilde{m}}$\\
$\displaystyle{\tilde{C}^{**}[i'] \in \left\{ 0, 1 \right\},}$ & $\displaystyle{i' = 1, \ldots, \tilde{m}}$\\
$\displaystyle{\delta_I[i,j] \in \left\{ 0, 1 \right\},}$ & $\displaystyle{i,j = 1, \ldots, q, i < j}$\\
$\displaystyle{\delta_K[i,j] \in \left\{ 0, 1 \right\},}$ & $\displaystyle{i,j =1 ,\ldots, q, i < j}$\\
$\displaystyle{c[i,j] \geq 0}$ & $\displaystyle{i,j = 1, \ldots, q, i \not= j}$\\
$\displaystyle{\gamma[i',\theta] \in \left\{ 0, 1 \right\},}$ & $\displaystyle{i' = 1, \ldots, \tilde{m}, \theta=1, \ldots, \tilde{m}-s}$\\
$\displaystyle{\gamma_K[i,\theta] \in \left\{ 0, 1 \right\},}$ & $\displaystyle{i =1, \ldots, q, \theta=1, \ldots, \tilde{m}-s}$\\
$\displaystyle{\gamma^*[i',\theta] \geq 0,}$ & $\displaystyle{i'=1,\ldots, \tilde{m}, \theta = 1, \ldots, \tilde{\ell}^*}$\\
$\displaystyle{L[\theta] \in [0, 1]},$ & $\displaystyle{\theta=1, \ldots, \tilde{m}-s}$\\
$\displaystyle{\tilde{L}^*[\theta] \in \{ 0, 1 \}}$ & $\displaystyle{\theta = 1, \ldots, \tilde{\ell}^*}$
\end{array}
\right\}
&
\end{flalign} }

\textbf{Objective functions:}

{\footnotesize \begin{equation}\tag{\textbf{MinDef}}\label{deficiency}
\mbox{minimize} \; \; \sum_{\theta=1}^{\tilde{m}-s} L[\theta]
\end{equation}}

{\footnotesize \begin{equation}\tag{\textbf{MinC}}
\label{optimal}
\mbox{minimize} \; \; \; \; \epsilon  \sum_{i'=1}^{\tilde{m}} (\tilde{C}^*[i']+\tilde{C}^{**}[i']) 
\end{equation}}

\textbf{Constraint Sets:}


\noindent {\footnotesize \begin{flalign}\tag{\textbf{Trl1}}\label{constraint1}
&
\left\{
\begin{array}{p{8cm}p{5.5cm}}
$\displaystyle{\mathop{\sum_{j' = 1}^{\tilde{m}}}_{j' \not= i'} \tilde{b}[i',j'] \left( \tilde{Y}_{k,j'} - \tilde{Y}_{k,i'} \right) = \sum_{i=1}^q M_{k,i} \cdot H[i,i'],}$&$\displaystyle{k=1,\ldots,n,i' = 1, \ldots, \tilde{m}}$\\
$\displaystyle{\sum_{j'=1}^{\tilde{m}} H[i,j'] = 1,}$&$i=1, \ldots, q$\\
\end{array}
\right\}
&
\end{flalign} }

\noindent {\footnotesize \begin{flalign}\tag{\textbf{Trl2}}\label{constraint2}
&
\left\{ \begin{array}{p{8cm}p{5.5cm}} $\displaystyle{\sum_{i=1}^{q} H[i,j'] \leq 1,}$&$j'=1, \ldots, \tilde{m}$ \end{array} \right\}
&
\end{flalign} }

\noindent {\footnotesize \begin{flalign}\tag{\textbf{Trl3}}\label{constraint}
&
\left\{ \begin{array}{p{8cm}p{5.5cm}} $\lambda[i,j'] \leq (1/\epsilon) (1 - H[i,j']),$&$i=1,\ldots,q,j'=1,\ldots,\tilde{m}$ \\
$-\lambda[i,j'] \leq (1/\epsilon)H[i,j'],$&$i=1,\ldots,q,j'=1,\ldots,\tilde{m}$\\
$\displaystyle{\sum_{j' =1}^{\tilde{m}}\lambda[i,j'] = 0,}$&$i=1,\ldots,q$\\
$\displaystyle{\sum_{j'=1}^{\tilde{m}} \lambda[i,j'] \tilde{Y}_{k,j'} = M_{k,i}},$&$k=1,\ldots,n,i=1,\ldots,q$
\end{array} \right\}
&
\end{flalign} }

\noindent {\footnotesize \begin{flalign}\tag{\textbf{WR}}\label{weaklyreversible}
&
\left\{
\begin{array}{p{8cm}p{5.5cm}} $\displaystyle{\mathop{\sum_{j'=1}^{\tilde{m}}}_{j' \not= i'} \tilde{w}[i',j'] = \mathop{\sum_{j'=1}^{\tilde{m}}}_{j' \not= i'} \tilde{w}[j',i']}$, &$ i' = 1, \ldots, \tilde{m}$ \\
$-\tilde{w}[i',j'] \leq \displaystyle{-\epsilon \tilde{b}[i',j']},$ & $i',j' = 1, \ldots, \tilde{m}, i' \not= j$ \\
$\tilde{w}[i',j'] \leq \displaystyle{ (1/ \epsilon) \tilde{b}[i',j']},$ & $i',j' = 1, \ldots, \tilde{m}, i' \not= j'$
\end{array}
\right\}
&
\end{flalign} }

\noindent {\footnotesize \begin{flalign}\tag{\textbf{Def}}\label{partition}
&
\left\{ \begin{array}{p{8cm}p{5.5cm}} $\displaystyle{\sum_{\theta=1}^{\tilde{m}-s} \gamma[i',\theta] = 1,}$ & $i'=1, \ldots, \tilde{m}$ \\ $\displaystyle{\sum_{i'=1}^{\tilde{m}} \gamma[i',\theta] \leq (1/\epsilon) L[\theta],}$& $\theta = 1, \ldots, \tilde{m}-s$ \\
$\displaystyle{-\sum_{i'=1}^{\tilde{m}} \gamma[i',\theta] \leq - \epsilon L[\theta]},$ & $\theta = 1 ,\ldots, \tilde{m}-s$ \\
$\displaystyle{\tilde{b}[i',j'] \leq (1/\epsilon) (\gamma[i',\theta]-\gamma[j',\theta]+1)}$, & $i',j' = 1 ,\ldots, \tilde{m}, i' \not= j', \newline \theta = 1, \ldots, \tilde{m}-s$\\
$\displaystyle{\sum_{l = \theta+1}^{\tilde{m}-s} \gamma[i',l] \leq \sum_{j' = 1}^{i'-1} \gamma[j',\theta],}$ & $i' = 1, \ldots, \tilde{m}, \theta=1, \ldots, \tilde{m}-s, \theta \leq i'$ \end{array} \right\}
&
\end{flalign}}

\noindent {\footnotesize \begin{flalign}\tag{\textbf{Rsl1}}\label{constraint3}
&
\left\{ \begin{array}{p{8cm}p{5.5cm}}
$-\delta_I[i,j] \leq 1-H[i,k']-H[j,k'], $&$ i,j = 1, \ldots, q, i< j, \; k' = 1, \ldots, \tilde{m}$\\
$\delta_I[i,j] \leq 1-H[i,k']+H[j,k'], $&$ i,j =1, \ldots, q, i < j, \; k' = 1, \ldots, \tilde{m}$\\
$\delta_I[j,i] \leq 1-H[i,k']+H[j,k'], $&$ i,j =1, \ldots, q, j < i, \; k' = 1, \ldots, \tilde{m}$\\
$\delta_K[i,j] \leq (1/\epsilon)(c[i,j] + c[j,i]), $&$ i, j = 1, \ldots, q, i < j$\\
$-\delta_K[i,j] \leq -\epsilon (c[i,j] + c[j,i]), $&$ i,j = 1, \ldots, q, i < j$\\
$\gamma[i,\theta]-\gamma_K[k',\theta] \leq 1-H[i,k'], $&$ i=1, \ldots, q, k' = 1, \ldots, \tilde{m}, \newline \theta = 1, \ldots, \tilde{m} - s$\\
$\gamma_K[k',\theta]-\gamma[i,\theta] \leq 1-H[i,k'], $&$ i=1, \ldots, q, k' = 1, \ldots, \tilde{m}, \newline\theta = 1, \ldots, \tilde{m} - s$\\
$\delta_K[i,j] \leq 1-\gamma_K[j,k]+\gamma_K[i,k], $&$ i, j = 1, \ldots, q, i < j, \theta = 1, \ldots, \tilde{m} - s$\\
$\displaystyle{\mathop{\sum_{i,j=1}^{q}}_{i < j} v[i,j] \delta_I[i,j] \left( Y_{k,i}-Y_{k,j} \right) = \mathop{\sum_{i,j = 1}^{q}}_{i \not= j} c[i,j] \left( Y_{k,i} - Y_{k,j} \right),}$ & $k=1,\ldots,n$
\end{array} \right\}
&
\end{flalign}}

\noindent {\footnotesize \begin{flalign}\tag{\textbf{Rsl2}}\label{constraint4}
&
\left\{ \begin{array}{p{8cm}p{5.5cm}}
$-\tilde{C}^*[k'] \leq 1-H[i,k']-H[j,k'], $&$ i, j =1, \ldots, q, i < j, \; k' =1, \ldots, \tilde{m}$\\
$\delta_K[i,j] \leq 2-H[i,k']-\tilde{C}^*[k'], $&$ i, j = 1, \ldots, q, i < j, \; k' = 1, \ldots, \tilde{m}$\\
$\delta_K[j,i] \leq 2-H[i,k']-\tilde{C}^*[k'], $&$ i, j = 1, \ldots, q, j < i, \; k' = 1, \ldots, \tilde{m}$\\
$\tilde{b}^*[i',j'] \leq (1/\epsilon) \tilde{C}^*[j'], $&$ i',j' = 1 ,\ldots, \tilde{m}, i' \not= j'$ \\
$\tilde{b}^*[i',j'] \leq (1/\epsilon) \tilde{b}[i',j'], $&$ i',j' = 1, \ldots, \tilde{m}, i' \not= j'$ \\
$-\tilde{b}^*[i',j'] \leq \epsilon (1-\tilde{b}[i',j']-\tilde{C}^*[j']), $&$ i', j' = 1, \ldots, \tilde{m}, i' \not= j'$
\end{array} \right\}
&
\end{flalign}}

\noindent {\footnotesize \begin{flalign}{\tag{\textbf{Rsl3}}}\label{constraint5}
&
\left\{
\begin{array}{p{8cm}p{5.5cm}}
$\displaystyle{\tilde{C}^*[i'] + \tilde{C}^{**}[i'] \leq 1,} $&$i'=1, \ldots, \tilde{m}$\\
$\displaystyle{b^{**}[i',j'] \leq (1/\epsilon)\tilde{C}^{**}[i'],} $&$i',j'=1,\ldots, \tilde{m}, i'\not= j'$ \\
$\displaystyle{\sum_{\theta =1}^{\tilde{\ell}^*} \gamma^*[i',\theta] = \tilde{C}^*[i'] + \tilde{C}^{**}[i'],} $&$ i' = 1, \ldots, \tilde{m}$\\
$\tilde{b}^*[i',j'] + \tilde{b}^{**}[i',j'] \leq (1/\epsilon) (\gamma^*[i',\theta]-\gamma^*[j',\theta]+1), $&$ i', j' = 1 , \ldots, \tilde{m}, i \not= j, \theta = 1, \ldots, \tilde{\ell}^*$\\
$\displaystyle{\sum_{i' =1}^{\tilde{m}} \gamma^*[i',\theta] \leq  (1/\epsilon) \tilde{L}^*[\theta],} $&$ \theta = 1, \ldots, \tilde{\ell}^*$\\
$\displaystyle{-\sum_{i' = 1}^{\tilde{m}} \gamma^*[i',\theta]  \leq -\epsilon \tilde{L}^*[\theta],} $&$ \theta = 1, \ldots, \tilde{\ell}^*$\\
$\sum_{i'=1}^{\tilde{m}} \displaystyle{\tilde{C}^{**}[i'] = \sum_{\theta=1}^{\tilde{\ell}^*} \tilde{L}^*[\theta]}, $&$ i',j' = 1, \ldots, \tilde{m}, i' \not= j'$\\
$ \displaystyle{\mathop{\sum_{j' =1}^{\tilde{m}}}_{j' \not= i'} (\tilde{b}^{*}[i',j']+\tilde{b}^{**}[i',j']) = \mathop{\sum_{j' = 1}^{\tilde{m}}}_{j' \not= i'} (\tilde{b}^*[j',i'] + \tilde{b}^{**}[j',i']),} $&$ i' = 1, \ldots, \tilde{m}$\\
$\displaystyle{\sum_{l = \theta+1}^{\tilde{\ell}^*} \gamma[i',l] \leq \sum_{j' = 1}^{i'-1} \gamma[j',\theta],}$ & $i' = 1, \ldots, \tilde{m}, \theta=1, \ldots, \tilde{\ell}^*, \theta \leq i'.$
\end{array}
\right\}
&
\end{flalign}}


\begin{thebibliography}{10}

\bibitem{Boros2013}
Balasz Boros.
\newblock On the dependence of the existence of the positive steady states on
  the rate coefficients for deficiency-one mass action systems: single linkage
  class.
\newblock {\em J. Math. Chem.}, 51(9):2455--2490, 2013.

\bibitem{C-S}
Carsten Conradi and Anne Shiu.
\newblock A global convergence result for processive multisite phosphorylation
  systems.
\newblock 2014.
\newblock Available on the ArXiv at arxiv:1404.5524.

\bibitem{Dasgupta}
T.~Dasgupta, D.H. Croll, M.G.~Vander Heiden, J.W. Locasale, U.~Alon, L.C.
  Cantley, and J.~Gunawardena.
\newblock A fundamental trade off in covalent switching and its circumvention
  in glucose homeostasis.
\newblock {\em J. Biol. Chem.}, 289(19):13010--13025, 2014.

\bibitem{F1}
Martin Feinberg.
\newblock Complex balancing in general kinetic systems.
\newblock {\em Arch. Ration. Mech. Anal.}, 49:187--194, 1972.

\bibitem{Fe2}
Martin Feinberg.
\newblock Chemical reaction network structure and the stability of complex
  isothermal reactors: {I.} the deficiency zero and deficiency one theorems.
\newblock {\em Chem. Eng. Sci.}, 42(10):2229--2268, 1987.

\bibitem{Fe4}
Martin Feinberg.
\newblock Chemical reaction network structure and the stability of complex
  isothermal reactors: {II.} multiple steady states for networks of deficiency
  one.
\newblock {\em Chem. Eng. Sci.}, 43(1):1--25, 1988.

\bibitem{F2}
Martin Feinberg.
\newblock The existence and uniqueness of steady states for a class of chemical
  reaction networks.
\newblock {\em Arch. Ration. Mech. Anal.}, 132:311--370, 1995.

\bibitem{H-F}
Martin Feinberg and Fritz Horn.
\newblock Chemical mechanism structure and the coincidence of the
  stoichiometric and kinetic subspaces.
\newblock {\em Arch. Rational Mech. Anal.}, 66:83--97, 1977.


\bibitem{Hi}
Archibald Hill.
\newblock The possible effects of the aggregation of the molecules of
  haemoglobin on its dissociation curves.
\newblock {\em J. Physiol.}, 40(4), 1910.

\bibitem{H}
Fritz Horn.
\newblock Necessary and sufficient conditions for complex balancing in chemical
  kinetics.
\newblock {\em Arch. Ration. Mech. Anal.}, 49:172--186, 1972.

\bibitem{J1}
Matthew~D. Johnston.
\newblock Translated chemical reaction networks.
\newblock {\em Bull. Math. Bio.}, 76(5):1081--1116, 2014.

\bibitem{J-S5}
Matthew~D. Johnston, David Siegel, and G\'{a}bor Szederk\'{e}nyi.
\newblock Dynamical equivalence and linear conjugacy of chemical reaction
  networks: New results and methods.
\newblock {\em MATCH Commun. Math. Comput. Chem.}, 68(2), 2012.

\bibitem{J-S4}
Matthew~D. Johnston, David Siegel, and G\'{a}bor Szederk\'{e}nyi.
\newblock A linear programming approach to weak reversibility and linear
  conjugacy of chemical reaction networks.
\newblock {\em J. Math. Chem.}, 50(1):274--288, 2012.
  arxiv:1107.1659.

\bibitem{J-S6}
Matthew~D. Johnston, David Siegel, and G\'{a}bor Szederk\'{e}nyi.
\newblock Computing weakly reversible linearly conjugate chemical reaction
  networks with minimal deficiency.
\newblock {\em Math. Biosci.}, 241(1), 2013.

\bibitem{Karp}
Robert L. Karp, Mercedes~P\'{e}rez Mill\'{a}n,  Tathagata~Dasgupta, Alicia~Dickenstein, and
  Jeremy~Gunawardena.
\newblock Complex-linear invariants of biochemical networks.
\newblock {\em J. Theor. Biol.}, 311:130--138, 2012.

\bibitem{M-M}
Leonor Michaelis and Maud Menten.
\newblock Die kinetik der invertinwirkung.
\newblock {\em Biochem. Z.}, 49:333--369, 1913.

\bibitem{M-D-S-C}
Mercedes~P\'{e}rez Mill\'{a}n, Alicia Dickenstein, Anne Shiu, and Carsten
  Conradi.
\newblock Chemical reaction systems with toric steady states.
\newblock {\em Bull. Math. Biol.}, 74(5):1027--1065, 2012.

\bibitem{M-F-R-C-S-D}
Stefan M\"{u}ller, Elisenda Feliu, Georg Regensburger, Carsten Conradi, Anne
  Shiu, and Alicia Dickenstein.
\newblock Sign conditions for injectivity of generalized polynomial maps with
  applications to chemical reaction networks and real algebraic geometry.
\newblock 2013.

\bibitem{M-R}
Stefan M\"{u}ller and Georg Regensburger.
\newblock Generalized mass action systems: Complex balancing equilibria and
  sign vectors of the stoichiometric and kinetic-order subspaces.
\newblock {\em SIAM J. Appl. Math.}, 72(6):1926--1947, 2012.

\bibitem{R-S-H}
J\'{a}nos Rudan, G\'{a}bor Szederk\'{e}nyi, and Katalin Hangos.
\newblock Efficient computation of alternative structures for large kinetic
  systems using linear programming.
\newblock {\em MATCH Commun. Math. Comput. Chem.}, 71(1):71--92, 2014.

\bibitem{Sa}
Michael~A. Savageau.
\newblock Biochemical systems analysis {II}. the steady state solutions for an
  $n$-pool system using a power-law approximation.
\newblock {\em J. Theoret. Biol.}, 25:370--379, 1969.

\bibitem{Sh-F}
Guy Shinar and Martin Feinberg.
\newblock Structural sources of robustness in biochemical reaction networks.
\newblock {\em Science}, 327(5971):1389--1391, 2010.


\bibitem{Sz2}
Gabor Szederk\'{e}nyi.
\newblock Computing sparse and dense realizations of reaction kinetic systems.
\newblock {\em J. Math. Chem.}, 47:551--568, 2010.

\bibitem{Sz-H}
Gabor Szederk\'{e}nyi and Katalin Hangos.
\newblock Finding complex balanced and detailed balanced realizations of
  chemical reaction networks.
\newblock {\em J. Math. Chem.}, 49:1163--1179, 2011.

\bibitem{Sz-H-P}
Gabor Szederk\'{e}nyi, Katalin Hangos, and Tam\'{a}s P\'{e}ni.
\newblock Maximal and minimal realizations of chemical kinetics systems:
  computation and properties.
\newblock {\em MATCH Commun. Math. Comput. Chem.}, 65:309--332, 2011.

\bibitem{Sz-H-T}
Gabor Szederk\'{e}nyi, Katalin Hangos, and Zsolt Tuza.
\newblock Finding weakly reversible realizations of chemical reaction networks
  using optimization.
\newblock {\em MATCH Commun. Math. Comput. Chem.}, 67:193--212, 2012.

\bibitem{V-H}
Aizik~I. Vol'pert and Sergei~I. Hudjaev.
\newblock {\em Analysis in Classes of Discontinuous Functions and Equations of
  Mathematical Physics}.
\newblock Martinus Nijhoff Publishers, Dordrecht, Netherlands, 1985.

\end{thebibliography}
\end{document}


\newtheorem{assumption}{Assumption}[section]
\newtheorem{definition}{Definition}[section]
\newtheorem{lemma}{Lemma}[section]
\newtheorem{proposition}{Proposition}[section]
\newtheorem{theorem}{Theorem}[section]
\newtheorem{corollary}{Corollary}[section]
\newtheorem{remark}{Remark}[section]

\small

\title{Supplemental Material for:\\A Computational Approach to Steady State Correspondence of Regular and Generalized Mass Action Systems}
\author{Matthew D. Johnston \bigskip \\
\small Department of Mathematics\\
\small University of Wisconsin-Madison\\
\small 480 Lincoln Dr., Madison, WI 53706\\
\small email: mjohnston3@wisc.edu}
\date{}
\maketitle

\tableofcontents

\vspace{0.5in}

\noindent In this Supplemental Material, we provide detailed analysis of the applications contained in Section 5 of the main text. We present a more detailed overview of the application of Lemma 6.1 and Theorem 3.1 and elaborate on the mixed-integer linear programming (MILP) algorithm presented in Section 4. We conclude with an example which presents an avenue for continued research into improper translations.

Computations were performed with the GNU Linear Programming Kit (GLPK) on the author's personal use Toshiba Satellite laptop (AMD Quad-Core A6-Series APU, 6GB RAM). The corresponding MathProg code is contained in a separate file. Unless otherwise indicated, references to Definitions, Lemmas, and Theorems are from the main text. 

\section{Application I: EnvZ-OmpR Mechanism}

Consider the following reaction-weighted chemical reaction network $\mathcal{N}(\mathcal{K}) = (\mathcal{S},\mathcal{C},\mathcal{R},\mathcal{K})$ with reaction-weight set $\mathcal{K} = \{ k_i >0 \; | \; i=1, \ldots, 14 \}$:
\begin{equation}
\label{system3}
\begin{split}
& \displaystyle{XD \mathop{\stackrel{k_1}{\rightleftarrows}}_{k_2} X \mathop{\stackrel{k_3}{\rightleftarrows}}_{k_4} XT \stackrel{k_5}{\rightarrow} X_p} \\
& \displaystyle{X_p + Y \mathop{\stackrel{k_6}{\rightleftarrows}}_{k_7} X_pY \stackrel{k_8}{\rightarrow} X+Y_p} \\
& \displaystyle{XT + Y_p \mathop{\stackrel{k_9}{\rightleftarrows}}_{k_{10}} XTY_p \stackrel{k_{11}}{\rightarrow} XT+Y} \\
& \displaystyle{XD + Y_p \mathop{\stackrel{k_{12}}{\rightleftarrows}}_{k_{13}} XDY_p \stackrel{k_{14}}{\rightarrow} XD +Y.}
\end{split}
\end{equation}
This network corresponds to a hypothetical mechanism for the EnvZ/OmpR signaling system in \emph{Escherichia coli} which was introduced by Shinar and Feinberg in the Supporting Online Material of \cite{Sh-F} (with $X =$ EnvZ and $Y=$ OmpR). The corresponding mass action system was shown to possess \emph{absolute concentration robustness} in the concentration of phosphorylated OmpR, $Y_p$. The model was furthermore analyzed by P\'{e}rez Mill\'{a}n \emph{et al.} in \cite{M-D-S-C}, where it was shown to have toric steady states for all reaction-weights.

The model was also a primary example of network translation by Johnston in \cite{J1}. In the supplemental material of that paper, it was shown that, with the indexing
\[\begin{split} & X_1 = XD, \; \; \; X_2 = X, \; \; \; X_3 = XT, \; \; \; X_4 = X_p, \; \; \; X_5 = Y, \\ & X_6 = X_pY, \; \; \; X_7 = Y_p, \; \; \; X_8 = XTY_p, \; \; \; X_9 = XDY_p
\end{split}\] 
the network could be corresponded by the translation scheme
\begin{equation}
\label{234}
\begin{split}
& \displaystyle{X_1 \rightleftarrows X_2 \rightleftarrows X_3 \rightarrow X_4} \; \; \; \; \; \; \; \; \; \; \; \; \; \; \; \; \; \; \; \; (+ X_1 + X_3 + X_5) \\
& \displaystyle{X_4 + X_5 \rightleftarrows X_6 \rightarrow X_2 + X_7} \; \; \; \; \; \; \; \; \; \; \; \; (+ X_1 + X_3) \\
& \displaystyle{X_3 + X_7 \rightleftarrows X_8 \rightarrow X_3 + X_5} \; \; \; \; \; \; \; \; \; \; \; \; (+ X_1 + X_2)\\
& \displaystyle{X_1 + X_7 \rightleftarrows X_9 \rightarrow X_1 + X_5 \; \; \; \; \; \; \; \; \; \; \; \; (+ X_2 + X_3).}
\end{split}
\end{equation}
to the following weakly reversible generalized chemical reaction network $\tilde{\mathcal{N}}(\tilde{\mathcal{B}})$, with reaction weight set $\tilde{\mathcal{B}} = \{\tilde{b}_i = k_i\; | \; i = 1, \ldots, 14 \}$:
\begin{equation}
\label{system4}
\begin{split} & 2X_1 + X_3 + X_5 \mathop{\stackrel{\tilde{b}_1}{\rightleftarrows}}_{\tilde{b}_2} X_1 + X_2 + X_3 + X_5 \mathop{\stackrel{\tilde{b}_3}{\rightleftarrows}}_{\tilde{b}_4} X_1 + 2X_3 + X_5 \\ & \; \; \; \; \; \; \; \; \; \; \; \; \; \; \; \; \; \; \; \; \; \; \; \; \nearrow \hspace{-0.1cm} {}_{\tilde{b}_{14}} \; \; \; \; \; \; \; \; \; \; \; \; \uparrow_{\tilde{b}_{11}} \; \; \; \; \; \; \; \; \; \; \; \; \; \; \; \; \; \; \; \; \; \; \; \; \; \; \downarrow_{\tilde{b}_5} \\ & \; \; X_2 + X_3 + X_9 \; \; \; \; \; \; \; \; \; \; \; X_1 + X_2 + X_8 \; \; \; \; \; \; \; X_1 + X_3 + X_4 + X_5 \\ & \; \; \; \; \; \; \; \; \; \; \; \; \; \; \; \; \; \; \; {}_{\tilde{b}_{12}} \hspace{-0.1cm} \nwarrow \hspace{-0.15cm} \searrow \hspace{-0.1cm} {}^{\tilde{b}_{13}} \; \; \; \; \; \; \; \; {}_{\tilde{b}_9} \uparrow \downarrow {}_{\tilde{b}_{10}} \; \; \; \; \; \; \; \; \; \; \; \; \; \; \; \; \; \; \; {}_{\tilde{b}_7} \uparrow \downarrow {}_{\tilde{b}_6} \\ & \; \; \; \; \; \; \; \; \; \; \; \; \; \; \; \; \; \; \; \; \; \; \; \; \; \; \; \; X_1 + X_2 + X_3 + X_7 \stackrel{\tilde{b}_8}{\leftarrow} X_1 + X_3 + X_6. \end{split}
\end{equation}

It was shown in \cite{J1} that the steady states of the mass action equations corresponding to (\ref{system3}) and the generalized mass action equations corresponding to (\ref{system4}) coincide for the reaction weight set $\tilde{\mathcal{K}}$ with $\tilde{k}_i = k_i$, $i=1, \ldots, 14$, $i \not= 12$, and 
\begin{equation}
\label{91}
\tilde{k}_{12} = \left( \frac{k_2(k_4+k_5)}{k_1k_3} \right) k_{12}.
\end{equation}
Note that this set does not correspond to the reaction-weight set $\tilde{\mathcal{B}}$ consistent with Definition 3.2 and (\ref{system4}). In situations where the translation is improper but steady state resolvable, it is a scaling of these rate constants which produces the steady state equivalent generalized network $\tilde{\mathcal{N}}(\tilde{\mathcal{K}})$.  

\begin{figure}[t]
\centering
\includegraphics[width=12cm]{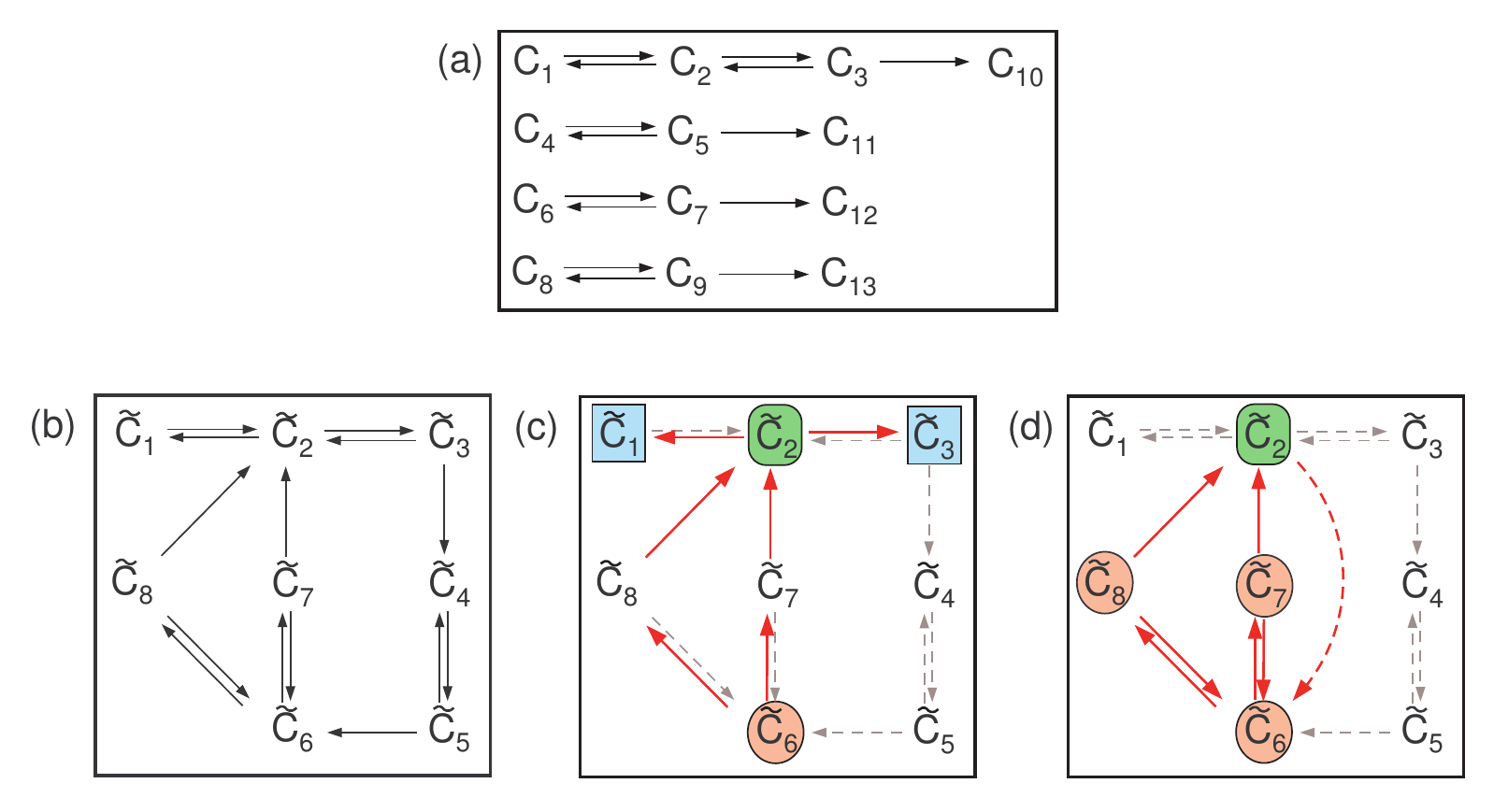}
\caption{\small Correspondence between (a) original network $\mathcal{N}(\mathcal{K})$ (\ref{system3}) and (b) computationally-determined translation $\tilde{\mathcal{N}}(\tilde{\mathcal{B}})$ (\ref{system4}). In (c), we identify the improper complex set $\tilde{\mathcal{C}}_I$ (pink) and a resolving complex set $\tilde{\mathcal{C}}_R$ (blue) which has the property that every path (red arrow) from an element of $\tilde{\mathcal{C}}_I$ to an element of $\tilde{\mathcal{C}}_R$ passes through the common complex $\tilde{C}_2$ (green). That is, the network satisfies the technical assumptions of Lemma 6.1. In (d), we identify the sets $\tilde{\mathcal{C}}^* = \{ 6, 7, 8 \}$ (pink), $\tilde{\mathcal{C}}^{**} = \{ 2 \}$ (green), $\tilde{\mathcal{R}}^* = \{ (6,7),(6,8),(7,2),(7,6),(8,2),(8,6)\}$ (solid red arrows), and $\tilde{\mathcal{R}}^{**} = \{ (2,6)\}$ (dashed red arrow). It is clear that the subnetwork $(\tilde{\mathcal{S}},\tilde{\mathcal{C}}^* \cup \tilde{\mathcal{C}}^{**},\tilde{\mathcal{R}}^* \cup \tilde{\mathcal{R}}^{**})$ is weakly reversible. It follows that the network satisfies the technical conditions of Theorem 3.1.}
\label{figure1}
\end{figure}

For completeness, we re-iterate here the analysis summarized in the main text. We start by re-indexing the complexes in (\ref{system3}) and (\ref{system4}). We index the complex set $\mathcal{C}$ according to:
\[\begin{split} &C_1 = X_1, \; C_2 = X_2, \; C_3 = X_3, \; C_4 =X_4 + X_5, \; C_5 = X_6,\\
&C_6 = X_3 + X_7, \; C_7 = X_8, \; C_8 =X_1 + X_7 , \; C_9 = X_9,\\
&C_{10} = X_4 , \; C_{11} =X_2 + X_7, \; C_{12} = X_3 + X_5 , \; C_{13} = X_1 + X_5,\end{split}\]
the translated complex set $\tilde{\mathcal{C}}$ according to:
\[\begin{split}& \tilde{C}_1 = 2X_1 + X_3 + X_5, \; \tilde{C}_2 = X_1 + X_2 + X_3 + X_5, \; \tilde{C}_3 = X_1 + 2X_3 + X_5,\\
& \tilde{C}_4 = X_1 + X_3 + X_4 + X_5, \; \tilde{C}_5 = X_1 + X_3 + X_6, \; \tilde{C}_6 = X_1 + X_2 + X_3 + X_7,\\
&\tilde{C}_7 = X_1 + X_2 + X_8, \; \tilde{C}_8 = X_2 + X_3 + X_9,\end{split}\]
and the translated kinetic complex set $\tilde{\mathcal{C}}_K$ according to:
\begin{equation}
\label{kinetic}
\begin{split} & (\tilde{C}_K)_1 = X_1, \; (\tilde{C}_K)_2 = X_2, \; (\tilde{C}_K)_3 = X_3, \; (\tilde{C}_K)_4 = X_4 + X_5, \\ & (\tilde{C}_K)_5 = X_6, \; (\tilde{C}_K)_6 = X_3 + X_7, \; (\tilde{C}_K)_7 = X_8, \; (\tilde{C}_K)_8 = X_9.\end{split}
\end{equation}
The corresponding network structures are given by Figure \ref{figure1}(a) and (b), respectively.

Notice first of all that $C_6 = X_3 + X_7$ and $C_8 = X_1 + X_7$ are both translated to $\tilde{C}_6 = X_1 + X_2 + X_3 + X_7$ by (\ref{234}) but only $X_3 + X_7$ appears as a kinetic complex (specifically $(\tilde{C}_K)_6 = X_3 + X_7$). It follows that the translation is improper and, specifically, that $h(6) = 6$ and $h(8) = 6$. In accordance with Definition 3.5, we have: 
\[\begin{split}
\tilde{\mathcal{C}}_I & = \{ 6 \},\\
h^{-1}(6) & = \{ 6, 8 \},\\
\tilde{S}_I & = \mbox{span} \{ y_8 - y_6 \} \\ & = \mbox{span} \{ (1,0,-1,0,0,0,0,0,0) \}.
\end{split}\]
To determine the resolving complex set $\tilde{\mathcal{C}}_R$ according to Definition 3.6, we consider $h^{-1}(6) = \{ 6,8 \}$. We can see that
\begin{equation}
\label{11111}
y_8 - y_6 = (1,0,0,0,0,0,0,0) - (0,0,1,0,0,0,0,0) = (\tilde{y}_K)_1 - (\tilde{y}_K)_3.
\end{equation}
We may therefore choose $c(1,3) = 1$ and $c(i',j') = 0$ for all other $i',j' = 1, \ldots, 8$, so that $\tilde{\mathcal{C}}_R = \{ 1,3 \}$ is our resolving complex set.

Notice that (\ref{11111}) is equivalent to the identity
\begin{equation}
\label{11112}
 \mathbf{x}^{y_8} =x_1 x_7  = \left( \frac{x_1}{x_3} \right) x_3 x_7= \left( \frac{\mathbf{x}^{(\tilde{y}_K)_1}}{\mathbf{x}^{(\tilde{y}_K)_3}} \right) \mathbf{x}^{(\tilde{y}_K)_6}.
\end{equation}
That is we are able to relate the untranslated monomial $x_1 x_7$ to $x_1$, $x_3$, and $x_3x_7$, all of which correspond to kinetic complexes in the set (\ref{kinetic}). We may therefore think of any reaction corresponding to $x_1 x_7$ as a reaction from source $x_3 x_7$ with the additional \emph{state dependent} rate $x_1/x_3$. It can be checked that the mass action system corresponding to (\ref{system3}) and the generalized mass action system corresponding to (\ref{system4}) differ in only the monomials $k_{12} x_1x_7$ and $\tilde{k}_{12} x_3 x_7$. It was shown in \cite{J1} that at any steady state of the generalized mass action system corresponding to (\ref{system4}) we have
\[\frac{x_1}{x_3} = \frac{k_2(k_4+k_5)}{k_1k_3}\]
so that we may ``resolve'' the state dependent term to get
\[\tilde{k}_{12} x_3 x_7 = k_{12} x_1 x_7 = k_{12} \left( \frac{x_1}{x_3} \right) x_3 x_7 = k_{12} \left( \frac{k_2(k_4+k_5)}{k_1k_3} \right) x_3 x_7. \]
This gives an explicit equation in the undetermined rate constant $\tilde{k}_{12}$ which can be solved for directly to get (\ref{91}).

We now apply Lemma 6.1 (see Figure \ref{figure1}(a)). We identify the set of paths from the improper complex $\tilde{C}_6$ (pink) to the resolving complexes $\tilde{C}_1$ and $\tilde{C}_3$ (blue). The relevant paths $\tilde{P} = \{ \tilde{\mathcal{C}}_{\tilde{P}}, \tilde{\mathcal{R}}_{\tilde{P}} \}$ (red arrows) are given as follows:
\[\begin{split}
& (1) \; \; \tilde{P} = \{ \tilde{\mathcal{C}}_{\tilde{P}},\tilde{\mathcal{R}}_{\tilde{P}} \} \mbox{ with } \tilde{\mathcal{C}}_{\tilde{P}} = \{ 1,2,6,7 \}, \tilde{\mathcal{R}}_{\tilde{P}} = \{(2,1),(6,7),(7,2)\} \\
& (2) \; \; \tilde{P} = \{ \tilde{\mathcal{C}}_{\tilde{P}},\tilde{\mathcal{R}}_{\tilde{P}} \} \mbox{ with } \tilde{\mathcal{C}}_{\tilde{P}} = \{ 1,2,6,8 \}, \tilde{\mathcal{R}}_{\tilde{P}} = \{(2,1),(6,8),(8,2)\} \\
& (3) \; \; \tilde{P} = \{ \tilde{\mathcal{C}}_{\tilde{P}},\tilde{\mathcal{R}}_{\tilde{P}} \} \mbox{ with } \tilde{\mathcal{C}}_{\tilde{P}} = \{ 2,3,6,7 \}, \tilde{\mathcal{R}}_{\tilde{P}} = \{(2,3),(6,7),(7,2)\} \\
& (4) \; \; \tilde{P} = \{ \tilde{\mathcal{C}}_{\tilde{P}},\tilde{\mathcal{R}}_{\tilde{P}} \} \mbox{ with } \tilde{\mathcal{C}}_{\tilde{P}} = \{ 2,3,6,8 \}, \tilde{\mathcal{R}}_{\tilde{P}} = \{(2,3),(6,8),(8,2)\}.\end{split}\]
We can see that $\{2\} \subset \tilde{\mathcal{C}}_{\tilde{P}}$ for all such paths, so that every path goes through $\tilde{C}_2$. It follows that we have $\tilde{C}_{k'} = \tilde{C}_2$ (green) and, since $\tilde{\delta} = 0$, it follows by Lemma 6.1 that $\tilde{\mathcal{N}}(\tilde{\mathcal{B}})$ and $\mathcal{N}(\mathcal{K})$ are steady state resolvable.

We now apply Theorem 3.1 (see Figure \ref{figure1}(d)). Consider the network $(\tilde{\mathcal{S}},\tilde{\mathcal{C}}^{*} \cup \tilde{\mathcal{C}}^{**},\tilde{\mathcal{R}}^{*} \cup \tilde{\mathcal{R}}^{**})$ where $\tilde{\mathcal{C}}^* = \{ 6, 7, 8 \}$ (pink), $\tilde{\mathcal{C}}^{**} = \{ 2 \}$ (green), $\tilde{\mathcal{R}}^* = \{ (6,7),(6,8),\\(7,2),(7,6),(8,2),(8,6)\}$ (solid red arrows), and $\tilde{\mathcal{R}}^{**} = \{ (2,6) \}$ (dashed red arrows). It is clear that $| \tilde{\mathcal{C}}^{**} | = | \tilde{\mathcal{L}}^* |$ and that $(\tilde{\mathcal{S}},\tilde{\mathcal{C}}^{*} \cup \tilde{\mathcal{C}}^{**},\tilde{\mathcal{R}}^{*} \cup \tilde{\mathcal{R}}^{**})$ is weakly reversible. Notice that $\tilde{\mathcal{C}}^{*}$ does contain any complex in $\tilde{\mathcal{C}}_R$ (although it is permissible to have $\mathcal{C}^{**}$ contain such a complex). Also notice that reactions in $\tilde{\mathcal{R}}^{**}$ need not be in the original network, nor be singletons, but that, by the construction of $\tilde{\mathcal{C}}^{**}$, they do need to originate at a sink of linkage class in $\tilde{\mathcal{L}}^{*}$. Since $\tilde{\delta} = 0$, it follows by Theorem 3.1 that $\tilde{\mathcal{N}}(\tilde{\mathcal{B}})$ and $\mathcal{N}(\mathcal{K})$ are steady state resolvable.

We now apply the mixed-integer linear programming algorithm presented in Section 5 of the main text. We set $\epsilon = 0.1$, $\tilde{\ell} = \tilde{m}-s = 10$, and $\tilde{\ell}^* = 2$. We initialize the relevant matrices as follows, where the notation is taken from the main text:
\[{\footnotesize Y = \left[
\begin{array}{ccccccccc}
1 & 0 & 0 & 0 & 0 & 0 & 0 & 1 & 0 \\
0 & 1 & 0 & 0 & 0 & 0 & 0 & 0 & 0 \\
0 & 0 & 1 & 0 & 0 & 1 & 0 & 0 & 0 \\
0 & 0 & 0 & 1 & 0 & 0 & 0 & 0 & 0 \\
0 & 0 & 0 & 1 & 0 & 0 & 0 & 0 & 0 \\
0 & 0 & 0 & 0 & 1 & 0 & 0 & 0 & 0 \\
0 & 0 & 0 & 0 & 0 & 1 & 0 & 1 & 0 \\
0 & 0 & 0 & 0 & 0 & 0 & 1 & 0 & 0 \\
0 & 0 & 0 & 0 & 0 & 0 & 0 & 0 & 1 
\end{array} \right]},\]
\[{\footnotesize \tilde{Y} = \left[ \begin{array}{ccccccccccccccccc}
1 & 0 & 0 & 0 & 0 & 0 & 0 & 1 & 0 & 2 & 1 & 1 & 1 & 1 & 1 & 1 & 0 \\
0 & 1 & 0 & 0 & 0 & 0 & 0 & 0 & 0 & 0 & 1 & 0 & 0 & 0 & 1 & 1 & 1 \\
0 & 0 & 1 & 0 & 0 & 1 & 0 & 0 & 0 & 1 & 1 & 2 & 1 & 1 & 1 & 0 & 1 \\
0 & 0 & 0 & 1 & 0 & 0 & 0 & 0 & 0 & 0 & 0 & 0 & 1 & 0 & 0 & 0 & 0 \\
0 & 0 & 0 & 1 & 0 & 0 & 0 & 0 & 0 & 1 & 1 & 1 & 1 & 0 & 0 & 0 & 0 \\
0 & 0 & 0 & 0 & 1 & 0 & 0 & 0 & 0 & 0 & 0 & 0 & 0 & 1 & 0 & 0 & 0 \\
0 & 0 & 0 & 0 & 0 & 1 & 0 & 1 & 0 & 0 & 0 & 0 & 0 & 0 & 1 & 0 & 0 \\
0 & 0 & 0 & 0 & 0 & 0 & 1 & 0 & 0 & 0 & 0 & 0 & 0 & 0 & 0 & 1 & 0 \\
0 & 0 & 0 & 0 & 0 & 0 & 0 & 0 & 1 & 0 & 0 & 0 & 0 & 0 & 0 & 0 & 1
\end{array} \right]},\]
and
\[{\tiny M = \left[ \begin{array}{ccccccccc}
-k[1] & k[2] & 0 & 0 & 0 & 0 & 0 & -k[12] & 0 \\
k[1] & -k[2,3] & k[4] & 0 & k[8] & 0 & 0 & 0 & 0 \\
0 & k[3] & -k[4,5] & 0 & 0 & -k[9] & k[10,11] & 0 & k[13,14] \\
0 & 0 & k[5] & -k[6] & k[7] & 0 & 0 & 0 & 0 \\
0 & 0 & 0 & -k[6] & k[7] & 0 & k[11] & 0 & k[14] \\
0 & 0 & 0 & k[6] & -k[7,8] & 0 & 0 & 0 & 0 \\
0 & 0 & 0 & 0 & k[8] & -k[9] & k[10] & -k[12] & k[13] \\
0 & 0 & 0 & 0 & 0 & k[9] & -k[10,11] & 0 & 0 \\
0 & 0 & 0 & 0 & 0 & 0 & 0 & k[12] & -k[13,14]
\end{array} \right]},\]
where $k[i,j] := k[i]+k[j]$ and $k[i]$ are determined stochastically from the range $[\sqrt{\epsilon},1/\sqrt{\epsilon}]$ for all $i=1,\ldots, 14$. The parameters $v[i,j]$, $i,j = 1, \ldots, 9$, $i \not= j$, are also determined stochastically from the range $[\sqrt{\epsilon},1/\sqrt{\epsilon}]$. The algorithm was successfully run to completion $25$ times. The mean time to completion in the sample was $2.788$ seconds with a standard deviation of $1.4898$ seconds. In each realization in the sample, the algorithm determined the generalized reaction-weighted chemical reaction network given by (\ref{system4}).

\begin{remark}
Note that computational efficiency depends on the user-determined parameter values $\epsilon$ and $\tilde{\ell}^*$. In general, smaller $\epsilon$ values increase the numerical instability of the optimizer and larger $\tilde{\ell}^*$ values increase the required computation time.
\end{remark}

\section{Application II: PFK-2/FBPase-2 Mechanism}

\begin{figure}[t]
\centering
\includegraphics[width=13cm]{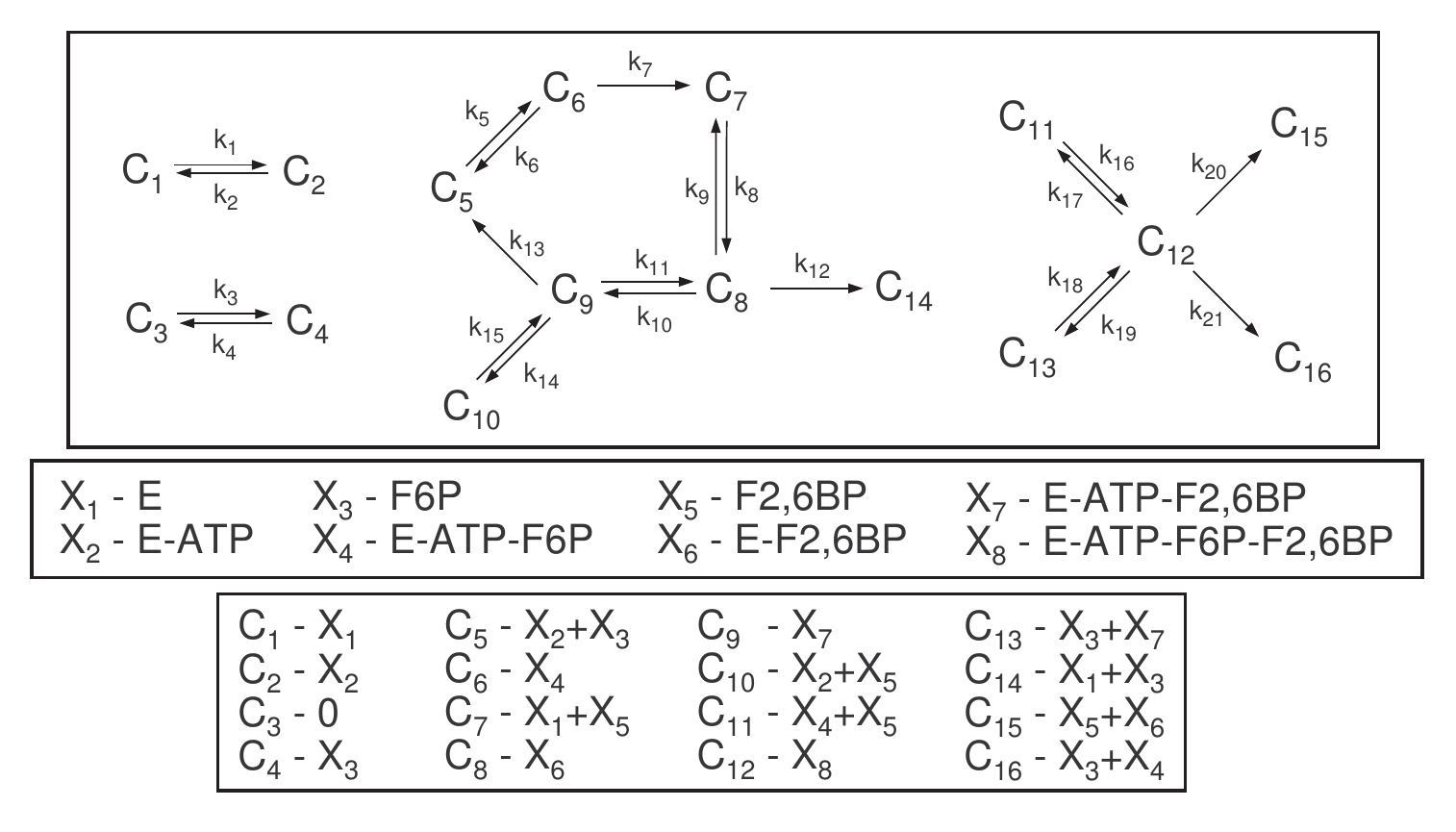}
\caption{Hypothetical PFK-2/FBPase-2 mechanism in mammalian cells.}
\label{figure3}
\end{figure}

Consider the PFK-2/FBPase-2 mechanism given in Figure \ref{figure3}. This model is slightly modified from the one proposed by Dasgupta \emph{et al.} in \cite{Dasgupta,Karp} to include a reversible reaction pair $C_3 \rightleftarrows C_4$ corresponding to $0 \rightleftarrows X_3$. That is, our mechanism allows for inflow and outflow of Fructose 6-phosphate ($F6P$). We defer biochemical justification and analysis of this mechanism to \cite{Dasgupta,Karp}.

The network in Figure \ref{figure3} is not weakly reversible, has a seven-dimensional stoichiometric subspace (i.e. $s = 7$), and a deficiency of five (i.e. $\delta = 5$). It is notable that, for some choices of reaction weights, the kinetic subspace may be smaller than the stoichiometric subspace as a result of there being two terminal strongly linked components in the fourth linkage class (see \cite{H-F}). 

We now apply the computational algorithm of Section 5 of the main text. We set $\epsilon = 0.1$, $\tilde{\ell} = \tilde{m}-s = 2$, and $\tilde{\ell}^* = 2$, and initialize the relevant matrices as follows, where the notation is taken from the main text:
\[{\footnotesize Y = \left[ \begin{array}{ccccccccccccc}
1&	0&	0&	0&	0&	0&	1&	0&	0&	0&	0&	0&	0\\
0&	1&	0&	0&	1&	0&	0&	0&	0&	1&	0&	0&	0\\
0&	0&	0&	1&	1&	0&	0&	0&	0&	0&	0&	0&	1\\
0&	0&	0&	0&	0&	1&	0&	0&	0&	0&	1&	0&	0\\
0&	0&	0&	0&	0&	0&	1&	0&	0&	1&	1&	0&	0\\
0&	0&	0&	0&	0&	0&	0&	1&	0&	0&	0&	0&	0\\
0&	0&	0&	0&	0&	0&	0&	0&	1&	0&	0&	0&	1\\
0&	0&	0&	0&	0&	0&	0&	0&	0&	0&	0&	1&	0\end{array} \right]},\]
\[{\footnotesize \tilde{Y} = \left[ \begin{array}{ccccccccccc}
0&	0&	1&	0&	0&	1&	0&	0&	0&	0&	0\\
0&	0&	0&	1&	0&	0&	0&	0&	1&	0&	0\\
0&	1&	2&	2&	1&	1&	1&	1&	1&	0&	0\\
0&	0&	0&	0&	1&	0&	0&	0&	0&	1&	0\\
0&	0&	0&	0&	0&	1&	0&	0&	1&	1&	0\\
0&	0&	0&	0&	0&	0&	1&	0&	0&	0&	0\\
0&	0&	0&	0&	0&	0&	0&	1&	0&	0&	0\\
0&	0&	0&	0&	0&	0&	0&	0&	0&	0&	1\end{array} \right]
}\]
and
\[\begin{split} & {\tiny M = \left[ \begin{array}{cccccccc}
0&	0&	-k[3]&	k[4]&	0&	k[7]&	-k[8]&	k[9,11]\\
0&	0&	k[3]&	-k[4]&	-k[5]&	k[6]&	0&	0\\
k[1]&	-k[2]&	0&	0&	-k[5]&	k[6]&	0&	k[11]\\
0&	0&	0&	0&	k[5]&	-k[6,7]&	0&	0\\
0&	0&	0&	0&	0&	k[7]&	-k[8]&	k[9]\\
0&	0&	0&	0&	0&	0&	k[8]&	-k[9,10,11]\\
0&	0&	0&	0&	0&	0&	0&	k[10]\\
0&	0&	0&	0&	0&	0&	0&	0 \end{array} \right. \hspace{1in}}
\\ &
{\tiny \hspace{2in} \left. \begin{array}{ccccc}
0 & 0 & 0 & 0 & 0 \\
k[12,14]& -k[15]&	0&	0&	0 \\
k[12]&	0&	0&	k[18,19]&	-k[21]\\
0&	0&	-k[16]&	k[17,19]&	0\\
k[14]&	-k[15]&	-k[16]&	k[17,20]&	0\\
k[13]&	0&	0&	k[20]&	0\\
-k[12,13,14]&	k[15]&	0&	k[18]&	-k[21]\\
0&	0&	k[16]&	-k[17,18,19,20]&	k[21]
\end{array} \right]} \end{split}\]
where $k[i_1,\ldots,i_n] = k[i_1] + \cdots + k[i_n]$. We choose the parameters $k[i] \in [\sqrt{\epsilon}, 1/\sqrt{\epsilon}]$, $i=1, \ldots, 21$, to be random parameters chosen uniformly from their range and further impose that $k[19]=k[20]$. The parameters $v[i,j]$, $i,j = 1, \ldots, 13$, $i \not= j$, are determined stochastically from the range $[\sqrt{\epsilon},1/\sqrt{\epsilon}]$. The algorithm was run to completion $25$ times. The mean time to completion in the sample was $6.604$ seconds with a standard deviation of $2.6871$ seconds.

\begin{figure}[h!]
\centering
\includegraphics[width=12cm]{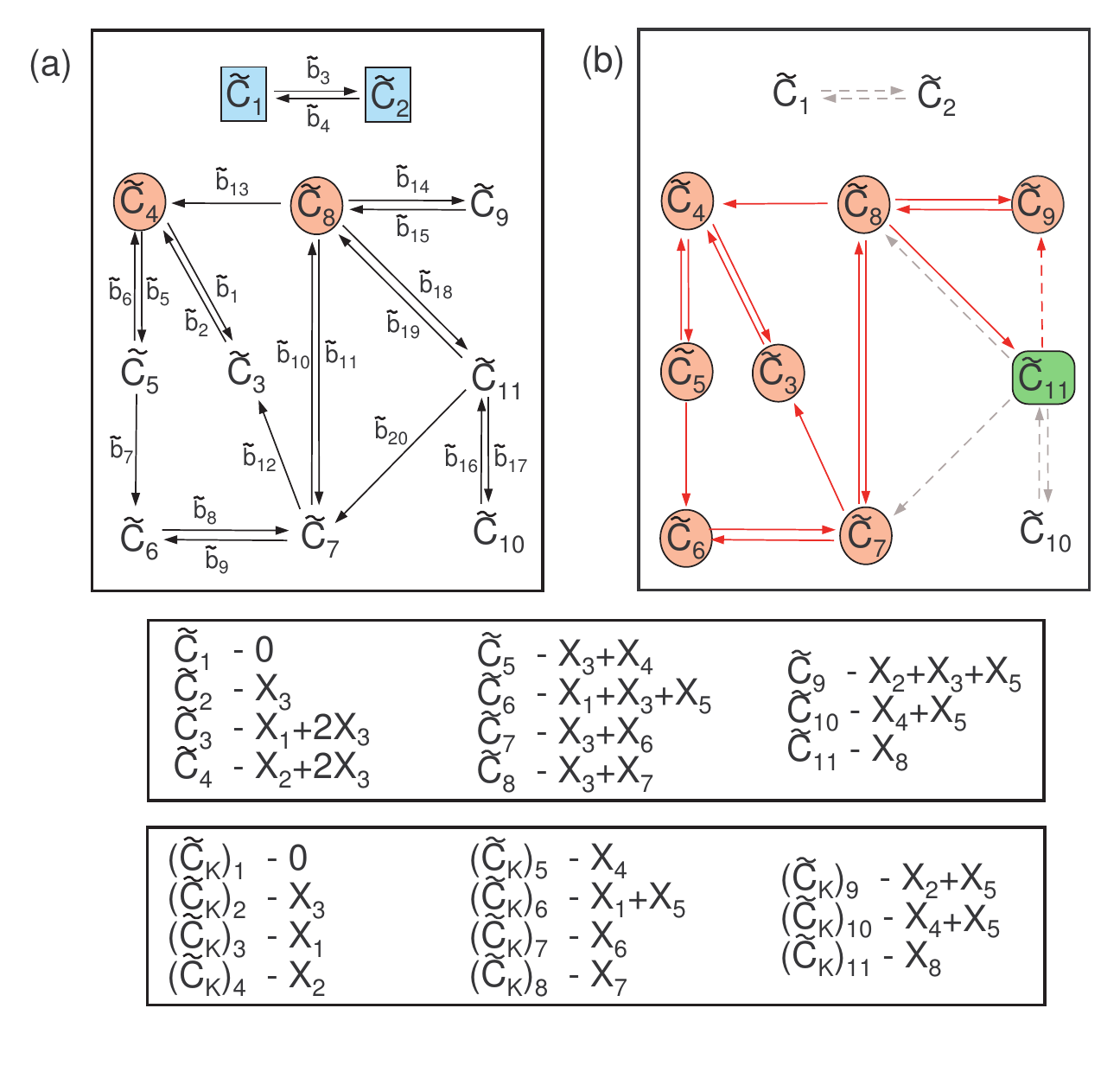}
\caption{Reaction-weighted translation $\tilde{\mathcal{N}}(\tilde{\mathcal{B}})= (\tilde{\mathcal{S}},\tilde{\mathcal{C}},\tilde{\mathcal{C}}_K,\tilde{\mathcal{R}},\tilde{\mathcal{B}})$ corresponding $\mathcal{N}(\mathcal{K})$. In (a), we identify the improper complex set $\tilde{\mathcal{C}}_I$ (pink) and the resolving complex set $\tilde{\mathcal{C}}_R$ (blue). In (b), a computationally-determined network $(\tilde{\mathcal{S}},\tilde{\mathcal{C}}^* \cup \tilde{\mathcal{C}}^{**}, \tilde{\mathcal{R}}^* \cup \tilde{\mathcal{R}}^{**})$ with $\tilde{\mathcal{C}}^*$ (pink), $\tilde{\mathcal{C}}^{**}$ (green), $\tilde{\mathcal{R}}^*$ (solid red arrows), and $\tilde{\mathcal{R}}^{**}$ (dashed red arrows) identified.}
\label{figure4}
\end{figure}

A recurring network structure is the one contained in Figure \ref{figure4}. It can be verified directly that the network structure is valid for all reaction-weights with the reaction-weights given in Table \ref{table2}. The translation $\tilde{\mathcal{N}}(\tilde{\mathcal{B}})$ is improper since both $C_2 = X_2$ and $C_5 = X_2 + X_3$ are translated to $\tilde{C}_4 = X_2 + 2X_3$, and both $C_9 = X_7$ and $C_{13} = X_3 + X_7$ are translated to $\tilde{C}_8 = X_3 + X_7$. With the selection of kinetic complexes given in Figure \ref{figure4}, we have that the source complexes $C_{5} = X_2 + X_3$ and $C_{13} = X_3 + X_7$ are not translated to the new network.

In order to see whether the corresponding monomials may be ``resolved'', we check Definition 3.5. We have that
\[\begin{split} \tilde{\mathcal{C}}_I & = \{ 4, 8\} \\ h^{-1}(4) & = \{ 2, 5\} \\ h^{-1}(8) & = \{ 9, 13 \} \\ \tilde{S}_I & = \mbox{span} \{ y_5 - y_2, y_{13} - y_9 \} \\ & = \mbox{span} \{ [0,0,1,0,0,0,0,0] \}. \end{split}\]
To construct the resolving complex set $\tilde{\mathcal{C}}_R$ according to Definition 3.6, we first verify that $\tilde{S}_I \subseteq \tilde{S}_K$. For this example, we have that
\begin{equation}
\label{432}
\begin{split} y_{5} - y_{2} & = [0,0,1,0,0,0,0,0] = (\tilde{y}_K)_2 - (\tilde{y}_K)_1 \\ y_{13} - y_9 & = [0,0,1,0,0,0,0,0] = (\tilde{y}_K)_2 - (\tilde{y}_K)_1.\end{split}
\end{equation}
It follows that we may take $\tilde{\mathcal{C}}_R = \{ 1, 2\}$. Notice that (\ref{432}) corresponds to the trivial identities
\begin{equation}
\label{23}
\begin{split} & \mathbf{x}^{y_5} = x_2x_3 = \left( \frac{x_3}{1} \right) x_2 = \left( \frac{\mathbf{x}^{(\tilde{y}_K)_2}}{\mathbf{x}^{(\tilde{y}_K)_1}} \right) \mathbf{x}^{(\tilde{y}_K)_4}, \\ & \mathbf{x}^{y_{13}} = x_3x_7  = \left( \frac{x_3}{1} \right) x_7= \left( \frac{\mathbf{x}^{(\tilde{y}_K)_2}}{\mathbf{x}^{(\tilde{y}_K)_1}} \right) \mathbf{x}^{(\tilde{y}_K)_8}.\end{split}
\end{equation}
That is, we may relate the untranslated monomials $x_2x_3$ and $x_3x_7$ to monomials which appear in $\tilde{\mathcal{N}}(\tilde{\mathcal{B}})$ (see Figure \ref{figure4}). 

\begin{table}[t]
\caption{Reaction-weights for correspondence of $\mathcal{N}(\mathcal{B})$ and $\tilde{\mathcal{N}}(\tilde{\mathcal{K}})$ according to Definition 3.2.}
\label{table2}
\begin{center}
\begin{tabular}[h]{|l|l|l|l|}
\hline
$\tilde{b}_1 = k_1$ & $\tilde{b}_{6} = k_6$ & $\tilde{b}_{11} = k_{11}$ & $\tilde{b}_{16} = k_{16}$ \\
$\tilde{b}_2 = k_2$ & $\tilde{b}_7 = k_7$ & $\tilde{b}_{12} = k_{12}$ & $\tilde{b}_{17} = k_{17} + k_{20}$ \\
$\tilde{b}_3 =k_3$ & $\tilde{b}_8 = k_8$ & $\tilde{b}_{13} = k_{13}$ & $\tilde{b}_{18} = k_{18}$ \\
$\tilde{b}_4 = k_4$ & $\tilde{b}_9 = k_9$ & $\tilde{b}_{14} = k_{14}$ & $\tilde{b}_{19} = k_{19}$ \\
$\tilde{b}_5 = k_5 $ & $\tilde{b}_{10} = k_{10}$ & $\tilde{b}_{15} = k_{15}$ & $\tilde{b}_{20} = k_{20}$ \\
\hline
\end{tabular}
\end{center}
\end{table}

For this example, it is trivial to verify that the technical conditions of Lemma 6.1 and Theorem 3.1 are satisfied since the improper complex set $\tilde{\mathcal{C}}_I = \{ 4, 8\}$ is contained entirely within a different linkage class of $\tilde{\mathcal{N}}(\tilde{\mathcal{B}})$ than the resolving complex set $\tilde{\mathcal{C}}_R = \{1, 2\}$ (see Figure \ref{figure4}(a)). Specifically, the technical conditions of Lemma 6.1 are satisfied because there are no paths from complexes in  $\tilde{\mathcal{C}}_I$ to complexes in $\tilde{\mathcal{C}}_R$, and the technical conditions of Theorem 3.1 may be satisfied by choosing $\tilde{\mathcal{C}}^* \cup \tilde{\mathcal{C}}^{**}$ to coincide with the linkage class containing $\tilde{\mathcal{C}}_I$. A computationally-determined alternative way to satisfy the conditions, which requires fewer complexes, is given in Figure \ref{figure4}(b).

We may not, however, apply Lemma 6.1 or Theorem 3.1 directly to conclude steady state resolvability since $\tilde{\delta} = 2$ for $\tilde{\mathcal{N}}(\tilde{\mathcal{B}})$. To guarantee steady-state resolvability, however, it is sufficient to guarantee that, at steady-state, the ratio $\mathbf{x}^{(\tilde{y}_K)_2} / \mathbf{x}^{(\tilde{y}_K)_1}$ is a constant value. Since $\Psi_K(\mathbf{x}) \in$ ker$(\tilde{Y} \cdot A(\tilde{\mathcal{B}}))$ is a necessary and sufficient condition for steady state, we show that $\tilde{C}_1$ and $\tilde{C}_2$ are both contained on the support of exactly one element of ker$(\tilde{Y} \cdot A(\tilde{\mathcal{B}}))$. We have that
\[\begin{split} & {\tiny \tilde{Y} \cdot A(\tilde{\mathcal{B}}) = \left[
\begin{array}{ccccccc}
0 & 0 & -\tilde{b}_2 & \tilde{b}_1 & \tilde{b}_7 & -\tilde{b}_8 & \tilde{b}_{12}+\tilde{b}_9 \\
0 & 0 & \tilde{b}_2 & -\tilde{b}_1 - \tilde{b}_5 & \tilde{b}_6 & 0 & \tilde{b}_{12} \\
\tilde{b}_3 & -\tilde{b}_4 & 0 & -\tilde{b}_5 & \tilde{b}_6 & 0 & \tilde{b}_{13} \\
0 & 0 & 0 & \tilde{b}_5 & -\tilde{b}_6 - \tilde{b}_7 & 0 & 0 \\
0 & 0 & 0 & 0 & \tilde{b}_7 & -\tilde{b}_8 & \tilde{b}_9 \\
0 & 0 & 0 & 0 & 0 & \tilde{b}_8 & -\tilde{b}_9 - \tilde{b}_{10} -\tilde{b}_{12} \\
0 & 0 & 0 & 0 & 0 & 0 & \tilde{b}_{10} \\
0 & 0 & 0 & 0 & 0 & 0 & 0 \\
\end{array} \right.} \\
& \hspace{2in}
{\tiny \left. \begin{array}{cccc}
0 & 0 & 0 & 0 \\ 
\tilde{b}_{13}-\tilde{b}_{14} & -\tilde{b}_{15} & 0 & 0 \\
\tilde{b}_{13}-\tilde{b}_{18} & 0 & 0 & \tilde{b}_{20} + \tilde{b}_{19} \\
0 & 0 & -\tilde{b}_{16} & \tilde{b}_{17} \\
\tilde{b}_{14} & -\tilde{b}_{15} & -\tilde{b}_{16} & \tilde{b}_{17} \\
\tilde{b}_{11} & 0 & 0 & \tilde{b}_{20} \\
-\tilde{b}_{11}-\tilde{b}_{13} -\tilde{b}_{14}-\tilde{b}_{18} & \tilde{b}_{15} & 0 & \tilde{b}_{19} \\
\tilde{b}_{18} & 0 & \tilde{b}_{16} & -\tilde{b}_{17}-\tilde{b}_{19}-\tilde{b}_{20}
\end{array}
\right].} \end{split}\]
It can be easily computed that the dimension of ker$(\tilde{Y} \cdot A(\tilde{\mathcal{B}}))$ is four and that every vector which has support on the first and second component has the form
\[t[\tilde{b}_4, \tilde{b}_3, *, *, *, *, *, *, *, *, *]\]
for some $t \not= 0$. It follows that at every steady state we have
\[\begin{split} \mathbf{x}^{(y_K)_1} & = t \tilde{b}_4 \\ \mathbf{x}^{(y_K)_2} & = t \tilde{b}_3\end{split}\]
so that
\[x_3 = \frac{\mathbf{x}^{(\tilde{y}_K)_2}}{\mathbf{x}^{(\tilde{y}_K)_1}} = \frac{\tilde{b}_3}{\tilde{b}_4}.\]
Consequently, from (\ref{23}) and the reaction-weights given in Table \ref{table2}, we have that at steady state
\[x_2x_3 = \frac{k_3}{k_4} x_2 \; \; \; \mbox{ and } \; \; \; x_3x_7 =  \frac{k_3}{k_4} x_7.\]
It follows that, if we choose $(\tilde{C}_K)_4 = C_2 = X_2$ and $(\tilde{C}_K)_8 = C_9 = X_7$, we may relate the untranslated monomials $x_2x_3$ and $x_3x_7$ to $x_2$ and $x_7$, respectively, with a rescaling of the corresponding reaction weight. The complete list of required rate constants is given in Table \ref{table1}. The generalized mass action system corresponding to $\tilde{\mathcal{N}}(\tilde{\mathcal{K}})$ is steady state equivalent to the  mass action system corresponding to $\mathcal{N}(\mathcal{K})$.

\begin{table}[t]
\caption{Reaction-weights for steady state equivalence of $\mathcal{N}(\mathcal{K})$ and $\tilde{\mathcal{N}}(\tilde{\mathcal{K}})$.}
\label{table1}
\begin{center}
\begin{tabular}[h]{|l|l|l|l|}
\hline
$\tilde{k}_1 = k_1$ & $\tilde{k}_{6} = k_6$ & $\tilde{k}_{11} = k_{11}$ & $\tilde{k}_{16} = k_{16}$ \\
$\tilde{k}_2 = k_2$ & $\tilde{k}_7 = k_7$ & $\tilde{k}_{12} = k_{12}$ & $\tilde{k}_{17} = k_{17} + k_{20}$ \\
$\tilde{k}_3 =k_3$ & $\tilde{k}_8 = k_8$ & $\tilde{k}_{13} = k_{13}$ & $\tilde{k}_{18} = \frac{k_3}{k_4} k_{18}$ \\
$\tilde{k}_4 = k_4$ & $\tilde{k}_9 = k_9$ & $\tilde{k}_{14} = k_{14}$ & $\tilde{k}_{19} = k_{19}$ \\
$\tilde{k}_5 = \frac{k_3}{k_4} k_5 $ & $\tilde{k}_{10} = k_{10}$ & $\tilde{k}_{15} = k_{15}$ & $\tilde{k}_{20} = k_{20}$ \\
\hline
\end{tabular}
\end{center}
\end{table}


\begin{remark}
It is worth noting that, although the mixed-integer linear programming algorithm requires numerical reaction-weights, it can inform subsequent symbolic analysis. In particular, even though the program only determines a reaction-weighted translation for the given numerical parameter values, it was possible for this example to then verify that the translation works for all parameter values.
\end{remark}

\begin{remark}
Despite the successful runs of the algorithm, numerical stability remains an issue. The $25$ successful runs were produced from a sample of $27$ runs. The two unsuccessful runs produced no feasible solution for the linear relaxation. The author suspects this is a result of the interplay between the chosen value of the parameter $\epsilon$ and the ranges of the stochastic parameters $k[i]$ and $v[i,j]$ (which depend on $\epsilon$). Tightening and simplifying the algorithm remains a significant priority moving forward.
\end{remark}

\section{Motivating Example for Future Work}

We now present an example of a network which falls beyond the scope of Theorem 3.1, and the underlying theory in \cite{J1}, but for which an improper reaction-weighted translation exists which is steady-state resolvable to the original network. That is, we show that the resolvability conditions which exist in the literature to date are sufficient but not necessary.

Consider the reaction-weighted chemical reaction network $\mathcal{N}(\mathcal{K})$ given by:
\begin{equation}
\label{1}
\begin{split} & X_1 \hspace{-0.1cm} \mathop{\stackrel{k_1}{\begin{array}{c} \vspace{-0.25cm} \longleftarrow \vspace{-0.3cm} \\ \longrightarrow \end{array}}}_{k_5} \hspace{-0.1cm} X_2 \stackrel{k_2}{\longrightarrow} X_3 \\ & X_1 + X_3 \stackrel{k_3}{\longrightarrow} 2X_1 \\ & X_2 + X_3 \stackrel{k_4}{\longrightarrow} 2X_2.\end{split}
\end{equation}
In order to characterize the steady states, we wish to determine a network translation. It can be quickly verified that the translation scheme
\[\begin{array}{ll}
X_1 \hspace{-0.1cm} \begin{array}{c} \vspace{-0.25cm} \longleftarrow \vspace{-0.3cm} \\ \longrightarrow \end{array} \hspace{-0.1cm} X_2 \longrightarrow X_3 \hspace{0.5cm} \vspace{0.1cm} & (+ \emptyset) \\ X_1 + X_3 \longrightarrow 2X_1 & (-X_1) \vspace{0.1cm} \\ X_2 + X_3 \longrightarrow 2X_2 & (-X_2)
\end{array}\]
yields the following generalized network $\tilde{\mathcal{N}}(\tilde{\mathcal{B}})$ with reaction-weight set $\tilde{\mathcal{B}} = \{ \tilde{b}_i = k_i \; | \; i=1, \ldots, 5 \}$:
\begin{equation}
\label{2}
\begin{array}{c} \displaystyle{X_1 \hspace{-0.1cm} \mathop{\stackrel{\tilde{b}_1}{\begin{array}{c} \vspace{-0.25cm} \longleftarrow \vspace{-0.3cm} \\ \longrightarrow \end{array}}}_{\tilde{b}_5} \hspace{-0.1cm} X_2} \vspace{-0.1cm} \\ \hspace{0.2cm} {}_{\tilde{b}_3} \hspace{-0.1cm} \nwarrow \hspace{0.1cm} {}_{\tilde{b}_4} \hspace{-0.1cm} \nnearrow \hspace{-0.1cm} \sswarrow \hspace{-0.1cm} {}_{\tilde{b}_2} \\ \hspace{0.2cm} X_3 \end{array}
\end{equation}
The network is improper since both $X_1+X_3$ and $X_2 + X_3$ are translated to $X_3$. It can be quickly determined that, $\tilde{S}_I \subseteq \tilde{S}_K$, $\tilde{\delta}=0$, $\tilde{\mathcal{C}}_I = \{ X_3 \}$, and $\tilde{\mathcal{C}}_R = \{ X_1, X_2 \}$, but that the technical conditions of Lemma 6.1 and Corollary 3.1 may not be satisfied. The slightly more general algebraic conditions for resolvability presented in \cite{J1} can also be shown to fail. (This is shown directly by (\ref{9429}).)

We instead consider ``resolving'' the monomials $x_1x_3$ and $x_2x_3$ directly. We choose $X_2+X_3$ to be the kinetic complex corresponding to the complex $X_3$ in (\ref{2}). It is clear that we have
\[x_1x_3 = \left( \frac{x_1}{x_2} \right) x_2 x_3\]
so that the untranslated monomial $x_1x_3$ is related to the translated monomials $x_1$, $x_2$, and $x_2x_3$. To show steady state resolvability, it is sufficient to show that the ratio $x_1/x_2$ is constant at every steady state of the generalized mass action system corresponding to (\ref{2}). To accomplish this, we allow the reaction-weight for every reaction not associated with $x_1x_3$ to correspond to its pre-translation value. That is, we set $\tilde{k}_1 = k_1$, $\tilde{k}_2 = k_2$, $\tilde{k}_4 = k_4$, and $\tilde{k}_5 = k_5$, and allow $\tilde{k}_3$ to remain undetermined. With these values, the steady state equations corresponding to (\ref{1}) and (\ref{2}), respectively, are given by
\begin{equation}
\label{3}
\begin{split}
0 & = -k_1 x_1 + k_5 x_2 + k_3 x_1 x_3 \\
0 & = k_1 x_1 - (k_2 + k_5)x_2 + k_4 x_2 x_3 \\
0 & = k_2 x_2 - k_3 x_1 x_3 - k_4 x_2 x_3
\end{split}
\end{equation}
and
\begin{equation}
\label{4}
\begin{split}
0 & = -k_1 x_1 + k_5 x_2 + \tilde{k}_3 x_2 x_3 \\
0 & = k_1 x_1 - (k_2 + k_5)x_2 + k_4 x_2 x_3 \\
0 & = k_2 x_2 - (\tilde{k}_3 + k_4 )x_2 x_3.
\end{split}
\end{equation}

These systems differ only in the monomials $k_3 x_1 x_3$ and $\tilde{k}_3 x_2 x_3$. We wish to explicitly relate these monomials at steady state. We notice, first of all, that the generalized network (\ref{2}) is weakly reversible and $\tilde{\delta}=0$. It follows that
\[\tilde{Y} \; A(\tilde{\mathcal{K}}) \; \Psi_K(\mathbf{x}) = \mathbf{0} \; \; \; \Longleftrightarrow \; \; \; A(\tilde{\mathcal{K}}) \; \Psi_K(\mathbf{x}) = \mathbf{0} \; \; \; \Longleftrightarrow \; \; \; \Psi_K(\mathbf{x}) \in \mbox{ker}(A(\tilde{\mathcal{K}})).\]
We can compute that
\[\mbox{ker}(A(\tilde{\mathcal{K}})) = \mbox{span} \{ [K_1,K_2,K_3] \}\]
where we have the tree constants $K_1 = k_2 \tilde{k}_3 + \tilde{k}_3 k_5 + k_4 k_5$, $K_2 = k_1 \tilde{k}_3 + k_1 k_4$, and $K_3 = k_1 k_2$ (see Appendix A of the main text). From the structure of $\tilde{\Psi}_K(\mathbf{x}) = (x_1,x_2,x_2x_3),$ it follows that, at steady state, we have
\begin{equation}
\label{9429}
\frac{x_1}{x_2} = \frac{\tilde{k}_3(k_2 + k_5) + k_4 k_5}{k_1 \tilde{k}_3 + k_1 k_4}.
\end{equation}

Returning to the steady state conditions (\ref{3}) and (\ref{4}), we require
\[\tilde{k}_3 x_2 x_3 = k_3 x_1 x_3 = k_3 \left( \frac{x_1}{x_2} \right) x_2 x_3 = k_3 \left( \frac{\tilde{k}_3(k_2 + k_5) + k_4 k_5}{k_1 \tilde{k}_3 + k_1 k_4} \right) x_2 x_3.\]
To determine $\tilde{k}_3$, we solve the quadratic
\[k_1 \tilde{k}_3^2 + (k_1 k_4 - (k_2 + k_5)k_3) \tilde{k}_3 - k_3 k_4 k_5 = 0\]
to get
\begin{equation}
\label{5}
\tilde{k}_3 = \frac{(k_2 + k_5)k_3 - k_1 k_4 \pm \sqrt{((k_2+k_5)k_3-k_1k_4)^2 + 4k_1k_3k_4k_5}}{2k_1}.
\end{equation}
The value under the root is strictly positive and greater in magnitude than the value outside. It follows that, for all rate values, we may pick a positive value for $\tilde{k}_3$ so that the monomials $k_3 x_1 x_3$ and $\tilde{k}_3 x_2 x_3$ coincide at steady state.

It is worth emphasizing where the process used to obtain (\ref{9429}) differs from what is allowed by the process given in \cite{J1} (which provided the basis for Theorem 3.1). In \cite{J1}, we also aspired to resolve ratios of monomials at steady state as ratios of tree constants, as in (\ref{9429}); however, in \cite{J1} the ratios of tree constants were not allowed to depend upon any rate constant corresponding to a reaction attached to an ``unresolved'' monomial. That is, $\tilde{k}_3$ could not appear \emph{implicitly} on the right-hand side of (\ref{9429}). Nevertheless, we have succeeded in determining that there is a value of $\tilde{k}_3$ for which the correspondence between the original and generalized mass action systems can be made.


As a more concrete illustration, we now take the following reaction-weights for the original reaction-weighted network $\mathcal{N}(\mathcal{K})$:
\begin{equation}
\label{6}
\begin{split} k_1 = 1, \; k_2 = \frac{3}{4}, \; k_3 = 2, \; k_4 = 2, \; k_5 = \frac{1}{4}.\end{split}
\end{equation}
This choice simplifies (\ref{5}) since $k_1k_3k_4k_5 =1$ and $(k_2+k_5)k_3-k_1k_4 = 0$ so that $\tilde{k}_3 = \pm 1$. That is, we may consider the reaction-weighted translation $\tilde{\mathcal{N}}(\tilde{\mathcal{K}})$ with the reaction-weights:
\begin{equation}
\label{7}
\tilde{k}_1 = 1, \; \tilde{k}_2 = \frac{3}{4}, \; \tilde{k}_3 = 1, \; \tilde{k}_4 = 2, \; \tilde{k}_5 = \frac{1}{4}.
\end{equation}
We now return to the steady state equations (\ref{3}) and (\ref{4}) for $\mathcal{N}(\mathcal{K})$ and $\tilde{\mathcal{N}}(\tilde{\mathcal{K}})$, respectively. It can be computed that the Gr\"{o}bner basis for (\ref{3}) with reaction-weights (\ref{6}) is given by
\[I = \langle x_2 (4x_3 -1 )(4x_3 - 3),x_1 + 2x_2 x_3 - x_2 \rangle\]
while the Gr\"{o}bner basis of (\ref{4}) with reaction-weights (\ref{7}) is given by
\[I = \langle x_2(4x_3-1),2x_1-x_2 \rangle.\]

We can see that there are three steady state possibilities for (\ref{3}):
\[\begin{split} \mbox{(i)} & \; \; x_1 = 0, x_2 = 0, x_3 \mbox{ arbitrary} \\ \mbox{(ii)} & \; \; x_2 = 2x_1, x_3 = \frac{1}{4} \\ \mbox{(iii)} & \; \; x_2 = -2x_1, x_3 = \frac{3}{4}.\end{split}\]
The steady state manifold for (iii) is not chemically interesting, since one of $x_1$ or $x_2$ must be negative unless $x_1=x_2$ (which intersects case (i)). The steady state manifold of (\ref{7}) therefore captures the two chemical interesting steady state conditions of (\ref{3}), specifically, cases (i) and (ii). That is, as expected, the steady states of the system (\ref{3}) with reaction-weights (\ref{6}) correspond to the steady states of the system (\ref{4}) with reaction-weights (\ref{7}) on the nonnegative orthant $\mathbb{R}_{\geq 0}^3$.

This example, while preliminary and lacking full symbolic analysis, demonstrates that the conditions given in \cite{J1} and refined in this paper are \emph{sufficient} but \emph{not necessary} for ``resolving'' untranslated monomials. Further exploration of when implicit equations of the form (\ref{5}) may be obtained and solved to yield a chemically-meaningful reaction-weights for a translated chemical reaction network will be explored further in future work.

\begin{remark}
It can be checked that the non-chemical choice of reaction-weights
\begin{equation}
\label{8}
\tilde{k}_1 = 1, \; \tilde{k}_2 = \frac{3}{4}, \; \tilde{k}_3 = -1, \; \tilde{k}_4 = 2, \; \tilde{k}_5 = \frac{1}{4}.
\end{equation}
yields a polynomial system (\ref{4}) which has the Gr\"{o}bner basis
\[I = \langle 4x_2x_3 - 3x_2, 2x_1 + x_2 \rangle.\]
This is easily recognizable as corresponding to the steady state conditions (i) and (iii) for $\mathcal{N}(\mathcal{K})$. That is to say, we are able to recover all of the steady state information from the original mass action system by considering \emph{two} different reaction-weighted chemical reaction networks---one chemically feasible, and one not.
\end{remark}
